\documentclass[b5paper,openany]{thesis}
\advance\textwidth1cm
\usepackage{amsmath,amssymb,euscript,latexsym,fancyheadings,theorem,calc}
\renewcommand{\chaptermark}[1]{\markboth{\chaptername\ \thechapter.\ \ #1}{}}

\renewcommand\theta{\vartheta}
\renewcommand\phi{\varphi}
\newcommand\smod[1]{\operatorname{mod} #1}
\newcommand\Z{\mathbf{Z}}
\newcommand\R{\mathbf{R}}
\newcommand\C{\mathbf{C}}
\newcommand\Q{\mathbf{Q}}
\newcommand\Ha{\EuScript{H}}
\newcommand\ef{\EuScript{F}}
\newcommand\sign{\operatorname{sgn}}
\newcommand\im{\operatorname{Im}}

\newcommand\qed{\hfill $\Box$\smallskip}
\newcommand\proof[1][]{\par\noindent{\bf Proof#1:}\ \ignorespaces}
\newcommand\Res{\operatornamewithlimits{Res}}
\newcommand\spa{\operatorname{span}}

\newtheorem{proposition}{Proposition}[chapter]
\newtheorem{corollary}[proposition]{Corollary}
\newtheorem{lemma}[proposition]{Lemma}
\newtheorem{theorem}[proposition]{Theorem}
{\theorembodyfont{\rmfamily}
\newtheorem{definition}[proposition]{Definition}
\newtheorem{remark}[proposition]{Remark} 
\newtheorem{example}[proposition]{Example}
}
\begin{document}
\frontmatter
\thispagestyle{empty}
\begin{center}
{\huge Mock Theta Functions}\\
\vspace{1cm}
{\Large Mock Thetafuncties}\\
\vspace{5mm}
(met een samenvatting in het Nederlands)\\
\vfill
\textsc{Proefschrift}\\
\vspace{3cm}
\parbox[c]{9cm}{\textsc{Ter verkrijging van de graad van doctor aan de Universiteit Utrecht op gezag van de Rector Magnificus, Prof.\ dr.\ W.H. Gispen, ingevolge het besluit van het College voor Promoties in het openbaar te verdedigen op woensdag 30 oktober 2002 des ochtends te 10.30 uur}}\\
\vspace{1cm}
\textsc{door}\\
\vspace{3cm}
\textsc{Sander Pieter Zwegers}\\
\vspace{5mm}
\textsc{geboren op 16 april 1975, te Oosterhout}
\end{center}

\thispagestyle{empty}
\noindent\begin{tabular}{ll}
Promotor: & Prof.\ dr.\ D.B. Zagier\\
& Max-Planck-Institut f\"ur Mathematik, Bonn\\
& oud-hoogleraar aan de Universiteit Utrecht\\
Copromotor: & Dr.\ R.W. Bruggeman\\
& Faculteit Wiskunde en Informatica\\
& Universiteit Utrecht
\end{tabular}
\vfill
\noindent\rule[1mm]{\textwidth}{.1mm}
\rule[1mm]{\textwidth}{.1mm}
Zwegers, Sander Pieter\\
Mock Theta Functions\\
Proefschrift Universiteit Utrecht --- Met samenvatting in het Nederlands\\
\rule[1mm]{\textwidth}{.1mm}
ISBN 90-393-3155-3\\
\rule[1mm]{\textwidth}{.1mm}
Printed by PrintPartners Ipskamp, Enschede

\tableofcontents
\clearemptydoublepage
\mainmatter
\chapter*{Introduction}
\addcontentsline{toc}{chapter}{\numberline{}Introduction}
\markboth{Introduction}{}
\section*{Mock $\theta$-functions}
\addcontentsline{toc}{section}{\numberline{}Mock $\theta$-functions}
\markright{Mock $\theta$-functions}

Early in 1920, three months before his death, S. Ramanujan wrote his last letter to G.H. Hardy. For the mathematical part of this letter see \cite[pp. 127--131]{lost} (also reproduced in \cite{andrews1}). In the course of it he said: ``I discovered very interesting functions recently which I call `Mock' $\theta$-functions. Unlike the `False' $\theta$-functions (studied partially by Prof.\ Rogers in his interesting paper \cite{rogers}) they enter into mathematics as beautifully as the ordinary $\theta$-functions. I am sending you with this letter some examples.'' He then provided a long list of mock $\theta$-functions, together with identities satisfied by them. The first three pages in which Ramanujan explained what he meant by a mock $\theta$-function are very obscure. 

G.N. Watson wrote the first papers (\cite{watson1} and \cite{watson2}) to elucidate the mock $\theta$-functions. The first of these is Watson's Presidential Address to the London Mathematical Society in 1935. He entitled it ``The Final Problem: An Account of the Mock Theta Functions.'' In it he writes: ``I make no apologies for my subject being what is now regarded as old-fashioned, because, as a friend remarked to me a few months ago, I am an old-fashioned mathematician.'' His methods may have been a bit old-fashioned, but looking at the number of articles on mock $\theta$-functions that have appeared since 1935, or even in the last ten years, we must conclude that the subject is still up-to-date. 

In these two papers, Watson proves most of the assertions found in the letter of Ramanujan. The first paper considers only the third-order functions. It provides three new mock $\theta$-functions not mentioned in the letter. The bulk of the paper is devoted to the modular transformation properties of these functions. 
To get these transformations, he first proves certain identities. For example, for the third order mock $\theta$-function 
\begin{equation*}
f(q):= 1+\frac{q}{(1+q)^2}+\frac{q^4}{(1+q)^2 (1+q^2)^2} + \frac{q^9}{(1+q)^2 (1+q^2)^2 (1+q^3)^2}  + \cdots
\end{equation*}
he finds:
\begin{equation}\label{efje}
f(q)= \frac{2}{(q)_\infty} \sum_{n\in\Z} \frac{ (-1)^n q^{\frac{3}{2}n^2+\frac{1}{2}n}}{1+q^n},
\end{equation}
with $q=e^{2\pi i\tau}$, $\tau\in\Ha:=\{\tau\in\C\mid \im(\tau)>0\}$, and $(q)_\infty=\prod_{n=1}^\infty (1-q^n) =q^{-\frac{1}{24}} \eta(\tau)$, where $\eta$ is the Dedekind eta-function.

In Watson's second paper on mock $\theta$-functions, he moves on to the fifth order functions. He manages to prove all of the identities given by Ramanujan in his letter. However, he is unable to find the modular transformation properties, simply because he is unable to find identities like \eqref{efje} for the fifth order functions. He even expressed his doubts about finding anything comparable to \eqref{efje}. However Andrews (see \cite{andrews3}) was able to find comparable results for most of the fifth order functions. For example, for the fifth order function which Watson denotes by $f_0$ one finds
\begin{equation*}
f_0(q) = \frac{1}{(q)_\infty} \sum_{n\geq0} \sum_{|j|\leq n} (-1)^j q^{\frac{5}{2}n^2+\frac{1}{2}n -j^2}(1-q^{4n+2}),
\end{equation*} 
which may be rewritten as
\begin{equation}\label{efnul}
f_0(q) = \frac{1}{(q)_\infty} \left( \sum_{n+j\geq 0, n-j \geq 0} - \sum_{n+j<0, n-j<0} \right) (-1)^j q^{\frac{5}{2}n^2+\frac{1}{2}n -j^2}.
\end{equation}

The seventh order functions were mostly neglected by Watson, perhaps because Ramanujan makes no positive assertions about them. However A. Selberg (see \cite{selberg}) provides a full account of the behaviour of the seventh order functions near the unit circle. In \cite[pp. 666]{hick2} we find identities similar to \eqref{efnul} for the seventh order mock $\theta$-functions $\ef_0$, $\ef_1$ and $\ef_2$. For example (slightly rewritten)
\begin{equation}\label{kringelefnul}
\ef_0(q) = \frac{1}{(q)_\infty}\left( \sum_{r,s\geq0} -\sum_{r,s<0}\right) (-1)^{r+s} q^{\frac{3}{2} r^2 +4rs+\frac{3}{2}s^2+\frac{1}{2}r+\frac{1}{2}s}.
\end{equation}

In \cite{gottsche} L. G\"ottsche and D. Zagier consider sums like the ones in \eqref{efnul} and \eqref{kringelefnul}. They call them theta-functions for indefinite lattices. For some special cases they find modular transformation properties for these functions. However, these results do not include the sums in \eqref{efnul} and \eqref{kringelefnul}. In \cite{polishchuk}, a theorem about the modularity of a certain family of $q$-series associated with indefinite binary quadratic forms is given. Again, the results do not include the sums in \eqref{efnul} and \eqref{kringelefnul}.

\section*{This thesis}
\addcontentsline{toc}{section}{\numberline{}This thesis}
\markright{This thesis}

This thesis is the result of my research on the following two questions, both posed by Don Zagier:
\begin{enumerate}
\item How do the mock $\theta$-functions fit in the theory of modular forms?
\item Is there a theory of indefinite theta functions?
\end{enumerate} 
Since most of the mock $\theta$-functions had been related to sums like the one in \eqref{efje}, I first considered this type of sum. The result of this research is Chapter 1. In it we consider the series
\begin{equation*}
\sum_{n\in\Z} \frac{(-1)^n e^{\pi i(n^2+n)\tau+2\pi inv}}{1-e^{2\pi in\tau + 2\pi iu}}\qquad (\tau\in\Ha, v\in\C, u\in\C\setminus (\Z\tau+\Z)).
\end{equation*}
This function was also studied by Lerch in \cite{lerch1} (see \cite{lerch2} for an abstract). Therefore we call this a Lerch sum. This sum is of the same type as the sum in \eqref{efje}. The function does not transform like a Jacobi form. However, we find that on addition of a (relatively easy) correction term the function does transform like a Jacobi form. This correction term is real-analytic.

In Chapter 2 we consider certain indefinite $\theta$-functions, in an attempt to give a partial answer to the second question. These indefinite $\theta$-functions are modified versions of the sums considered by G\"ottsche and Zagier in \cite{gottsche}. We find elliptic and modular transformation properties for these functions. Because of the modifications the indefinite $\theta$-functions are no longer holomorphic (in general). Although the results in this chapter are more general than the results in \cite{gottsche}, the second question is far from being solved. This is because we only consider indefinite quadratic forms of type $(r-1,1)$. It remains a problem of considerable interest to develop a theory of theta-series for quadratic forms of arbitrary type.

In \cite{coeff} Andrews gives most of the fifth order mock theta functions as Fourier coefficients of meromorphic Jacobi forms, namely certain quotients of ordinary Jacobi theta-series. This is the motivation for the study of the modularity of Fourier coefficients of meromorphic Jacobi forms, in Chapter 3. We find that modularity follows on adding a real-analytic correction term to the Fourier coefficients.

In Chapter 4 we use the results from Chapter 2, together with \eqref{efnul}, \eqref{kringelefnul} and similar identities for other mock $\theta$-functions, to get the modular transformation properties of the seventh-order mock $\theta$-functions and most of the fifth-order functions. The final result is that we can write each of these mock $\theta$-functions as the sum of two functions $H$ and $G$, where:
\begin{itemize}
\item $H$ is a real-analytic modular form of weight 1/2 and is an eigenfunction of the appropriate Casimir operator with eigenvalue 3/16 (this is also the eigenvalue of holomorphic modular forms of this weight; for the theory of real-analytic modular forms see for example \cite[Ch.\ IV]{maass}); and
\item $G$ is a theta series associated to a negative definite unary quadratic form, i.e.\ has the form $\sum \sign(f) \beta(2f^2 y) e^{-\pi if^2\tau-2\pi ifb}$, where $f$ ranges over a certain arithmetic progression $a\Z+b$ ($a,b\in\Q$), $\tau=x+iy\in\Ha$ and $\beta(x)= \int_{x}^{\infty} u^{-\frac{1}{2}} e^{-\pi u}du$. Moreover $G$ is bounded as $\tau$ tends vertically to any rational limit.
\end{itemize}

This decomposition is thus similar to the one found in \cite{zagier} for an Eisenstein series of weight 3/2, the holomorphic part of the series, in that case, having class numbers as Fourier coefficients. 

Many of the results of Chapter 4 could also be deduced using the methods from Chapter 1 or Chapter 3 instead of Chapter 2, i.e.\ we have actually given 3 approaches to proving modularity properties of the mock $\theta$-functions.

\clearemptydoublepage
\chapter{Lerch Sums}
\section{Introduction}
In this chapter, we first study a function $h$, which is essentially the function $\phi$ studied by Mordell in \cite{mordell1} and \cite{mordell2}. We reproduce some of the results of Mordell. In the next section we study a function $\mu$, which is essentially a function studied by Lerch in \cite{lerch1}. We find elliptic and modular transformation properties of this function. One of these transformations involves the function $h$:
\begin{equation*}
\frac{1}{\sqrt{-i\tau}}\ e^{\pi i(u-v)^2/\tau} \mu\Bigl(\frac{u}{\tau},\frac{v}{\tau};-\frac{1}{\tau}\Bigr)+ \mu(u,v;\tau) =\frac{1}{2i} h(u-v;\tau).
\end{equation*}
In section 4 we find a real-analytic function $R$ with essentially the same elliptic and modular transformation properties as $\mu$. Combining the properties of $\mu$ and $R$ we find a real-analytic function $\tilde{\mu}$, which transforms like a Jacobi form.

In the last section we relate $h$ to a period integral of a unary theta function of weight 3/2.

\section{The Mordell integral}
 
In this section, we present results of Mordell found in \cite{mordell1} and \cite{mordell2}, in a form suitable for the purpose of this chapter. The function $h$ defined in Definition \ref{def1} is essentially the function $\phi$ studied by Mordell: $\phi(x;\tau)=-\frac{1}{2}\tau e^{-\pi i\tau/4 +\pi ix} h(x-\frac{\tau}{2} +\frac{1}{2};\tau)$. The same function was used earlier by Riemann (as described by Siegel \cite{siegel2}) to prove the functional equation for the Riemann zeta function. Mordell was the first to analyze the behaviour of this integral relative to modular transformations. Consequently, we shall refer to this integral as the Mordell integral.

\begin{definition}\label{def1}
For $z\in\C$ and $\tau\in\Ha$ set
\begin{equation*}
h(z)=h(z;\tau):= \int_{\R} \frac{e^{\pi i\tau x^2 -2\pi zx}}{\cosh\pi x}\,dx.
\end{equation*}
\end{definition}

\begin{proposition} \label{prop1}
The function $h$ has the following properties:
\begin{description}
\item[(1)] $h(z)+h(z+1)=\frac{2}{\sqrt{-i\tau}}\, e^{\pi i(z+\frac{1}{2})^2/\tau}$,
\item[(2)] $h(z)+e^{-2\pi iz-\pi i\tau} h(z+\tau) =2e^{-\pi iz-\pi i\tau/4}$,
\item[(3)] $z\mapsto h(z;\tau)$ is the unique holomorphic function satisfying (1) and (2),
\item[(4)] $h$ is an even function of $z$,
\item[(5)] $h(\frac{z}{\tau};-\frac{1}{\tau})=\sqrt{-i\tau}\, e^{-\pi iz^2/
\tau}\, h(z;\tau)$,
\item[(6)] $\displaystyle{h(z;\tau) = e^{\frac{\pi i}{4}} h(z;\tau +1) +e^{-\frac{\pi i}{4}} \frac{e^{\pi iz^2/(\tau+1)}}{\sqrt{\tau+1}}\ h\Bigl(\frac{z}{\tau+1};\frac{\tau}{\tau+1}\Bigr)}$.
\end{description}
\end{proposition}
\proof
(1) We have
\begin{equation*}
h(z)+h(z+1)=\int_{\R} \frac{e^{\pi i\tau x^2 -2\pi zx}}{\cosh\pi x}(1+e^{-2\pi x})\,dx= 2 \int_{\R} e^{\pi i\tau x^2 -2\pi x(z+\frac{1}{2})}\,dx.
\end{equation*}
This last integral is well known and equals $\frac{1}{\sqrt{-i\tau}} e^{\pi i (z+\frac{1}{2})^2/\tau}$.\\
(2) If we change $x$ into $x+i$ we find
\begin{equation*}
\int_{\R+i} \frac{e^{\pi i\tau x^2 -2\pi zx}}{\cosh\pi x}\,dx= \int_{\R} \frac{e^{\pi i\tau (x+i)^2 -2\pi z(x+i)}}{\cosh\pi (x+i)}\,dx= -e^{-2\pi iz-\pi i\tau} h(z+\tau).
\end{equation*}
Now using Cauchy's theorem we find
\begin{equation*}
\begin{split}
h(z)+e^{-2\pi iz-\pi i\tau} h(z+\tau)&=\left(\int_{\R}-\int_{\R+i}\right) \frac{e^{\pi i\tau x^2 -2\pi zx}}{\cosh\pi x}\,dx\\
 &= 2\pi i \Res_{x=i/2} \frac{e^{\pi i\tau x^2 -2\pi zx}}{\cosh\pi x}= 2e^{-\pi iz-\pi i\tau/4}.
\end{split}
\end{equation*}
(3) If we have two holomorphic functions $h_1$ and $h_2$ both satisfying the two equations, their difference $f=h_1-h_2$ is a holomorphic function satisfying:
\begin{equation*}
\begin{cases} f(z)+f(z+1)=0&\\
f(z) +e^{-2\pi iz -\pi i\tau} f(z+\tau)=0.&
\end{cases}
\end{equation*}
So $f(z_0+m\tau+n)=(-1)^{m+n} e^{\pi im^2\tau +2\pi imz_0} f(z_0)$. Letting $z_0$ vary over a fundamental parallelogram $[0,1)\Z+[0,1)$ and $m,n$ vary over $\Z$, we see that $f(z)$ is bounded and tends to 0 as $\im(z)\rightarrow \infty$, so $f\equiv 0$ by Liouville's theorem.\\
(4) Replace $x$ by $-x$ in the integral.\\
(5) Let $g(x)=\frac{1}{\cosh \pi x}$. We first compute the Fourier transform $\ef g$ of $g$: Using Cauchy's formula we get
\begin{equation*}
\left( \int_{\R} -\int_{\R+i} \right) \frac{e^{2\pi izx}}{\cosh \pi x}\, dx =2\pi i \operatornamewithlimits{Res}_{x=i/2} \frac{e^{2\pi izx}}{\cosh \pi x} =2e^{-\pi z},
\end{equation*}
but 
\begin{equation*}
\int_{\R+i} \frac{e^{2\pi izx}}{\cosh \pi x}\, dx =\int_{\R} \frac{e^{2\pi iz(x+i)}}{\cosh \pi (x+i)}\, dx
=-e^{-2\pi z} \int_{\R} \frac{e^{2\pi izx}}{\cosh \pi x}\, dx,
\end{equation*}
so we find 
\begin{equation*}
(\ef g)(z) := \int_{\R} \frac{e^{2\pi izx}}{\cosh \pi x}\, dx
= \frac{2e^{-\pi z}}{1+e^{-2\pi z}} = g(z).
\end{equation*}
We see that $g$ is its own Fourier transform! (Note the unusual plus sign in the definition of the Fourier transform).

Let $f_{\tau} (x) = e^{\pi i \tau x^2}$, $\tau\in\Ha$. The Fourier transform of $f_{\tau}$ is given by
\begin{equation*}
\ef f_{\tau} = \frac{1}{\sqrt{-i\tau}} f_{-\frac{1}{\tau}}.
\end{equation*}
We now see
\begin{equation*}
\begin{split}
\int_{\R} \frac{e^{\pi i\tau x^2 +2\pi izx}}{\cosh \pi x}\, dx &= \ef (f_{\tau} \cdot g)(z) = (\ef f_{\tau} ) \ast (\ef g)(z) \\
&= \frac{1}{\sqrt{-i\tau}} f_{-\frac{1}{\tau}} \ast g(z)
= \frac{1}{\sqrt{-i\tau}} \int_{\R} \frac{e^{\pi i\frac{-1}{\tau} (z-x)^2}}{\cosh \pi x}\, dx.
\end{split}
\end{equation*}
This identity holds for $z\in\R$. Since both sides are analytic functions of $z$, the identity holds for all $z\in\C$. If we replace $z$ by $iz$ we get the desired result. 

We may also prove the identity of part (5) by using (1) and (2) to show that $z\mapsto \frac{1}{\sqrt{-i\tau}} e^{\pi iz^2/\tau} h(\frac{z}{\tau};-\frac{1}{\tau})$ also satisfies the two equations (1) and (2). By uniqueness we get the equation.\\
(6) Using (1) and (2) we can show that the right hand side, considered as a function of $z$, also satisfies (1) and (2). The equation now follows from (3).
\qed

\section{Lerch sums}
In this section we will study the function
\begin{equation*}
\sum_{n\in\Z} \frac{(-1)^n e^{\pi i(n^2+n)\tau+2\pi inv}}{1-e^{2\pi in\tau + 2\pi iu}}\qquad (\tau\in\Ha, v\in\C, u\in\C\setminus (\Z\tau+\Z)).
\end{equation*}
This function was also studied by Lerch. The original paper \cite{lerch1} is in Czech and is not very easy to obtain. See \cite{lerch2} for an abstract in German. We will prove elliptic and modular transformation properties of this function in Proposition \ref{prop3} and \ref{prop4} respectively. These results are equivalent to the results found by Lerch. 

It is more convenient to normalize the above sum by dividing by the classical Jacobi theta function $\theta$. (Lerch did this too.) We will first give, without proof, some standard properties of $\theta$. For the theory of $\theta$-functions see \cite{mumford}. 
\begin{proposition} \label{prop2}
For $z\in\C$ and $\tau\in\Ha$ define
\begin{equation*}
\theta(z)=\theta(z;\tau):=\sum_{\nu\in\frac{1}{2}+\Z} e^{\pi i\nu^2 \tau +2\pi i\nu(z+\frac{1}{2})}.
\end{equation*}
Then $\theta$ satisfies:
\begin{description}
\item[(1)] $\theta(z+1)=-\theta(z)$.
\item[(2)] $\theta(z+\tau)= -e^{-\pi i\tau -2\pi iz} \theta(z)$.
\item[(3)] Up to a multiplicative constant, $z\mapsto \theta(z)$ is the unique holomorphic function satisfying  (1) and (2).
\item[(4)] $\theta(-z)=-\theta(z)$.
\item[(5)] The zeros of $\theta$ are the points $z=n\tau+m$, with $n,m\in\Z$. These are simple zeros.
\item[(6)] $\theta(z;\tau+1)=e^{\frac{\pi i}{4}} \theta(z;\tau)$.
\item[(7)] $\theta(\frac{z}{\tau};-\frac{1}{\tau})=-i\sqrt{-i\tau} e^{\pi i z^2/\tau}\theta(z;\tau)$.
\item[(8)] $\displaystyle{\theta(z;\tau)=-iq^{\frac{1}{8}} \zeta^{-\frac{1}{2}} \prod_{n=1}^{\infty} (1-q^n)(1-\zeta q^{n-1})(1-\zeta^{-1} q^n)}$, with $q=e^{2\pi i\tau}$, $\zeta=e^{2\pi iz}$. This is the Jacobi triple product identity.
\item[(9)] $\theta'(0;\tau+1)= e^{\frac{\pi i}{4}} \theta'(0;\tau)$ and $\theta'(0;-\frac{1}{\tau})=(-i\tau)^{3/2}\ \theta'(0;\tau)$.
\item[(10)] $\theta'(0;\tau)=-2\pi\eta(\tau)^3$, with $\eta$ as in the introduction.
\end{description}
\end{proposition}

We now turn to the normalized version of Lerch's function:

\begin{proposition} \label{prop3}
For $u,v\in\C\setminus(\Z\tau+\Z)$ and $\tau\in\Ha$, define
\begin{equation*}
\mu(u,v)=\mu(u,v;\tau):= \frac{e^{\pi iu}}{\theta(v;\tau)} \sum_{n\in\Z} \frac{(-1)^n e^{\pi i(n^2+n)\tau+2\pi inv}}{1-e^{2\pi in\tau + 2\pi iu}}. 
\end{equation*}
Then $\mu$ satisfies:
\begin{description}
\item[(1)] $\mu(u+1,v)=-\mu(u,v),$
\item[(2)] $\mu(u,v+1)=-\mu(u,v),$
\item[(3)] $\mu(u,v)+e^{-2\pi i(u-v) -\pi i\tau} \mu(u+\tau,v)=-ie^{-\pi i(u-v)-\pi i\tau/4},$
\item[(4)] $\mu(u+\tau,v+\tau)=\mu(u,v),$
\item[(5)] $\mu(-u,-v)=\mu(u,v),$
\item[(6)] $u\mapsto \mu(u,v)$ is a meromorphic function, with simple poles in the points $u=n\tau+m$ $(n,m\in\Z)$, and residue $\frac{-1}{2\pi i} \frac{1}{\theta(v)}$ in $u=0$,
\item[(7)] $\displaystyle{\mu(u+z,v+z)-\mu(u,v)=  \frac{1}{2\pi i} \frac{\theta'(0)\theta(u+v+z) \theta(z)}{\theta(u) \theta(v) \theta(u+z) \theta(v+z)}}$,\\
for $u,v,u+z,v+z\not\in\Z\tau +\Z$,
\item[(8)] $\mu(v,u)=\mu(u,v).$
\end{description}
\end{proposition}
\proof
(1) is trivial and (2) follows from (1) of Proposition \ref{prop2}.\\
(3) The definition of $\theta$ gives the following:
\begin{equation*}
\begin{split}
i &e^{-\pi i\tau/4+\pi iv} \theta(v) = \sum_{n\in\Z} (-1)^n e^{\pi i(n^2-n)\tau +2\pi inv} \\
&= \sum_{n\in\Z} \frac{(-1)^n e^{\pi i(n^2-n)\tau +2\pi inv}}{1- e^{2\pi in\tau+2\pi iu}} \left(1- e^{2\pi in\tau+2\pi iu}\right)\\
&= -e^{2\pi iv} \sum_{n\in\Z} \frac{(-1)^n e^{\pi i(n^2+n)\tau +2\pi inv}}{1- e^{2\pi in\tau+2\pi i(u+\tau)}}-e^{2\pi iu} \sum_{n\in\Z} \frac{(-1)^n e^{\pi i(n^2+n)\tau +2\pi inv}}{1- e^{2\pi in\tau+2\pi iu}}.
\end{split}
\end{equation*}
Dividing both sides by $-e^{\pi iu} \theta(v)$, we get the desired result.\\
(4) Part (2) of Proposition \ref{prop2} gives
\begin{equation*}
\begin{split}
\mu(u+\tau,v+\tau) &= \frac{e^{\pi i(u+\tau)}}{\theta(v+\tau)} \sum_{n\in\Z} \frac{(-1)^n e^{\pi i(n^2+n)\tau+2\pi in(v+\tau)}}{1-e^{2\pi in\tau + 2\pi i(u+\tau)}}\\
&= -\frac{e^{\pi i(u+\tau) +\pi i\tau +2\pi iv}}{\theta(v)} \sum_{n\in\Z} \frac{(-1)^n e^{\pi i(n^2+3n)\tau+2\pi inv}}{1-e^{2\pi i(n+1)\tau + 2\pi iu}}.
\end{split}
\end{equation*}
Replace $n$ by $n-1$ in the last sum to get the desired result.\\
(5) If we replace $n$ by $-n$ in the definition of $\mu$ we see
\begin{equation*}
\mu(u,v) = \frac{e^{\pi iu}}{\theta(v)} \sum_{n\in\Z} \frac{(-1)^n e^{\pi i(n^2-n)\tau-2\pi inv}}{1-e^{-2\pi in\tau + 2\pi iu}}.
\end{equation*}
We multiply by $\frac{-e^{2\pi in\tau -2\pi iu}}{-e^{2\pi in\tau -2\pi iu}}$ to find
\begin{equation*}
\mu(u,v) = -\frac{e^{-\pi iu}}{\theta(v)} \sum_{n\in\Z} \frac{(-1)^n e^{\pi i(n^2+n)\tau-2\pi inv}}{1-e^{2\pi in\tau - 2\pi iu}}.
\end{equation*}
Now using (4) of Proposition \ref{prop2} we find
\begin{equation*}
\mu(u,v)=\mu(-u,-v).
\end{equation*}
(6) From the definition we see that $u\mapsto \mu(u,v)$ has a simple pole if $1-e^{2\pi in\tau+2\pi iu}=0$, for some $n\in\Z$. So $u\mapsto \mu(u,v)$ has simple poles in the points $u=-n\tau+m$ $(n,m\in\Z)$.

The pole in $u=0$ comes from the term $n=0$. We see
\begin{equation*}
\lim_{u\rightarrow 0} u\,\mu(u,v) = \frac{1}{\theta(v)} \lim_{u\rightarrow 0} \frac{u}{1-e^{2\pi iu}} =\frac{-1}{2\pi i} \frac{1}{\theta(v)}.
\end{equation*}
(7) Consider $f(z)=\theta(u+z)\theta(v+z)\left( \mu(u+z,v+z)-\mu(u,v)\right)$. Using (1), (2) and (5) of Proposition \ref{prop2}, and (1), (2), (4) and (6) of this proposition, we see that $f$ has no poles, a zero for $z=0$, and satisfies
\begin{equation*}
\begin{cases}
f(z+1)=f(z)&\\
f(z+\tau)= e^{-2\pi i\tau-2\pi i(u+v+2z)} f(z).&
\end{cases}
\end{equation*}
It follows that the quotient $f(z)/\theta(z) \theta(u+v+z)$ is a double periodic function with at most one simple pole in each fundamental parallelogram, and hence constant:
\begin{equation}
f(z)= C(u,v) \theta(z) \theta(u+v+z).
\label{form}
\end{equation}
To compute $C$ we consider $z=-u$. If we take $z=-u$ in \eqref{form} we find
\begin{equation}\label{f1}
f(-u)=C(u,v)\theta(-u)\theta(v) = -C(u,v)\theta(u)\theta(v)
\end{equation}
by (4) of Proposition \ref{prop2}.

By definition we have 
\begin{equation}\label{f2}
\begin{split}
f(-u)&= \lim_{z\rightarrow -u} \theta(u+z)\theta(v+z)\left(\mu(u+z,v+z)-\mu(u,v)\right)\\
&= \theta(v-u) \cdot \lim_{z\rightarrow 0} \theta(z)\mu(z,v-u)\\
&= \theta(v-u) \cdot \lim_{z\rightarrow 0} \frac{\theta(z)}{z} \cdot \lim_{z\rightarrow 0} z \mu(z,v-u)= -\frac{1}{2\pi i} \theta'(0),
\end{split}
\end{equation}
where we have used (6).

Combining \eqref{f1} and \eqref{f2} gives the desired result.\\
(8) Take $z=-u-v$ in (7) and use (5) of Proposition \ref{prop2} to find
\begin{equation*}
\mu(-v,-u)=\mu(u,v).
\end{equation*}
If we now use (5), we get the desired result. \qed
 
\begin{proposition} \label{prop4}
Let $\mu$ be as in Proposition \ref{prop3}. Then $\mu$ satisfies the following  modular transformation properties:
\begin{description}
\item[(1)] $\displaystyle{\mu(u,v;\tau+1)=e^{-\frac{\pi i}{4}} \mu(u,v;\tau),}$
\item[(2)] $\displaystyle{\frac{1}{\sqrt{-i\tau}}\ e^{\pi i(u-v)^2/\tau} \mu\left(\frac{u}{\tau},\frac{v}{\tau};-\frac{1}{\tau}\right)+ \mu(u,v;\tau) =\frac{1}{2i} h(u-v;\tau)}$,\\
with $h$ as in Definition \ref{def1}. 
\end{description}
\end{proposition}
\proof
(1) Use (6) of Proposition \ref{prop2}.\\
(2) Replacing $(u,v,z,\tau)$ by $(\frac{u}{\tau},\frac{v}{\tau},\frac{z}{\tau},-\frac{1}{\tau})$ in (7) of Proposition \ref{prop3} and using (7) and (9) of Proposition \ref{prop2} we see that the left hand side depends only on $u-v$, not on $u$ and $v$ separately. Call it $\frac{1}{2i}\tilde{h}(u-v;\tau)$. Using (1) and (3) of Proposition \ref{prop3} we see that $\tilde{h}$ satisfies the two identities (1) and (2) of Proposition \ref{prop1}, so if we can prove that $\tilde{h}$ is a holomorphic function, then we may conclude that $\tilde{h}=h$, as desired.

The poles of both $u\mapsto \mu(u,v)$ and $u\mapsto \mu(\frac{u}{\tau},\frac{v}{\tau};-\frac{1}{\tau})$ are simple, and occur at $u\in\Z\tau+\Z$, so the only poles of $u\mapsto \tilde{h}(u-v)$ could be simple poles for $u\in\Z\tau+\Z$. Since this is a function of $u-v$ it has no poles at all, and hence is holomorphic. 

Alternatively, we can check, using (6) of Proposition \ref{prop3} and (7) of Proposition \ref{prop2}, that the residue at $u=0$ vanishes. By (1) and (3) of Proposition \ref{prop3} the residues vanish for all $u\in\Z\tau +\Z$, hence $\tilde{h}$ is holomorphic. \qed

\section{A real-analytic Jacobi form?}

\begin{definition}\label{defE} For $z\in\C$ we define
\begin{equation*}
E(z)=2\int_0^z e^{-\pi u^2} du=\sum_{n=0}^{\infty} \frac{(-\pi)^n}{n!} \frac{z^{2n+1}}{n+1/2}.
\end{equation*}
\end{definition}
This is an odd entire function of $z$.
\begin{lemma} \label{lem131} For $z\in\R$ we have
\begin{equation*}
E(z)= \sign(z)\left(1-\beta(z^2)\right),
\end{equation*}
where
\begin{equation*}
\beta(x) = \int_x^\infty u^{-\frac{1}{2}} e^{-\pi u} du \qquad (x\in\R_{\geq 0}).
\end{equation*}
\end{lemma}
\proof
Write $\int_0^z e^{-\pi u^2} du$ as $\sign(z) \int_0^{|z|} e^{-\pi u^2} du$ and substitute $u=\sqrt{v}$. \qed

We consider for $u\in\C$ and $\tau\in\Ha$ the series
\begin{equation*}
R(u;\tau)=\sum_{\nu\in\frac{1}{2}+\Z} \Bigl\{ \sign(\nu) - E\Bigl((\nu+a)\sqrt{2y}\Bigr)\Bigr\}\ (-1)^{\nu-\frac{1}{2}} e^{-\pi i \nu^2\tau -2\pi i\nu u},
\end{equation*}
with $y=\im(\tau)$ and $a=\frac{\im(u)}{\im(\tau)}$.
\begin{lemma}\label{defR}   
For all $c,\epsilon>0$, this series converges absolutely and uniformly on the set $\{ u\in\C,\tau\in\Ha \mid |a|<c, y>\epsilon\}$. The function $R$ it defines is real-analytic and satisfies
\begin{equation}\label{dif1}
\frac{\partial R}{\partial \overline{u}}(u;\tau)= \sqrt{2} y^{-1/2} e^{-2\pi a^2 y} \theta(\overline{u};-\overline{\tau})
\end{equation}
and 
\begin{equation}\label{dif2}
\frac{\partial}{\partial \overline{\tau}}R(a\tau-b;\tau)=-\frac{i}{\sqrt{2y}} e^{-2\pi a^2y} \sum_{\nu\in\frac{1}{2}+\Z} (-1)^{\nu-\frac{1}{2}} (\nu+a) e^{-\pi i \nu^2\overline{\tau} -2\pi i\nu (a\overline{\tau}-b)}.
\end{equation}
\end{lemma}
\proof
We split $\sign(\nu) - E((\nu+a)\sqrt{2y})$ into the sum of $\sign(\nu)-\sign((\nu+a)\sqrt{2y})$ and $\sign((\nu+a)\sqrt{2y})\beta(2(\nu+a)^2 y)$. We see that $\sign(\nu)-\sign((\nu+a)\sqrt{2y})$ is nonzero for only a finite number of values $\nu\in\frac{1}{2}+\Z$ (this number depends on $a$, but since $a$ is bounded, so is this number). Hence the series
\begin{equation*}
\sum_{\nu\in\frac{1}{2}+\Z} \Bigl\{ \sign(\nu) - \sign\Bigl((\nu+a)\sqrt{2y}\Bigr)\Bigr\} (-1)^{\nu-\frac{1}{2}} e^{-\pi i \nu^2\tau -2\pi i\nu u}
\end{equation*}
converges absolutely and uniformly. 

We can easily see that $0\leq \beta(x)\leq e^{-\pi x}$ for all $x\in\R_{\geq 0}$, hence
\begin{equation*}
\begin{split}
&\left|\left\{ \sign\Bigl((\nu+a)\sqrt{2y}\Bigr)\beta\Bigl(2(\nu+a)^2 y\Bigr)\right\} (-1)^{\nu-\frac{1}{2}} e^{-\pi i \nu^2\tau -2\pi i\nu u}\right| \\
&\leq e^{-2\pi (\nu+a)^2 y} \left| e^{-\pi i \nu^2\tau -2\pi i\nu u}\right| \\ 
&= e^{-\pi (\nu+a)^2 y -\pi a^2 y} \leq e^{-\pi (\nu+a)^2 \epsilon}.
\end{split}
\end{equation*}
We have the inequality
\begin{equation*}
(\nu +a)^2 \geq \frac{1}{2} \nu^2,
\end{equation*}
for $|\nu|\geq\nu_0$, for some $\nu_0\in\R$ which depends only on $c$ ($a$ is bounded by $c$). Hence we see that the series
\begin{equation*}
\sum_{\nu\in\frac{1}{2}+\Z} \left\{ \sign\Bigl((\nu+a)\sqrt{2y}\Bigr)\beta\Bigl(2(\nu+a)^2 y\Bigr)\right\} (-1)^{\nu-\frac{1}{2}} e^{-\pi i \nu^2\tau -2\pi i\nu u}
\end{equation*}
converges absolutely and uniformly on the given set. 

Since $R$ is the (infinite) sum of real-analytic functions, and the series converges absolutely and uniformly, it is real-analytic.

We fix $\tau\in\Ha$, and determine $u=a\tau-b$ by the coordinates $a,b\in\R$. We see
{\allowdisplaybreaks
\begin{align*}
&\left(\frac{\partial}{\partial a}+ \tau\frac{\partial}{\partial b} \right) R(a\tau-b;\tau)\\
&= \left(\frac{\partial}{\partial a}+ \tau\frac{\partial}{\partial b} \right) \sum_{\nu\in\frac{1}{2}+\Z} \left\{ \sign(\nu) - E\Bigl((\nu+a)\sqrt{2y}\Bigr)\right\} (-1)^{\nu-\frac{1}{2}} e^{-\pi i \nu^2\tau -2\pi i\nu (a\tau-b)}\\
&= -\sqrt{2y} \sum_{\nu\in\frac{1}{2}+\Z} E'\Bigl((\nu+a)\sqrt{2y}\Bigr) (-1)^{\nu-\frac{1}{2}} e^{-\pi i \nu^2\tau -2\pi i\nu (a\tau-b)}\\
&= -2\sqrt{2y} \sum_{\nu\in\frac{1}{2}+\Z} e^{-2\pi(\nu+a)^2 y} (-1)^{\nu-\frac{1}{2}} e^{-\pi i \nu^2\tau -2\pi i\nu (a\tau-b)}\\
&= -2\sqrt{2y} e^{-2\pi a^2y} \sum_{\nu\in\frac{1}{2}+\Z} (-1)^{\nu-\frac{1}{2}} e^{-\pi i \nu^2\overline{\tau} -2\pi i\nu (a\overline{\tau}-b)}\\
&=-2i\sqrt{2y} e^{-2\pi a^2 y} \theta(a\overline{\tau} -b;-\overline{\tau}),
\end{align*}
with $\theta$ as in Proposition \ref{prop2} and the term-by-term differentiation being easily justified. Since $\frac{\partial}{\partial \overline{u}} =\frac{i}{2y} \left(\frac{\partial}{\partial a}+ \tau\frac{\partial}{\partial b} \right)$, this gives the differential equation \eqref{dif1}. Similarly
\begin{equation*}
\begin{split}
&\frac{\partial}{\partial \overline{\tau}}R(a\tau-b;\tau)\\
&= \frac{1}{2}\left(\frac{\partial}{\partial x}+ i\frac{\partial}{\partial y} \right) \sum_{\nu\in\frac{1}{2}+\Z} \left\{ \sign(\nu) - E\Bigl((\nu+a)\sqrt{2y}\Bigr)\right\} (-1)^{\nu-\frac{1}{2}} e^{-\pi i \nu^2\tau -2\pi i\nu (a\tau-b)}\\
&= -\frac{i}{2}\frac{1}{\sqrt{2y}} \sum_{\nu\in\frac{1}{2}+\Z} (\nu +a)\ E'\Bigl((\nu+a)\sqrt{2y}\Bigr)\ (-1)^{\nu-\frac{1}{2}} e^{-\pi i \nu^2\tau -2\pi i\nu (a\tau-b)}\\
&= -\frac{i}{\sqrt{2y}} e^{-2\pi a^2y} \sum_{\nu\in\frac{1}{2}+\Z} (-1)^{\nu-\frac{1}{2}} (\nu+a)\ e^{-\pi i \nu^2\overline{\tau} -2\pi i\nu (a\overline{\tau}-b)},
\end{split}
\end{equation*}
proving equation \eqref{dif2}. \qed}
\begin{proposition} \label{prop5}
The function $R$ has the following elliptic transformation properties:
\begin{description}
\item[(1)] $R(u+1)=-R(u)$,
\item[(2)] $R(u)+e^{-2\pi iu-\pi i\tau} R(u+\tau) = 2e^{-\pi iu-\pi i\tau/4}$,
\item[(3)] $R(-u)=R(u)$.
\end{description}
\end{proposition}
\proof
Part (1) is trivial, and for (3) we replace $\nu$ by $-\nu$ in the sum and use the fact that $E$ is an odd function. To prove (2), we start with
\begin{equation*}
\begin{split}
&e^{-2\pi iu-\pi i\tau} R(u+\tau) \\
&=e^{-2\pi iu-\pi i\tau} \sum_{\nu\in\frac{1}{2}+\Z} \left\{ \sign(\nu) - E\Bigl((\nu+a+1)\sqrt{2y}\Bigr)\right\} (-1)^{\nu-\frac{1}{2}} e^{-\pi i \nu^2\tau -2\pi i\nu (u+\tau)}\\
&= -\sum_{\nu\in\frac{1}{2}+\Z} \left\{ \sign(\nu-1) - E\Bigl((\nu+a)\sqrt{2y}\Bigr)\right\} (-1)^{\nu-\frac{1}{2}} e^{-\pi i \nu^2\tau -2\pi i\nu u},
\end{split}
\end{equation*}
where we have replaced $\nu$ by $\nu -1$. We now find
\begin{equation*}
\begin{split}
R(u)+&e^{-2\pi iu-\pi i\tau} R(u+\tau) \\
&=\sum_{\nu\in\frac{1}{2}+\Z} \Bigl\{ \sign(\nu) - \sign(\nu-1) \Bigl\} (-1)^{\nu-\frac{1}{2}} e^{-\pi i \nu^2\tau -2\pi i\nu u} = 2e^{-\pi iu-\pi i\tau/4},
\end{split}
\end{equation*}
since $\sign(\nu)-\sign(\nu-1)$ is zero for all $\nu\in \frac{1}{2}+\Z$ except for $\nu=\frac{1}{2}$. \qed

\begin{proposition} \label{prop6}
$R$ has the following  modular transformation properties:
\begin{description}
\item[(1)] $\displaystyle{R(u;\tau+1)=e^{-\frac{\pi i}{4}} R(u;\tau),}$
\item[(2)] $\displaystyle{\frac{1}{\sqrt{-i\tau}}\ e^{\pi iu^2/\tau} R\Bigl(\frac{u}{\tau};-\frac{1}{\tau}\Bigr)\ + R(u;\tau) = h(u;\tau).}$
\end{description}
\end{proposition}
\proof
Part (1) is trivial. The left hand side of (2) we call $\tilde{h}(u;\tau)$. Using (1) and (2) of Proposition \ref{prop5} we can see that $\tilde{h}$ satisfies:
\begin{equation*}
\begin{cases}
\tilde{h}(u)+\tilde{h}(u+1)=\frac{2}{\sqrt{-i\tau}}\, e^{\pi i(u+\frac{1}{2})^2/\tau},&\\
\tilde{h}(u)+e^{-2\pi iu-\pi i\tau} \tilde{h}(u+\tau) =2e^{-\pi iu-\pi i\tau/4}.&
\end{cases}
\end{equation*}
Part (3) of Proposition \ref{prop1} determines $h$ as the unique holomorphic function with these properties. This reduces the proof to showing that $\tilde{h}$ is a holomorphic function of $u$.

We fix $\tau\in\Ha$, and determine $u=a\tau-b$ by the coordinates $a,b\in\R$ (this implies $a=\frac{\im(u)}{\im(\tau)}$ as in Lemma \ref{defR}).
Since $\frac{\partial}{\partial \overline{u}} =\frac{i}{2y} \left(\frac{\partial}{\partial a}+ \tau\frac{\partial}{\partial b} \right)$, we have to show that 
\begin{equation*} 
\left(\frac{\partial}{\partial a}+ \tau\frac{\partial}{\partial b} \right) \tilde{h} (a\tau-b;\tau) =0.
\end{equation*}
According to Lemma \ref{defR} we have
\begin{equation}\label{form2}
\left(\frac{\partial}{\partial a}+ \tau\frac{\partial}{\partial b} \right) R(a\tau-b;\tau)= -2i\sqrt{2y} e^{-2\pi a^2 y} \theta(a\overline{\tau} -b;-\overline{\tau})
\end{equation}
We have
\begin{equation*}
\left(\frac{\partial}{\partial a}+ \tau\frac{\partial}{\partial b} \right) R\left(\frac{a\tau-b}{\tau};-\frac{1}{\tau}\right)= \tau \left(\frac{\partial}{\partial b}+\frac{1}{\tau} \frac{\partial}{\partial a} \right) R\left(a-\frac{b}{\tau};-\frac{1}{\tau}\right).
\end{equation*}
Up to a factor $\tau$ this is the same as $\left(\frac{\partial}{\partial a}+ \tau\frac{\partial}{\partial b} \right) R(a\tau-b;\tau)$, with $(a,b,\tau)$ replaced by $(b,-a,-\frac{1}{\tau})$. Hence by \eqref{form2} we find
\begin{equation*}
\begin{split}
\left(\frac{\partial}{\partial a}+ \tau\frac{\partial}{\partial b} \right) R\left(\frac{a\tau-b}{\tau};-\frac{1}{\tau}\right)&=-2i\tau\sqrt{2y'}e^{-2\pi b^2 y'} \theta\left(-\frac{b}{\overline{\tau}}+a;\frac{1}{\overline{\tau}}\right)\\
&= 2i\tau\sqrt{2y'}e^{-2\pi b^2 y'} \theta\left(-\frac{a\overline{\tau}-b}{\overline{\tau}};\frac{1}{\overline{\tau}}\right),
\end{split}
\end{equation*}
with $y'=\im(-\frac{1}{\tau}) = \frac{y}{\tau \overline{\tau}}$. In the last step we have used (4) of Proposition \ref{prop2}.\\
If we now use (7) of Proposition \ref{prop2}, with $z=a\overline{\tau}-b$ and $\tau$ replaced by $-\overline{\tau}$, we see that this equals
\begin{equation}
\begin{split}\label{form3}
2i&\tau\sqrt{2y'} e^{-2\pi b^2 y'} \cdot -i\sqrt{i\overline{\tau}} e^{-\pi i(a\overline{\tau}-b)^2/\overline{\tau}} \theta(a\overline{\tau}-b;-\overline{\tau})\\
&= 2i \sqrt{2y} \sqrt{-i\tau} e^{-\pi i(a\tau-b)^2/\tau} e^{-2\pi a^2 y} \theta(a\overline{\tau}-b;-\overline{\tau}).
\end{split}
\end{equation}
Using \eqref{form2} and \eqref{form3} we find
\begin{equation*}
\begin{split}
\left(\frac{\partial}{\partial a}+ \tau\frac{\partial}{\partial b} \right) &\tilde{h}(a\tau-b;\tau)\\
&= \frac{1}{\sqrt{-i\tau}} e^{\pi i (a\tau-b)^2/\tau} \left(\frac{\partial}{\partial a}+ \tau\frac{\partial}{\partial b} \right) R\left(\frac{a\tau-b}{\tau};-\frac{1}{\tau}\right)\\ &\quad+\left(\frac{\partial}{\partial a}+ \tau\frac{\partial}{\partial b} \right) R(a\tau-b;\tau)=0.
\end{split}
\end{equation*}
We have established the fact that $\tilde{h}$ is holomorphic, and hence equals $h$. \qed

In the next theorem we combine the properties of $\mu$ and $R$ to find a function $\tilde{\mu}$ which is no longer meromorphic, but has better elliptic and modular transformation properties than $\mu$.

\begin{theorem}\label{them1}
We set 
\begin{equation}\label{comb}
\tilde{\mu} (u,v;\tau) = \mu (u,v;\tau) +\frac{i}{2} R(u-v;\tau),
\end{equation}
then 
\begin{description}
\item[(1)] $\displaystyle{\tilde{\mu} (u+k\tau+l,v+m\tau+n)= (-1)^{k+l+m+n} e^{\pi i(k-m)^2\tau +2\pi i(k-m)(u-v)} \tilde{\mu} (u,v)}$,\\
for $k,l,m,n\in\Z$,
\item[(2)] $\displaystyle{\tilde{\mu} \left( \frac{u}{c\tau+d}, \frac{v}{c\tau+d}; \frac{a\tau +b}{c\tau+d}\right) = v(\gamma)^{-3} (c\tau+d)^{\frac{1}{2}} e^{-\pi i c(u-v)^2/(c\tau+d)} \tilde{\mu}(u,v;\tau)}$,\\
for $\gamma=\left(\begin{smallmatrix}a&b\\c&d\end{smallmatrix}\right)\in \operatorname{SL}_2(\Z)$, with $v(\gamma) = \eta(\frac{a\tau+b}{c\tau+d})/\Bigl((c\tau+d)^{\frac{1}{2}} \eta(\tau)\Bigr)$
\item[(3)] $\tilde{\mu}(-u,-v)=\tilde{\mu}(v,u)=\tilde{\mu}(u,v)$,
\item[(4)] $\displaystyle{\tilde{\mu}(u+z,v+z)-\tilde{\mu}(u,v)=  \frac{1}{2\pi i} \frac{\theta'(0)\theta(u+v+z) \theta(z)}{\theta(u) \theta(v) \theta(u+z) \theta(v+z)}}$,\\
for $u,v,u+z,v+z\not\in\Z\tau +\Z$,
\item[(5)] $u\mapsto \tilde{\mu}(u,v)$ has singularities in the points $u=n\tau+m$ $(n,m\in\Z)$. Furthermore we have $\lim_{u\rightarrow 0} u\tilde{\mu}(u,v) = \frac{-1}{2\pi i} \frac{1}{\theta(v)}$.
\end{description}
\end{theorem}
\begin{remark} Parts (1) and (2) of the theorem say that the function $\tilde{\mu}$ transforms like a two-variable Jacobi form of weight $\frac{1}{2}$ and index $\left(\begin{smallmatrix}-1&1\\1&-1\end{smallmatrix}\right)$ (for the theory of Jacobi forms, see \cite{eichler}, where, however, only Jacobi forms of one variable are considered). Furthermore we can find several differential equations satisfied by $\tilde{\mu}$. Therefore we would like to call this function a real-analytic Jacobi form. However, in the literature I haven't been able to find a satisfying definition of a real-analytic Jacobi form. I intend to return to this problem in the future. 
\end{remark}
\begin{remark}
All three function in \eqref{comb} have a property that the other two do not have: $\tilde{\mu}$ transforms well (like a Jacobi form), $\mu$ is meromorphic and $u,v\mapsto R(u-v)$ depends only on $u-v$.
\end{remark}
\proof
(1) Using the first four parts of Proposition \ref{prop3} and the first two of Proposition \ref{prop5} we find
\begin{equation*}
\begin{split}
\tilde{\mu} (u+1,v) &=-\tilde{\mu}(u,v),\\
\tilde{\mu} (u,v+1) &=-\tilde{\mu}(u,v),\\
\tilde{\mu} (u+\tau,v) &=-e^{2\pi i(u-v) +\pi i\tau} \tilde{\mu}(u,v),\\
\tilde{\mu} (u,v+\tau) &=-e^{2\pi i(v-u) +\pi i\tau} \tilde{\mu}(u,v).
\end{split}
\end{equation*}
Combining these equations we get the desired result.\\
(2) Using Proposition \ref{prop4} and Proposition \ref{prop6} we find
\begin{equation*}
\begin{split}
\tilde{\mu} (u,v;\tau+1) &= e^{-\frac{\pi i}{4}} \tilde{\mu}(u,v;\tau)\\
\tilde{\mu} \Bigl( \frac{u}{\tau},\frac{v}{\tau}; -\frac{1}{\tau} \Bigr) &= -\sqrt{-i\tau} e^{-\pi i (u-v)^2/\tau} \tilde{\mu} (u,v;\tau)
\end{split}
\end{equation*}
Set $m(u,v;\tau):=\theta(u-v;\tau) \tilde{\mu}(u,v;\tau)$. Using (6) and (7) of Proposition \ref{prop2} we see 
\begin{equation*}
\begin{split}
m(u,v;\tau+1)&= m(u,v;\tau)\\
m\left(\frac{u}{\tau},\frac{v}{\tau};-\frac{1}{\tau} \right) &= \tau\ m(u,v;\tau)
\end{split}
\end{equation*}
and so 
\begin{equation*}
m\left(\frac{u}{c\tau+d},\frac{v}{c\tau+d};\frac{a\tau+b}{c\tau+d}\right)=(c\tau+d)\ m(u,v;\tau),
\end{equation*}
for all $\left(\begin{smallmatrix}a&b\\c&d\end{smallmatrix}\right)\in\operatorname{SL}_2(\Z)$. Hence
\begin{equation}\label{muc}
\tilde{\mu} \left(\frac{u}{c\tau+d},\frac{v}{c\tau+d};\frac{a\tau+b}{c\tau+d}\right) = (c\tau+d) \frac{\theta(u-v;\tau)}{\theta\left( \frac{u-v}{c\tau+d};\frac{a\tau+b}{c\tau+d} \right)} \tilde{\mu}(u,v;\tau).
\end{equation}
From (6) and (7) of Proposition \ref{prop2} we find
\begin{equation}\label{tach}
\theta \left(\frac{z}{c\tau+d};\frac{a\tau+b}{c\tau+d}\right)= \chi (\gamma) \sqrt{c\tau+d}\ e^{\pi icz^2/(c\tau+d)}\theta(z;\tau),
\end{equation}
with $\chi(\gamma)$ some eighth root of unity. Applying $\left.\frac{d}{dz}\right|_{z=0}$ on both sides gives
\begin{equation*}
\theta'\left(0;\frac{a\tau+b}{c\tau+d}\right)= \chi(\gamma) (c\tau+d)^{\frac{3}{2}} \theta'(0;\tau).
\end{equation*}
Using (10) of Proposition \ref{prop2} we find
\begin{equation*}
\chi(\gamma)=v(\gamma)^3
\end{equation*}
If we combine this with \eqref{muc} and \eqref{tach} we get the desired result.\\
(3) Using (5) of Proposition \ref{prop3} and (3) of Proposition \ref{prop5} we find
\begin{equation*}
\tilde{\mu} (-u,-v) =\tilde{\mu} (u,v)
\end{equation*}
Using (8) of Proposition \ref{prop3} and (3) of Proposition \ref{prop5} we find
\begin{equation*}
\tilde{\mu} (v,u) = \tilde{\mu}(u,v)
\end{equation*}
(4) This follows directly from (7) of Proposition \ref{prop3}.\\
(5) R has no singularities, so the singularities come from $\mu$. The location and nature of these singularities is already given in (6) of Proposition \ref{prop3}. \qed

\section{Period integrals of weight 3/2 unary theta functions}

In this section we will rewrite $h$ in terms of the period integral of a unary theta function of weight 3/2.

To state the main result we need the following definition: 
\begin{definition}\label{defg}
Let $a,b\in\R$ and $\tau\in\Ha$ then
\begin{equation*}
g_{a,b}(\tau):= \sum_{\nu\in a+\Z} \nu e^{\pi i \nu^2\tau +2\pi i\nu b}.
\end{equation*}
\end{definition}
The function $g_{a,b}$ is a unary theta function. 

\begin{proposition} \label{prop141} $g_{a,b}$ satisfies:
\begin{description}
\item[(1)] $g_{a+1,b}(\tau)=g_{a,b}(\tau)$
\item[(2)] $g_{a,b+1}(\tau)=e^{2\pi ia} g_{a,b}(\tau)$
\item[(3)] $g_{-a,-b}(\tau)=-g_{a,b}(\tau)$
\item[(4)] $g_{a,b}(\tau+1)=e^{-\pi ia(a+1)} g_{a,a+b+\frac{1}{2}}(\tau)$
\item[(5)] $g_{a,b}(-\frac{1}{\tau})=ie^{2\pi iab} (-i\tau)^{\frac{3}{2}} g_{b,-a}(\tau)$
\end{description}
\end{proposition}
I will not prove these relations, since they are all easy. For (5) use Poisson summation.\\
From these relations it follows that if $a$ and $b$ are rational, the function $g_{a,b}$ is a modular form of weight 3/2.

\begin{theorem} \label{period} Let $\tau\in\Ha$, then
\begin{description}
\item[(1)] for $a\in(-\frac{1}{2},\frac{1}{2})$ and $b\in\R$ 
\begin{equation*}
\int_{-\overline{\tau}}^{i\infty} \frac{g_{a+\frac{1}{2},b+\frac{1}{2}}(z)}{\sqrt{-i(z+\tau)}}\, dz= -e^{-\pi ia^2 \tau +2\pi ia(b+\frac{1}{2})} R(a\tau-b),
\end{equation*}
with $R$ as in Lemma \ref{defR}.
\item[(2)] for $a,b\in(-\frac{1}{2},\frac{1}{2})$ 
\begin{equation*}
\int_0^{i\infty} \frac{g_{a+\frac{1}{2},b+\frac{1}{2}}(z)}{\sqrt{-i(z+\tau)}}\, dz = -e^{-\pi ia^2\tau +2\pi ia(b+\frac{1}{2})}h(a\tau-b), 
\end{equation*}
with $h$ as in Definition \ref{def1}.
\end{description}
\end{theorem}
\begin{remark} The left hand side in (2) is 1-periodic as a function of $a$ and, up to a factor, as a function of $b$, while the right hand side is not.
\end{remark}
For the proof we need the following two lemmas:

\begin{lemma} \label{lem1}
If $r\in\R$, $r\not=0$, and $\tau\in\Ha$ then
\begin{equation*}
\int_{-\infty}^{\infty} \frac{e^{\pi i\tau w^2}}{w+ir}\, dw = -\pi r \int_0^{i\infty} \frac{e^{\pi ir^2 z}}{\sqrt{-i(z+\tau)}}\,dz.
\end{equation*}
\end{lemma}

\proof
Both sides define a holomorphic function of $\tau\in\Ha$. So we only have to prove the identity for $\tau=it$, with $t\in\R_{>0}$. The identity then becomes
\begin{equation*}
\int_{-\infty}^{\infty} \frac{e^{-\pi tw^2}}{w+ir}\, dw = -\pi ir\int_0^\infty \frac{e^{-\pi r^2 u}}{\sqrt{u+t}}\,du,
\end{equation*}
where we have substituted $z=iu$ in the integral on the right. We can easily see that both sides, considered as functions of $t$, are solutions of $(-\frac{\partial}{\partial t} +\pi r^2) f(t) = -\frac{\pi ir}{\sqrt{t}}$:
\begin{equation*}
\begin{split}
\Bigl(-\frac{\partial}{\partial t} +\pi r^2\Bigr) \operatorname{LHS} &= \pi \int_{-\infty}^{\infty} (w-ir) e^{-\pi tw^2} dw = -\pi ir \int_{-\infty}^{\infty} e^{-\pi tw^2}dw = -\frac{\pi ir}{\sqrt{t}},\\
\left(-\frac{\partial}{\partial t} +\pi r^2\right) \operatorname{RHS}&= -e^{\pi r^2 t} \frac{\partial}{\partial t} (e^{-\pi r^2 t} \operatorname{RHS}) = \pi ir e^{\pi r^2 t} \frac{\partial}{\partial t} \int_t^\infty \frac{e^{-\pi r^2 x}}{\sqrt{x}}\, dx = -\frac{\pi ir}{\sqrt{t}}.
\end{split}
\end{equation*}
Both sides have the same limit 0 as $t\rightarrow\infty$, hence they are equal. \qed

\begin{lemma} \label{lem2}
Let $b\in(-\frac{1}{2},\frac{1}{2})$ and let $z\in\C$, such that $z\not\in\left(\frac{1}{2} +\Z\right)i$, then
\begin{equation*}
- \frac{e^{2\pi bz}}{\cosh \pi z} = \frac{1}{\pi} \sum_{\nu\in \frac{1}{2} +\Z}\frac{e^{2\pi i\nu(b+\frac{1}{2})}} {z-i\nu}.
\end{equation*}
\end{lemma}
\proof
We first show that the series on the right hand side converges (it doesn't converge absolutely). Since Dirichlet's test for convergence (see \cite[pp.\ 17]{whit}) is not directly applicable, we prove this by means of partial summation (we could also do this by comparing the series with the series $\sum_{n\in\Z_{\not= 0}} \frac{e^{2\pi i n(b+\frac{1}{2})}}{n})$: Define
\begin{equation*}
T_{\nu}  := \frac{e^{2\pi i(\nu+1)(b+\frac{1}{2})}}{e^{2\pi i(b+\frac{1}{2})}-1},
\end{equation*}
then $T_\nu - T_{\nu -1} = e^{2\pi i\nu (b+\frac{1}{2})}$ and $|T_\nu|= \frac{1}{|e^{2\pi i(b+\frac{1}{2})}-1|}$. Hence 
\begin{equation*}
\begin{split}
&\underset{\nu_0 \leq \nu \leq \nu_1}{\sum_{\nu\in \frac{1}{2} +\Z}}  \frac{e^{2\pi i\nu(b+\frac{1}{2})}} {z-i\nu}= \underset{\nu_0 \leq \nu \leq \nu_1}{\sum_{\nu\in \frac{1}{2} +\Z}} \frac{T_\nu - T_{\nu-1}} {z-i\nu}\\
&=-i\underset{\nu_0 \leq \nu \leq \nu_1}{\sum_{\nu\in \frac{1}{2} +\Z}} \frac{T_\nu} {(z-i\nu)(z-i(\nu+1))} -\frac{T_{\nu_0 -1}}{z-i\nu_0} +\frac{T_{\nu_1}}{z-i(\nu_1 +1)},
\end{split}
\end{equation*}
with $\nu_0,\nu_1\in\frac{1}{2}+\Z$. We have $\lim_{\nu_0 \rightarrow -\infty}\frac{T_{\nu_0 -1}}{z-i\nu_0} =0$, because $T_{\nu_0}$ is bounded, and also $\lim_{\nu_1 \rightarrow\infty} \frac{T_{\nu_1}}{z-i(\nu_1 +1)}=0$. Since $\sum_{\nu\in \frac{1}{2} +\Z} \frac{T_\nu} {(z-i\nu)(z-i(\nu+1))}$ converges (it converges absolutely), so does $\sum_{\nu\in \frac{1}{2} +\Z}\frac{e^{2\pi i\nu(b+\frac{1}{2})}} {z-i\nu}$.

For $z\not\in i\Z$ we consider the 1-periodic function given by $b\mapsto e^{2\pi bz}$ for $b\in(-\frac{1}{2},\frac{1}{2})$. An easy calculation shows that
\begin{equation*}
\int_{-\frac{1}{2}}^{\frac{1}{2}} e^{2\pi bz} e^{-2\pi inb} db= \frac{\sinh \pi z}{\pi} \frac{(-1)^n}{z-in}
\end{equation*}
The given function is continuous on the interval $(-\frac{1}{2},\frac{1}{2})$, hence we get from the theory of Fourier series that for $b\in(-\frac{1}{2},\frac{1}{2})$
\begin{equation*}
e^{2\pi bz} = \frac{\sinh \pi z}{\pi} \sideset{}{^*}\sum_{n\in\Z} \frac{(-1)^n}{z-in} e^{2\pi inb}= \frac{\sinh \pi z}{\pi} \sideset{}{^*}\sum_{n\in\Z} \frac{e^{2\pi in(b+\frac{1}{2})}}{z-in}, 
\end{equation*}
where $\sum_{n\in\Z}^{*}$ means $\lim_{m\rightarrow \infty} \sum_{n=-m}^m$.
If we replace $z$ by $z-\frac{1}{2}i$, substitute $\nu=n+\frac{1}{2}$ and multiply both sides by $-e^{\pi ib}/\cosh\pi z$, we find
\begin{equation*}
-\frac{e^{2\pi bz}}{\cosh \pi z} = \frac{1}{\pi} \sideset{}{^*}\sum_{\nu\in \frac{1}{2} +\Z}\frac{e^{2\pi i\nu(b+\frac{1}{2})}} {z-i\nu}.
\end{equation*}
We have shown that $\sum_{\nu\in \frac{1}{2} +\Z}\frac{e^{2\pi i\nu(b+\frac{1}{2})}} {z-i\nu}$ converges, hence we get the desired result. \qed
\proof[ of Theorem \ref{period}]
Let $a\in (-\frac{1}{2},\frac{1}{2})$, then $g_{a+\frac{1}{2},b+\frac{1}{2}}(z) = \mathcal{O} \left(e^{-\pi \nu_0^2 \im(z)}\right)$ for $(\im(z)\rightarrow \infty)$, for some $\nu_0>0$ (for $a\in(-\frac{1}{2},0)$ we can take $\nu_0=a+\frac{1}{2}$, for $a\in [0,\frac{1}{2})$ we can take $\nu_0= \frac{1}{2}-a$). From this estimate, the (absolute) convergence of both $\int_{-\overline{\tau}}^{i\infty} \frac{g_{a+\frac{1}{2},b+\frac{1}{2}}(z)}{\sqrt{-i(z+\tau)}}\, dz$ and $\int_0^{i\infty} \frac{g_{a+\frac{1}{2},b+\frac{1}{2}}(z)}{\sqrt{-i(z+\tau)}}\, dz$ follows.\\
(1) We see
\begin{equation*}
\begin{split}
\int_{-\overline{\tau}}^{i\infty} &\frac{g_{a+\frac{1}{2},b+\frac{1}{2}}(z)}{\sqrt{-i(z+\tau)}}\, dz =
\int_{2iy}^{i\infty} \frac{g_{a+\frac{1}{2},b+\frac{1}{2}}(z-\tau)}{\sqrt{-iz}}\, dz = i\int_{2y}^{\infty} \frac{g_{a+\frac{1}{2},b+\frac{1}{2}}(iu-\tau)}{\sqrt{u}}\, du\\
&= i\int_{2y}^{\infty} u^{-\frac{1}{2}} \sum_{\nu\in a+\frac{1}{2}+\Z} \nu e^{\pi i\nu^2(iu-\tau)+2\pi i\nu(b+\frac{1}{2})}du \\
&= i \sum_{\nu\in a+\frac{1}{2}+\Z} \nu e^{-\pi i\nu^2\tau +2\pi i\nu (b+\frac{1}{2})} \int_{2y}^\infty u^{-\frac{1}{2}}e^{-\pi\nu^2 u} du\\
&= i \sum_{\nu\in a+\frac{1}{2}+\Z} \sign(\nu) e^{-\pi i\nu^2\tau +2\pi i\nu (b+\frac{1}{2})} \int_{2y\nu^2}^\infty v^{-\frac{1}{2}}e^{-\pi v} dv\\
&=i \sum_{\nu\in a+\frac{1}{2}+\Z} \left\{ \sign(\nu) -E\Bigl(\nu\sqrt{2y}\Bigr)\right\} e^{-\pi i\nu^2\tau +2\pi i\nu (b+\frac{1}{2})}\\
&= -e^{-\pi ia^2\tau +2\pi ia(b+\frac{1}{2})} \cdot\\
& \qquad \qquad \sum_{\nu\in \frac{1}{2}+\Z} \left\{ \sign(\nu+a) - E\Bigl((\nu+a)\sqrt{2y}\Bigr)\right\} (-1)^{\nu-\frac{1}{2}} e^{-\pi i \nu^2\tau -2\pi i\nu (a\tau-b)} 
\end{split}
\end{equation*}
If we use that for $a\in(-\frac{1}{2},\frac{1}{2})$ we have $\sign(\nu+a)=\sign(\nu)$ for all $\nu\in \frac{1}{2}+\Z$, we get the desired result. \\
(2) We see by means of Cauchy's theorem that for $a\in (-\frac{1}{2},\frac{1}{2})$ 
\begin{equation*}
h(x;\tau)=\int_{\R} \frac{e^{\pi i\tau z^2 -2\pi xz}}{\cosh\pi z}\,dz = \int_{-ia+\R} \frac{e^{\pi i\tau z^2 -2\pi xz}}{\cosh\pi z}\,dz,
\end{equation*}
hence
\begin{equation*}
\begin{split}
-e^{-\pi ia^2 \tau +2\pi ia(b+\frac{1}{2})} h(a\tau-b) &= -e^{\pi ia} \int_{-ia+\R} \frac{e^{\pi i\tau (z+ia)^2 +2\pi b(z+ia)}}{\cosh\pi z}\, dz \\
&= -e^{\pi ia} \int_\R \frac{e^{\pi i\tau z^2 +2\pi bz}}{\cosh\pi (z-ia)}\, dz.
\end{split}
\end{equation*}
If we replace $z$ by $z-ia$ in Lemma \ref{lem2}, multiply both sides by $e^{2\pi ia(b+\frac{1}{2})}$ and replace $\nu$ by $\nu-a$, we find 
\begin{equation*}
-e^{\pi ia} \frac{e^{2\pi bz}}{\cosh \pi (z-ia)} = \frac{1}{\pi} \sum_{\nu\in a+\frac{1}{2} +\Z} \frac{e^{2\pi i\nu(b+\frac{1}{2})}}{z-i\nu},
\end{equation*}
with $b\in(-\frac{1}{2},\frac{1}{2})$. Hence we find
\begin{equation*}
-e^{\pi ia} \frac{e^{2\pi bz}}{\cosh \pi (z-ia)} = -\frac{e^{\pi ia}}{\cos \pi a} +\frac{1}{\pi} \sum_{\nu\in a+\frac{1}{2} +\Z} e^{2\pi i\nu(b+\frac{1}{2})} \left( \frac{1}{z-i\nu}+\frac{1}{i\nu} \right).
\end{equation*}
This identity also follows from the theory of partial fraction decompositions given in \cite[pp.\ 134--136]{whit}. Using it we see
\begin{equation*}
\begin{split}
-&e^{-\pi ia^2 \tau +2\pi ia(b+\frac{1}{2})} h(a\tau -b)\\
&=\int_{\R} e^{\pi i\tau z^2} \biggl(-\frac{e^{\pi ia}}{\cos \pi a} +\frac{1}{\pi} \sum_{\nu\in a+\frac{1}{2} +\Z} e^{2\pi i\nu(b+\frac{1}{2})} \left( \frac{1}{z-i\nu}+\frac{1}{i\nu} \right)\biggr) dz\\
&= -\frac{e^{\pi ia}}{\cos \pi a} \frac{1}{\sqrt{-i\tau}} + \frac{1}{\pi} \int_{\R} \sum_{\nu\in a+\frac{1}{2} +\Z} e^{\pi i\tau z^2 +2\pi i\nu(b+\frac{1}{2})} \left( \frac{1}{z-i\nu}+\frac{1}{i\nu} \right) dz.
\end{split}
\end{equation*}
We want to change the order of summation and integration. This is allowed if
\begin{equation} \label{form4}
\int_{\R} \sum_{\nu\in a+\frac{1}{2} +\Z} \left|\ e^{\pi i\tau z^2 +2\pi i\nu(b+\frac{1}{2})} \left( \frac{1}{z-i\nu}+\frac{1}{i\nu} \right)\right| dz 
\end{equation}
converges. We have 
\begin{equation*}
\left|\ e^{\pi i\tau z^2 +2\pi i\nu(b+\frac{1}{2})} \left( \frac{1}{z-i\nu}+\frac{1}{i\nu} \right)\right| = e^{-\pi y z^2} \left|\frac{1}{z-i\nu}+\frac{1}{i\nu} \right| \leq \frac{|z|}{\nu^2} e^{-\pi y z^2}, 
\end{equation*}
with $y=\im(\tau)$. Both $\int_\R |z| e^{-\pi y z^2} dz$ and $\sum_{\nu\in a+\frac{1}{2} +\Z} \frac{1}{\nu^2}$ converge and hence the expression in \eqref{form4} converges. Now making the change of order we find
\begin{equation*}
\begin{split}
-&e^{-\pi ia^2 \tau +2\pi ia(b+\frac{1}{2})} h(a\tau -b)\\
&=-\frac{e^{\pi ia}}{\cos \pi a} \frac{1}{\sqrt{-i\tau}} + \frac{1}{\pi} \sum_{\nu\in a+\frac{1}{2} +\Z} e^{2\pi i\nu(b+\frac{1}{2})} \int_{\R} e^{\pi i\tau z^2} \left( \frac{1}{z-i\nu}+\frac{1}{i\nu} \right) dz.
\end{split}
\end{equation*}
Using Lemma \ref{lem1} we see
\begin{equation*}
\int_{\R} e^{\pi i\tau z^2} \left( \frac{1}{z-i\nu}+\frac{1}{i\nu} \right) dz = \pi \nu \int_0^{i\infty} \frac{e^{\pi i\nu^2z}}{\sqrt{-i(z+\tau)}}dz +\frac{1}{i\nu} \frac{1}{\sqrt{-i\tau}}.
\end{equation*}
Using partial integration we see that this equals
\begin{equation*}
-\frac{1}{2\nu} \int_0^{i\infty} \frac{e^{\pi i\nu^2 z}}{(-i(z+\tau))^{\frac{3}{2}}} dz,
\end{equation*}
so
\begin{equation}\label{form5}
\begin{split}
-&e^{-\pi ia^2 \tau +2\pi ia(b+\frac{1}{2})} h(a\tau -b)\\
&=-\frac{e^{\pi ia}}{\cos \pi a} \frac{1}{\sqrt{-i\tau}} - \frac{1}{2\pi} \sum_{\nu\in a+\frac{1}{2} +\Z} \int_0^{i\infty} \frac{e^{2\pi i\nu(b+\frac{1}{2})}}{\nu}  \frac{e^{\pi i\nu^2 z}}{(-i(z+\tau))^{\frac{3}{2}}} dz.
\end{split}
\end{equation}
We have
\begin{equation*}
\left| \frac{e^{2\pi i\nu(b+\frac{1}{2})}}{\nu}  \frac{e^{\pi i\nu^2 z}}{(-i(z+\tau))^{\frac{3}{2}}}\right| = \frac{1}{|\nu|}  \frac{e^{\pi i\nu^2 z}}{|z+\tau|^{\frac{3}{2}}}\leq \frac{e^{\pi i\nu^2 z}}{|\nu||\tau|^{\frac{3}{2}}},
\end{equation*}
and 
\begin{equation*}
\sum_{\nu\in a+\frac{1}{2} +\Z} \int_0^{i\infty} \frac{e^{\pi i\nu^2 z}}{|\nu||\tau|^{\frac{3}{2}}} |dz| =\frac{1}{\pi |\tau|^{\frac{3}{2}}} \sum_{\nu\in a+\frac{1}{2} +\Z} \frac{1}{|\nu|^3} <\infty.
\end{equation*}
Hence we can interchange the order of summation and integration in \eqref{form5}:
\begin{equation*}
\begin{split}
-&e^{-\pi ia^2 \tau +2\pi ia(b+\frac{1}{2})} h(a\tau -b)\\
&=-\frac{e^{\pi ia}}{\cos \pi a} \frac{1}{\sqrt{-i\tau}} - \frac{1}{2\pi} \int_0^{i\infty} \frac{1}{(-i(z+\tau))^{\frac{3}{2}}} \sum_{\nu\in a+\frac{1}{2} +\Z} \nu^{-1} e^{\pi i\nu^2 z +2\pi i\nu(b+\frac{1}{2})} dz.
\end{split}
\end{equation*}
Using partial integration we find
\begin{equation}\label{bijna}
\begin{split}
-&e^{-\pi ia^2 \tau +2\pi ia(b+\frac{1}{2})} h(a\tau -b)\\
&=-\frac{e^{\pi ia}}{\cos \pi a} \frac{1}{\sqrt{-i\tau}}+\frac{i}{\pi}\left. \frac{\sum_{\nu\in a+\frac{1}{2} +\Z} \nu^{-1} e^{\pi i\nu^2 z +2\pi i\nu(b+\frac{1}{2})}}{\sqrt{-i(z+\tau)}}\right|_0^{i\infty}\\
&-\frac{i}{\pi} \int_0^{i\infty} \frac{1}{\sqrt{-i(z+\tau)}}\biggl(\frac{d}{dz} \sum_{\nu\in a+\frac{1}{2} +\Z} \nu^{-1} e^{\pi i\nu^2 z +2\pi i\nu(b+\frac{1}{2})}\biggr) dz\\
&= -\frac{e^{\pi ia}}{\cos \pi a} \frac{1}{\sqrt{-i\tau}}-\frac{i}{\pi} \frac{1}{\sqrt{-i\tau}} \lim_{z\downarrow i0} \sum_{\nu\in a+\frac{1}{2} +\Z} \nu^{-1} e^{\pi i\nu^2 z +2\pi i\nu(b+\frac{1}{2})}\\
& + \int_0^{i\infty} \frac{g_{a+\frac{1}{2},b+\frac{1}{2}}(z)}{\sqrt{-i(z+\tau)}}\, dz .
\end{split}
\end{equation}
Let $z\in i\R_{\geq 0}$, and $\nu_0\in a+\frac{1}{2} +\Z$, $\nu_0>0$. By partial summation we find (with $T_\nu$ as in the proof of Lemma \ref{lem2})
\begin{equation*}
\underset{\nu\geq\nu_0}{\sum_{\nu\in a+\frac{1}{2} +\Z}} \nu^{-1} e^{\pi i\nu^2 z +2\pi i\nu(b+\frac{1}{2})} =
-\frac{e^{\pi i\nu_0^2 z}}{\nu_0} T_{\nu_0 -1} + \underset{\nu\geq\nu_0}{\sum_{\nu\in a+\frac{1}{2} +\Z}} T_\nu \biggl( \frac{e^{\pi i\nu^2 z}}{\nu} -\frac{e^{\pi i(\nu+1)^2 z}}{\nu+1}\biggr),
\end{equation*}
so 
\begin{align*}
&\Biggl| \underset{\nu\geq\nu_0}{\sum_{\nu\in a+\frac{1}{2} +\Z}} \nu^{-1} e^{\pi i\nu^2 z +2\pi i\nu(b+\frac{1}{2})}  \Biggr|\\
&\leq \frac{e^{\pi i\nu_0^2 z}}{\nu_0}\frac{1}{\left| e^{2\pi i(b+\frac{1}{2})}-1\right|}+ \underset{\nu\geq\nu_0}{\sum_{\nu\in a+\frac{1}{2} +\Z}} \frac{1}{\left| e^{2\pi i(b+\frac{1}{2})}-1\right|} \biggl( \frac{e^{\pi i\nu^2 z}}{\nu} -\frac{e^{\pi i(\nu+1)^2 z}}{\nu+1}\biggr)\\
&= \frac{e^{\pi i\nu_0^2 z}}{\nu_0} \frac{2}{\left| e^{2\pi i(b+\frac{1}{2})}-1\right|} \leq \frac{2}{\nu_0 \left| e^{2\pi i(b+\frac{1}{2})}-1\right|}.
\end{align*}
Hence 
\begin{equation*}
\underset{\nu\geq a+\frac{1}{2}}{\sum_{\nu\in a+\frac{1}{2} +\Z}} \nu^{-1} e^{\pi i\nu^2 z +2\pi i\nu(b+\frac{1}{2})}
\end{equation*}
converges uniformly for $z\in i\R_{\geq 0}$.

If we replace $\nu$ by $-\nu$ we find
\begin{equation*}
\underset{\nu\leq a-\frac{1}{2}}{\sum_{\nu\in a+\frac{1}{2} +\Z}} \nu^{-1} e^{\pi i\nu^2 z +2\pi i\nu(b+\frac{1}{2})} =-e^{2\pi i(a-\frac{1}{2})} \underset{\nu\geq -a+\frac{1}{2}}{\sum_{\nu\in a+\frac{1}{2} +\Z}} \nu^{-1} e^{\pi i\nu^2 z +2\pi i\nu(-b+\frac{1}{2})}, 
\end{equation*}
so
\begin{equation*}
\underset{\nu\leq a-\frac{1}{2}}{\sum_{\nu\in a+\frac{1}{2} +\Z}} \nu^{-1} e^{\pi i\nu^2 z +2\pi i\nu(b+\frac{1}{2})} 
\end{equation*}
and hence
\begin{equation*}
\sum_{\nu\in a+\frac{1}{2} +\Z} \nu^{-1} e^{\pi i\nu^2 z +2\pi i\nu(b+\frac{1}{2})} 
\end{equation*}
converge uniformly for $z\in i\R_{\geq 0}$.

Since the summand is continuous on $i\R_{\geq 0}$, so is
\begin{equation*}
z\mapsto \sum_{\nu\in a+\frac{1}{2} +\Z} \nu^{-1} e^{\pi i\nu^2 z +2\pi i\nu(b+\frac{1}{2})}.
\end{equation*}
Hence
\begin{equation*}
\lim_{z\downarrow i0} \sum_{\nu\in a+\frac{1}{2} +\Z} \nu^{-1} e^{\pi i\nu^2 z +2\pi i\nu(b+\frac{1}{2})} = \sum_{\nu\in a+\frac{1}{2} +\Z} \nu^{-1} e^{2\pi i\nu(b+\frac{1}{2})} = \pi i\frac{e^{\pi ia}}{\cos\pi a},
\end{equation*}
by Lemma \ref{lem2} with $z=0$.
If we put this into \eqref{bijna}, we get the desired result. \qed
\begin{remark}
We may also prove part (2) by using (2) of Proposition \ref{prop6} together with part (1):
We split the integral $\int_0^{i\infty}$ into the sum of $\int_{-\overline{\tau}}^{i\infty}$ and $\int_0^{-\overline{\tau}}$. For the first integral we have
\begin{equation}\label{split1}
\int_{-\overline{\tau}}^{i\infty} \frac{g_{a+\frac{1}{2},b+\frac{1}{2}}(z)}{\sqrt{-i(z+\tau)}}\, dz= -e^{-\pi ia^2 \tau +2\pi ia(b+\frac{1}{2})} R(a\tau-b),
\end{equation}
by (1). In the second integral we replace $z$ by $-\frac{1}{z}$:
\begin{equation*}
\begin{split}
\int_0^{-\overline{\tau}} \frac{g_{a+\frac{1}{2},b+\frac{1}{2}}(z)}{\sqrt{-i(z+\tau)}}\, dz&=
\int_{\frac{1}{\overline{\tau}}}^{i\infty} \frac{g_{a+\frac{1}{2},b+\frac{1}{2}}\left( -\frac{1}{z} \right)}{\sqrt{-i\left(-\frac{1}{z}+\tau\right)}} \frac{1}{(-iz)^2} dz\\
&= \frac{i}{\sqrt{-i\tau}} e^{2\pi i(a-\frac{1}{2})(b+\frac{1}{2})} \int_{\frac{1}{\overline{\tau}}}^{i\infty} \frac{g_{b+\frac{1}{2},-a+\frac{1}{2}}(z)}{\sqrt{-i\left(z-\frac{1}{\tau}\right)}}\, dz,
\end{split}
\end{equation*}
by (5) and (2) of Proposition \ref{prop141}. Using part (1) of this proposition with $(a,b,\tau)$ replaced by $(b,-a,-1/\tau)$ we see that this equals
\begin{equation}\label{split2}
-\frac{1}{\sqrt{-i\tau}} e^{\pi ib^2/\tau+\pi ia} R\left( \frac{a\tau-b}{\tau};-\frac{1}{\tau}\right).
\end{equation}
Combining \eqref{split1} and \eqref{split2} we see
\begin{equation*}
\begin{split}
&\int_0^{i\infty} \frac{g_{a+\frac{1}{2},b+\frac{1}{2}}(z)}{\sqrt{-i(z+\tau)}}\, dz\\
 &= -e^{-\pi ia^2 \tau +2\pi ia(b+\frac{1}{2})} R(a\tau-b) -\frac{1}{\sqrt{-i\tau}} e^{\pi ib^2/\tau+\pi ia} R\Bigl( \frac{a\tau-b}{\tau};-\frac{1}{\tau}\Bigr)\\
&= -e^{-\pi ia^2 \tau +2\pi ia(b+\frac{1}{2})}  \Bigl\{ R(a\tau-b;\tau) +\frac{1}{\sqrt{-i\tau}} e^{\pi i(a\tau-b)^2/\tau} R \Bigl( \frac{a\tau-b}{\tau};-\frac{1}{\tau}\Bigr)\Bigr\}\\
&= -e^{-\pi ia^2\tau +2\pi ia(b+\frac{1}{2})}h(a\tau-b),
\end{split}
\end{equation*}
by (2) of Proposition \ref{prop6}. 
\end{remark}

\clearemptydoublepage
\chapter{Indefinite $\theta$-functions}
\section{Introduction}

The classical theta series associated to a positive definite quadratic form $Q:\R^r \longrightarrow \R$ and $B:\R^r\times\R^r \longrightarrow \R$, the associated bilinear form $B(x,y)=Q(x+y)-Q(x)-Q(y)$, is the series
\begin{equation}\label{the1}
\Theta(z;\tau) := \sum_{n\in\Z^r} e^{2\pi i Q(n)\tau + 2\pi iB(n,z)}.
\end{equation}
These theta series have well-known transformation properties. In particular $\Theta(0;\tau)$ is a modular form of weight $r/2$. 

In \cite{gottsche} G\"ottsche and Zagier define a theta function for the case when the type of $Q$ is $(r-1,1)$. The definition of these functions is almost the same as in \eqref{the1}, only here the sum doesn't run over $\Z^r$, but some appropriate subset. However, in general, these functions do not have nice modular transformation properties. 

In this chapter we give a modified definition. We find elliptic and modular transformation properties for these functions. The theta functions we define depend not only on $Q$, but also on two vectors $c_1,c_2\in\R^r$ with $Q(c_i)\leq0$, $i=1,2$. The case $Q(c_1)=Q(c_2)=0$ gives the same functions as in \cite{gottsche}.

There is a connection between the indefinite $\theta$-functions from this chapter and certain $\theta$-functions considered by Siegel (see \cite{siegel}). However, I will not give this connection here.

\section{Definition of $\theta$}

Let $A$ be a symmetric $r\times r$-matrix with integer coefficients, which is non-degenerate.
We consider the quadratic form $Q:\C^r \longrightarrow \C$, $Q(x)=\frac{1}{2}\left<x,Ax\right>$ and  the associated bilinear form $B(x,y)=\left<x,Ay\right>= Q(x+y)-Q(x)-Q(y)$.

The \emph{type} of $Q$ is the pair $(r-s,s)$, where $s$ is the largest dimension of a linear subspace of $\R^r$ on which $Q$ is negative definite. The \emph{signature} of $Q$ is the number $r-2s$.  

From now on we assume that $s=1$, i.e., that $Q$ has type $(r-1,1)$. Then the set of vectors $c\in \R^r$ with $Q(c)<0$ has two components. If $B(c_1,c_2)<0$ then $c_1$ and $c_2$ belong to the same component, while if $B(c_1,c_2)>0$ then $c_1$ and $c_2$ belong to opposite components. Let $C_Q$ be one of the two components. If $c_0$ is a vector in that component, then $C_Q$ is given by:
\begin{equation*}
C_Q:=\{ c\in\R^r\mid Q(c)<0,\ B(c,c_0)<0\}.
\end{equation*}
We further set
\begin{equation*}
S_Q:=\{c\in\Z^r\mid c\text{ primitive},\ Q(c)=0,\ B(c,c_0)<0\}.
\end{equation*}
($c$ primitive means that the greatest common divisor of the components of $c$ is 1). The $(r-1)$-dimensional \emph{hyperbolic space} $C_Q/\R_+$ is the natural domain of definition of automorphic forms with respect to $O_A^+(\Z)$ (see section 2.4 for the definition), and $S_Q$ is a set of representatives for the corresponding set of \emph{cusps} 
\begin{equation*}
\{ c\in\Q^r\mid Q(c)=0,\ B(c,c_0)<0\}/\Q_+.
\end{equation*}
Note that $S_Q$ is empty in some cases, for example if $A=\left(\begin{smallmatrix}1&0\\0&-3\end{smallmatrix}\right)$. Further we put $\overline{C}_Q := C_Q \cup S_Q$. This is a generalisation of the usual construction $\overline{\Ha} = \Ha \cup \mathbb{P}^1(\Q)$, which is the special cone $\overline{C}_Q = C_Q \cup S_Q$ for the quadratic form $Q(a,b,c)=\frac{1}{2}(b^2-4ac)$.

For $c\in \overline{C}_Q$ put
\begin{align*}
R(c)&:=
\begin{cases}
\R^r &\text{if } c\in C_Q\\
\{ a\in\R^r \mid B(c,a)\not\in \Z\} & \text{if } c\in S_Q
\end{cases}\\
\intertext{and}
D(c)&:= \{ (z,\tau) \in \C^r\times \Ha \mid \im(z)/\im(\tau) \in R(c)\}. 
\end{align*}

\begin{definition}\label{deftheta}
Let $c_1,c_2\in \overline{C}_Q$. We define the \emph{theta function} of $Q$ with \emph{characteristics} $a\in R(c_1) \cap R(c_2)$ and $b\in\R^r$, with respect to $(c_1,c_2)$ by
\begin{equation*}
\theta_{a,b}(\tau)=\theta_{a,b}^{c_1,c_2}(\tau) := \sum_{\nu\in a+\Z^r} \rho(\nu;\tau) e^{2\pi iQ(\nu)\tau+2\pi iB(\nu,b)},
\end{equation*}
where $\rho(\nu;\tau)$ is defined by
\begin{equation*}
\rho(\nu;\tau)= \rho^{c_1,c_2}_A (\nu;\tau):=\rho^{c_1}(\nu;\tau)-\rho^{c_2}(\nu;\tau),
\end{equation*}
with
\begin{equation*}
\rho^{c}(\nu;\tau)= 
\begin{cases} E\left(\frac{B(c,\nu)}{\sqrt{-Q(c)}} y^{1/2}\right)& \text{if } c\in C_Q,\\
\sign(B(c,\nu))& \text{if } c\in S_Q,\\
\end{cases}
\end{equation*}
with $y=\im(\tau)$ and $E$ as in Definition \ref{defE}.

For $(z,\tau)\in D(c_1)\cap D(c_2)$, we define the theta function of $Q$ with respect to $(c_1,c_2)$ by
\begin{equation*}
\begin{split}
\theta(z;\tau)=\theta^{c_1,c_2}_A (z;\tau) &:= e^{-2\pi iQ(a)\tau -2\pi iB(a,b)} \theta_{a,b}(\tau)\\
&= \sum_{n\in\Z^r} \rho(n+a;\tau) e^{2\pi iQ(n)\tau+2\pi iB(n,z)},
\end{split}
\end{equation*}
with $a,b\in\R^r$ defined by $z=a\tau+b$, so $a=\frac{\im(z)}{\im(\tau)}$, $b=\frac{\im(\overline{z}\tau)}{\im(\tau)}$.
\end{definition}

\begin{remark}
The definition doesn't change if we replace $c_i$ by $\lambda c_i$, with $\lambda\in\R_+$. Hence we could replace the condition $Q(c_i)<0$ by $Q(c_i)=-1$. This would simplify the definition of $\rho$.
\end{remark}

\begin{remark}
In some special cases $\theta_{a,b}$ is holomorphic: if $c_1,c_2\in S_Q$ and, as we will see in section 2.5, also for some special values of $c_1$, $c_2$, $a$ and $b$.
In general however, the functions $\theta$ and $\theta_{a,b}$ are not holomorphic. 
\end{remark}

Because $Q$ is indefinite, $e^{2\pi iQ(n)\tau}$ isn't bounded. Therefore it's not immediately clear that the series defining $\theta(z;\tau)$ converges absolutely. However, using an estimate for the growth of $\rho$ we shall find:
\begin{proposition}\label{lempje}
The series defining $\theta(z;\tau)$ converges absolutely. 
\end{proposition}
For the proof of this proposition, we need two lemmas
\begin{lemma}\label{sublem1}
Let $c\in C_Q$. The quadratic form $Q_c:\R^r \longrightarrow \R$, $Q_c(\nu):=Q(\nu) -\frac{B(c,\nu)^2}{2Q(c)}$ is positive definite, and we have
\begin{equation*}
Q_c(\nu) \geq \lambda_{c,c_0} Q_{c_0}(\nu) \qquad \forall \nu\in\R^r,
\end{equation*}
with
\begin{equation*}
\lambda_{c,c_0} =\frac{B(c,c_0)^2 -2Q(c)Q(c_0) - |B(c,c_0)| \sqrt{B(c,c_0)^2-4Q(c)Q(c_0)}}{2Q(c)Q(c_0)}>0.
\end{equation*}
\end{lemma}
\proof 
If $\nu\in\R^r$ is linearly independent of $c$, the quadratic form $Q$ has type $(1,1)$ on $\spa\{c,\nu\}$; hence the matrix
\begin{equation*}
\left(\begin{matrix} 2Q(c)&B(c,\nu)\\
B(c,\nu)&2Q(\nu)\\
\end{matrix}\right)
\end{equation*}
has determinant $<0$, so $4Q(\nu)Q(c)-B(c,\nu)^2<0$.  Rewriting gives 
\begin{equation*} 
 Q(\nu) -\frac{B(c,\nu)^2}{2Q(c)}> -\frac{B(c,\nu)^2}{4Q(c)}\geq0.
\end{equation*}
If $\nu=\lambda c$, with $\lambda\neq 0$, we get $Q(\nu) -\frac{B(c,\nu)^2}{2Q(c)}=-Q(c)\lambda^2>0$, which proves that $Q_c$ is positive definite.  

For the second part we consider the restriction of $Q_c$ to the ellipsoid $S=\{ \nu\in\R^r | Q_{c_0} (\nu) =1 \}$. Since $S$ is compact, $Q_c|_S$ assumes its absolute minimum at some point $\nu_0$. We compute that minimum with the method of Lagrange multipliers: There is a real number $\lambda$, such that
\begin{equation*}
\nabla Q_c (\nu_0) = \lambda \nabla Q_{c_0} (\nu_0),
\end{equation*}
or equivalently
\begin{equation}\label{lag}
A\left( \nu_0 -\frac{B(\nu_0,c)}{Q(c)} c\right) =\lambda A\left( \nu_0 -\frac{B(\nu_0,c_0)}{Q(c_0)} c_0\right).
\end{equation}
Taking the inner product with $c$ on both sides of \eqref{lag} we find
\begin{equation}\label{lagje1}
-B(\nu_0,c)=\lambda \left(B(\nu_0,c) -\frac{B(\nu_0,c_0)B(c,c_0)}{Q(c_0)}\right).
\end{equation}
Taking the inner product with $c_0$ on both sides of \eqref{lag} we find
\begin{equation}\label{lagje2}
B(\nu_0,c_0) -\frac{B(\nu_0,c)B(c,c_0)}{Q(c)}=-\lambda B(\nu_0,c_0).
\end{equation}
Combining \eqref{lagje1} and \eqref{lagje2} we find
\begin{equation}\label{pol}
Q(c)Q(c_0)(\lambda+1)^2 =\lambda B(c,c_0)^2
\end{equation}
or
\begin{equation*}
B(\nu_0,c)=B(\nu_0,c_0)=0.
\end{equation*}
If $B(\nu_0,c)=B(\nu_0,c_0)=0$, then \eqref{lag} reduces to $A\nu_0=\lambda A\nu_0$, from which we find $\lambda=1$.

The roots of \eqref{pol} are
\begin{equation*}
\lambda^\pm =\frac{B(c,c_0)^2 -2Q(c)Q(c_0) \pm |B(c,c_0)| \sqrt{B(c,c_0)^2-4Q(c)Q(c_0)}}{2Q(c)Q(c_0)}.
\end{equation*}

Taking the inner product with $\nu_0$ on both sides of \eqref{lag} and dividing by 2, we find
\begin{equation*}
Q_c(\nu_0) = \lambda Q_{c_0} (\nu_0) =\lambda.
\end{equation*}
Hence the absolute minimum of $Q_c|_S$ is the minimum of $\{1,\lambda^-,\lambda^+\}$ which is $\lambda^- = \lambda_{c,c_0}$. So we have
\begin{equation*}
Q_c(\nu) \geq \lambda_{c,c_0} \qquad \forall \nu\in S.
\end{equation*}
Let $\nu\in\R^r$, $\nu\not= 0$, then $\frac{\nu}{\sqrt{Q_{c_0}(\nu)}}\in S$, so
\begin{equation*}
\lambda_{c,c_0} \leq Q_c\left(\frac{\nu}{\sqrt{Q_{c_0}(\nu)}}\right) = \frac{Q_c(\nu)}{Q_{c_0}(\nu)}.
\end{equation*}
Multiplying both sides by $Q_{c_0}(\nu)$ we get the desired result. \qed
\begin{lemma}\label{sublem2}
Let $c_1,c_2\in C_Q$ be linearly independent. The quadratic form $Q^+:\R^r \longrightarrow \R$, $Q^+(\nu):=Q(\nu) +\frac{B(c_1,c_2)}{4Q(c_1)Q(c_2)-B(c_1,c_2)^2} B(c_1,\nu) B(c_2,\nu)$ is positive definite.
\end{lemma}
\proof
If $\nu\in\R^r$ is not a linear combination of $c_1$ and $c_2$, the quadratic form $Q$ has type $(2,1)$ on $\spa\{c_1,c_2,\nu\}$; so the matrix
\begin{equation}\label{mat}
\left(\begin{matrix} 2Q(c_1)&B(c_1,c_2)&B(c_1,\nu)\\
B(c_1,c_2)&2Q(c_2)&B(c_2,\nu)\\
B(c_1,\nu)&B(c_2,\nu)&2Q(\nu)
\end{matrix}\right)
\end{equation}
has determinant $<0$.  Rewriting gives 
\begin{equation*}
Q^+(\nu) > \frac{Q(c_2)B(c_1,\nu)^2 + Q(c_1) B(c_2,\nu)^2}{4Q(c_1)Q(c_2)-B(c_1,c_2)^2}\geq 0.
\end{equation*}
If $\nu\in\R^r$ is a linear combination of $c_1$ and $c_2$ the determinant of the matrix in \eqref{mat} is zero. Hence 
\begin{equation*}
Q^+(\nu) = \frac{Q(c_2)B(c_1,\nu)^2 + Q(c_1) B(c_2,\nu)^2}{4Q(c_1)Q(c_2)-B(c_1,c_2)^2}.
\end{equation*}
So if $Q^+(\nu)=0$, we have $B(c_1,\nu)=B(c_2,\nu)=0$, which implies $\nu=0$. Thus if $\nu\not=0$ then $Q^+(\nu)$ is strictly positive, so $Q^+$ is positive definite.
\proof[ of Proposition \ref{lempje}] If $c_1,c_2\in C_Q$, we write $\rho(\nu;\tau)$, using Lemma \ref{lem131}, as the sum of the three expressions
\begin{align}
-&\sign\Bigl(B(c_1,\nu)\Bigr) \beta\left( -\frac{B(c_1,\nu)^2}{Q(c_1)}y\right),\label{al1}\\
&\sign\Bigl(B(c_2,\nu)\Bigr) \beta\left( -\frac{B(c_2,\nu)^2}{Q(c_2)}y\right)\label{al2}\\
\intertext{and}
&\sign\Bigl(B(c_1,\nu)\Bigr) - \sign\Bigl(B(c_2,\nu)\Bigr),\label{al3}
\end{align}
with $\beta$ as in Lemma \ref{lem131}. If $c_1\in C_Q$ and $c_2\in S_Q$ we get only the sum of the first and the last expression. If $c_1\in S_Q$ and $c_2\in C_Q$ we get the sum of the last two expressions. If $c_1,c_2\in S_Q$ we have only the last expression. Hence the proof is reduced to showing that the series
\begin{equation}\label{221}
\sum_{\nu\in a+\Z^r} \sign\Bigl(B(c,\nu)\Bigr) \beta\left( -\frac{B(c,\nu)^2}{Q(c)}y\right)e^{2\pi iQ(\nu)\tau+2\pi iB(\nu,b)}
\end{equation}
converges absolutely for all $c$ with $Q(c)<0$, and that the series
\begin{equation}\label{222}
\sum_{\nu\in a+\Z^r} \Bigl\{\sign\Bigl(B(c_1,\nu)\Bigr) - \sign\Bigl(B(c_2,\nu)\Bigr)\Bigr\} e^{2\pi iQ(\nu)\tau+2\pi iB(\nu,b)}
\end{equation}
converges absolutely for all $c_1,c_2 \in \overline{C}_Q$.

We will first show that the series \eqref{221} converges absolutely for all $c$ with $Q(c)<0$:
We can easily see that $0\leq\beta(x)\leq e^{-\pi x}$ for all $x\in\R_{\geq0}$; hence if $Q(c)<0$
\begin{equation}\label{223}
\begin{split}
&\left| \sign\Bigl(B(c,\nu)\Bigr) \beta\left( -\frac{B(c,\nu)^2}{Q(c)}y\right)e^{2\pi iQ(\nu)\tau+2\pi iB(\nu,b)}\right|\\
&\leq e^{\pi \frac{B(c,\nu)^2}{Q(c)}y} \left| e^{2\pi iQ(\nu)\tau+2\pi iB(\nu,b)}\right|\\
&= e^{-2\pi \left( Q(\nu) -\frac{B(c,\nu)^2}{2Q(c)}\right)y}.
\end{split}
\end{equation}

Using Lemma \ref{sublem1}, we see that the series
\begin{equation*}
\sum_{\nu\in a+\Z^r} e^{-2\pi \left( Q(\nu) -\frac{B(c,\nu)^2}{2Q(c)}\right)y}
\end{equation*}
converges, and so the series \eqref{221} converges absolutely if $Q(c)<0$.

We will now show that the series \eqref{222} converges absolutely for all $c_1,c_2 \in \overline{C}_Q$:\\
If $c_1$ and $c_2$ are linearly dependent, we have $\rho^{c_1,c_2} =0$. Hence we can assume that they are linearly independent.

\noindent{\it Case $1$}: $c_1,c_2\in C_Q$.

If we have $B(c_1,\nu)B(c_2,\nu)>0$, then $\sign(B(c_1,\nu)) - \sign(B(c_2,\nu))=0$. If we have $B(c_1,\nu)B(c_2,\nu)\leq 0$, then (note that $4Q(c_1)Q(c_2)-B(c_1,c_2)^2<0$, as we saw before)
\begin{equation*}
Q(\nu) \geq Q^+(\nu),
\end{equation*}
with $Q^+$ as in Lemma \ref{sublem2}. Hence we find
\begin{equation}\label{nouja}
\begin{split}
&\left|\Bigl\{\sign\Bigl(B(c_1,\nu)\Bigr) - \sign\Bigl(B(c_2,\nu)\Bigr)\Bigr\} e^{2\pi iQ(\nu)\tau+2\pi iB(\nu,b)}\right|\\
&= \Bigl|\sign\Bigl(B(c_1,\nu)\Bigr) - \sign\Bigl(B(c_2,\nu)\Bigr)\Bigr|\  e^{-2\pi Q(\nu)y}\\
&\leq 2 e^{-2\pi Q^+(\nu)y}.
\end{split}
\end{equation}
Using Lemma \ref{sublem2}, we see that the series
\begin{equation*}
\sum_{\nu\in a+\Z^r} e^{-2\pi Q^+(\nu)y}
\end{equation*}
converges, and so the series \eqref{222} converges absolutely.

\noindent{\it Case $2$}: $c_1\in C_Q$ and $c_2\in S_Q$.

We can assume that $c_1\in C_Q\cap \Z^r$, since otherwise we pick any $c_1'\in C_Q\cap \Z^r$, write
\begin{equation*}
\begin{split}
\sum_{\nu\in a+\Z^r} &\Bigl\{\sign\Bigl(B(c_1,\nu)\Bigr) - \sign\Bigl(B(c_2,\nu)\Bigr)\Bigr\}\ e^{2\pi iQ(\nu)\tau+2\pi iB(\nu,b)}\\
&=\sum_{\nu\in a+\Z^r} \Bigl\{\sign\Bigl(B(c_1,\nu)\Bigr) - \sign\Bigl(B(c_1',\nu)\Bigr)\Bigr\}\ e^{2\pi iQ(\nu)\tau+2\pi iB(\nu,b)}\\
&+\sum_{\nu\in a+\Z^r} \Bigl\{\sign\Bigl(B(c_1',\nu)\Bigr) - \sign\Bigl(B(c_2,\nu)\Bigr)\Bigr\}\ e^{2\pi iQ(\nu)\tau+2\pi iB(\nu,b)},
\end{split}
\end{equation*}
and use that
\begin{equation*}
\sum_{\nu\in a+\Z^r} \Bigl\{\sign\Bigl(B(c_1,\nu)\Bigr) - \sign\Bigl(B(c_1',\nu)\Bigr)\Bigr\}\ e^{2\pi iQ(\nu)\tau+2\pi iB(\nu,b)}
\end{equation*}
converges absolutely.

We write $\nu=\mu+nc_2$ with $\mu\in a+\Z^r$ and $n\in\Z$, such that $\frac{B(c_1,\mu)}{B(c_1,c_2)}\in [0,1)$ (we see $n=\left[ \frac{B(c_1,\nu)}{B(c_1,c_2)}\right]$). Then $\sign(B(c_2,\nu))=\sign(B(c_2,\mu))$ (use $B(c_2,c_2)=0$) and $\sign(B(c_1,\nu))= -\sign(n+\frac{B(c_1,\mu)}{B(c_1,c_2)})$.
Hence
\begin{equation*}
\begin{split}
\sum_{\nu\in a+\Z^r} &\Bigl\{\sign\Bigl(B(c_1,\nu)\Bigr) - \sign\Bigl(B(c_2,\nu)\Bigr)\Bigr\}\ e^{2\pi iQ(\nu)\tau+2\pi iB(\nu,b)}\\
= &- \underset{\frac{B(c_1,\mu)}{B(c_1,c_2)}\in[0,1)}{\sum_{\mu\in a+\Z^r}} \sum_{n\in\Z} \Bigl\{\sign\Bigl(B(c_2,\mu)\Bigr) + \sign\Bigl(n+\frac{B(c_1,\mu)}{B(c_1,c_2)}\Bigr)\Bigr\}\ \cdot\\
& \qquad \qquad \qquad \qquad \qquad \qquad \cdot e^{2\pi iQ(\mu)\tau+2\pi iB(c_2,\mu)n\tau+2\pi iB(\mu,b)+2\pi iB(c_2,b)n}
\end{split}
\end{equation*}
Using 
\begin{equation*}
\frac{1}{1-x}= \begin{cases} \sum_{n=0}^{\infty} x^n &\text{ if}\ |x|<1\\
-\sum_{n=-\infty}^{-1} x^n &\text{ if}\ |x|>1,
\end{cases} 
\end{equation*}
we see
\begin{equation*}
\begin{split}
\sum_{n\in\Z} &\Bigl\{\sign\Bigl(B(c_2,\mu)\Bigr) + \sign\Bigl(n+\frac{B(c_1,\mu)}{B(c_1,c_2)}\Bigr)\Bigr\}\ e^{2\pi iB(c_2,\mu)n\tau+2\pi iB(c_2,b)n}\\
&= \frac{2}{1-e^{2\pi iB(c_2,\mu)\tau +2\pi iB(c_2,b)}} -\delta(B(c_1,\mu)).
\end{split}
\end{equation*}
Here we used that 
\begin{equation}\label{sch}
B(c_2,\mu)\geq B(c_2,\tilde{\mu})>0,
\end{equation}
for all $\mu\in a+\Z^r$, and for some $\tilde{\mu}\in a+\Z^r$. This is guaranteed by the fact that $(z,\tau)\in D(c_2)$.

Since $c_1,c_2\in\Z^r$ we have 
\begin{equation*}
\left\{ \mu\in a+\Z^r \left| \frac{B(c_1,\mu)}{B(c_1,c_2)} \in [0,1)\right.\right\} =
\bigcup_{\mu_0\in P_0} \left( \mu_0 + \left<c_1\right>_\Z^{\perp} \right),
\end{equation*}
for a suitable finite set $P_0$, with $\left<c_1\right>_\Z^{\perp}:=\{ \xi\in\Z^r \mid B(c_1,\xi)=0\}$. So
\begin{equation*}
\begin{split}
\sum_{\nu\in a+\Z^r} &\Bigl\{\sign\Bigl(B(c_1,\nu)\Bigr) - \sign\Bigl(B(c_2,\nu)\Bigr)\Bigr\}\ e^{2\pi iQ(\nu)\tau+2\pi iB(\nu,b)}\\
&= - \sum_{\mu_0\in P_0} \sum_{\xi\in\left<c_1\right>^\perp} \left\{\frac{2}{1-e^{2\pi iB(c_2,\xi+\mu_0)\tau +2\pi iB(c_2,b)}} -\delta(B(c_1,\mu_0))\right\} \cdot\\
&\hspace{8cm} e^{2\pi iQ(\xi+\mu_0)\tau +2\pi iB(\xi+\mu_0,b)}.
\end{split}
\end{equation*}
This series converges absolutely, since $Q$ is positive definite on $\left<c_1\right>_\Z^\perp$, and the term
\begin{equation*}
\frac{2}{1-e^{2\pi iB(c_2,\xi+\mu_0)\tau +2\pi iB(c_2,b)}} -\delta(B(c_1,\mu_0))
\end{equation*}
is bounded (use \eqref{sch}). 

\noindent{\it Case $3$}: $c_1\in S_Q$ and $c_2\in C_Q$.

Since $\theta^{c_1,c_2}=-\theta^{c_2,c_1}$, this follows directly from the previous case.

\noindent{\it Case $4$}: $c_1,c_2 \in S_Q$.

Since $\theta^{c_1,c_2}= \theta^{c_1,c_3}+\theta^{c_3,c_2}$, for arbitrary $c_3\in C_Q$, this follows directly from case 2 and 3. \qed

\section{Properties of the $\theta$-functions}

The theta functions in Definition \ref{deftheta} have some nice elliptic and modular transformation properties, similar to those of the theta functions associated to positive definite quadratic forms.

\begin{proposition}\label{thet}
The function $\theta$ satisfies:
\begin{itemize}
\item[(1)] For $c_1,c_2,c_3\in \overline{C}_Q$ and $(z,\tau)\in D(c_1)\cap D(c_2) \cap D(c_3)$ we have the cocycle conditions $\theta^{c_1,c_2}+\theta^{c_2,c_1}=0$ and $\theta^{c_1,c_2}+\theta^{c_2,c_3}+\theta^{c_3,c_1}=0$.
\item[(2)] $\theta(z+\lambda\tau+\mu;\tau)=e^{-2\pi iQ(\lambda)\tau -2\pi iB(z,\lambda)}\theta(z;\tau)$ for all $\lambda\in\Z^r$ and $\mu\in A^{-1} \Z^r$.
\item[(3)] $\theta(-z;\tau)=-\theta(z;\tau)$.
\item[(4)] The function $(c_1,c_2) \mapsto \theta^{c_1,c_2}$ is continuous on $C_Q \times C_Q$.
\item[(5)] Let $c_1,c_3\in C_Q$, $c_2\in S_Q$ and $(z,\tau)\in D(c_2)$. Set $c(t)=c_2+tc_3$. Then $c(t)\in C_Q$ for all $t\in(0,\infty)$ and $\lim_{t\downarrow 0} \theta^{c_1,c(t)}(z;\tau)=\theta^{c_1,c_2}(z;\tau)$.
\item[(6)] $\theta(z;\tau+1)=\theta(z+\frac{1}{2}A^{-1}A^*;\tau)$ with $A^*=(A_{11}\ldots A_{rr})^T\in\Z^r$, the vector of diagonal elements of $A$. In particular, $\theta(z;\tau+2)=\theta(z;\tau)$ and $\theta(z;\tau+1)=\theta(z;\tau)$ if the matrix $A$ is even.
\item[(7)] Let $D'(c):= \left\{(z,\tau)\in D(c) \mid \left(\frac{z}{\tau},-\frac{1}{\tau}\right)\in D(c) \right\} = \{(a\tau+b,\tau)\mid \tau\in \Ha,\ a,b\in \R^r,\ B(c,a)\not\in\Z,\ B(c,b)\not\in\Z\}$. If $(z,\tau)\in D'(c_1)\cap D'(c_2)$ then
\begin{equation*}
\theta\Bigl(\frac{z}{\tau};-\frac{1}{\tau}\Bigr)= \frac{i}{\sqrt{-\det A}} (-i\tau)^{r/2} \sum_{p\in A^{-1}\Z^r /\Z^r} e^{2\pi i Q(z+p\tau)/\tau}\theta(z+p\tau;\tau).
\end{equation*}
\end{itemize}
\end{proposition}
\proof (1) follows from the corresponding relations for $\rho^{c_1,c_2}$.\\
(2) The identity $\theta(z+\mu;\tau)=\theta(z;\tau)$ for $\mu\in A^{-1} \Z^r$ is easy, and we find $\theta(z+\lambda\tau;\tau)=e^{-2\pi iQ(\lambda)\tau -2\pi iB(z,\lambda)}\theta(z;\tau)$ for $\lambda\in\Z^r$ when we replace $n$ by $n+\lambda$ in the definition. 

For (3), replace $n$ by $-n$ in the definition and use that $E$ and $\sign$ are odd functions. \\
(4) We show that $c_1 \mapsto \theta^{c_1,c_2}$ is continuous on $C_Q$. The result then follows from (1).

Using the decomposition of $\rho$ as the sum of \eqref{al1}, \eqref{al2} and \eqref{al3} we see that it's sufficient to prove that
\begin{equation}\label{fun2}
c\mapsto \sum_{\nu\in a+\Z^r} \sign\Bigl(B(c,\nu)\Bigr) \beta\left( -\frac{B(c,\nu)^2}{Q(c)}y\right) e^{2\pi iQ(\nu)\tau+2\pi iB(\nu,b)},
\end{equation}
and
\begin{equation}\label{fun1}
c_1\mapsto \sum_{\nu\in a+\Z^r} \Bigl\{\sign\Bigl(B(c_1,\nu)\Bigr) - \sign\Bigl(B(c_2,\nu)\Bigr)\Bigr\}\ e^{2\pi iQ(\nu)\tau+2\pi iB(\nu,b)}
\end{equation}
are continuous on $C_Q$.

Using Lemma \ref{sublem1} and \eqref{223} we see 
\begin{equation*}
\begin{split}
\left| \sign\Bigl(B(c,\nu)\Bigr) \beta\left( -\frac{B(c,\nu)^2}{Q(c)}y\right)e^{2\pi iQ(\nu)\tau+2\pi iB(\nu,b)}\right| &\leq e^{-2\pi \left( Q(\nu) -\frac{B(c,\nu)^2}{2Q(c)}\right)y}\\
& \leq e^{-2\pi \lambda_{c,c_0} Q_{c_0} (\nu) y}
\end{split}
\end{equation*}
Since $c\mapsto \lambda_{c,c_0}$ is continuous and $\lambda_{c,c_0}>0$ for all $c\in C_Q$, we can find an neighbourhood $\EuScript{N}_c$ of $c$ such that $\lambda_{c,c_0} \geq \epsilon >0$ for all $c\in \EuScript{N}_c$. Hence on $\EuScript{N}_c$ we find
\begin{equation*}
\left| \sign\Bigl(B(c,\nu)\Bigr) \beta\left( -\frac{B(c,\nu)^2}{Q(c)}y\right)e^{2\pi iQ(\nu)\tau+2\pi iB(\nu,b)}\right| \leq e^{-2\pi \epsilon Q_{c_0} (\nu) y}.
\end{equation*}
The series
\begin{equation*}
\sum_{\nu\in a+\Z^r} e^{-2\pi \epsilon Q_{c_0} (\nu) y}
\end{equation*}
converges, and so the series in \eqref{fun2} converges uniformly for $c$ in $\EuScript{N}_c$. Hence the function in \eqref{fun2} is continuous on $\EuScript{N}_c$. Since this holds for all $c\in C_Q$, the function in \eqref{fun2} is continuous on $C_Q$.

In \eqref{nouja} we have seen that
\begin{equation*}
\left|\Bigl\{\sign\Bigl(B(c_1,\nu)\Bigr) - \sign\Bigl(B(c_2,\nu)\Bigr)\Bigr\}\ e^{2\pi iQ(\nu)\tau+2\pi iB(\nu,b)}\right| \leq 2 e^{-2\pi Q^+(\nu)y}.
\end{equation*}
The function $Q^+$ restricted to the sphere $S=\{\nu\in\R^r |\ \| \nu \| =1 \}$ assumes its absolute minimum $\lambda(c_1) >0$. Hence 
\begin{equation*}
Q^+(\nu)\geq \lambda(c_1) \qquad \forall \nu \in S,
\end{equation*}
and so 
\begin{equation*}
Q^+(\nu)\geq \lambda(c_1)\|\nu\|^2 \qquad \forall \nu \in \R^r.
\end{equation*}
Since $c_1\mapsto \lambda(c_1)$ is continuous and $\lambda(c_1)>0$ for all $c_1\in C_Q$, we can find an neighbourhood $\EuScript{N}_{c_1}$ of $c_1$ such that $\lambda(c_1) \geq \epsilon >0$ for all $c_1\in \EuScript{N}_{c_1}$. Hence on $\EuScript{N}_{c_1}$ we find
\begin{equation*}
\left|\Bigl\{\sign\Bigl(B(c_1,\nu)\Bigr) - \sign\Bigl(B(c_2,\nu)\Bigr)\Bigr\}\  e^{2\pi iQ(\nu)\tau+2\pi iB(\nu,b)}\right| \leq 2 e^{-2\pi \epsilon \|\nu\|^2 y}.
\end{equation*}
The series 
\begin{equation*}
\sum_{\nu\in a+\Z^r} e^{-2\pi \epsilon \|\nu\|^2 y}
\end{equation*}
converges, and so the series in \eqref{fun1} converges uniformly for $c_1$ in $\EuScript{N}_{c_1}$. Hence the function in \eqref{fun1} is continuous on $\EuScript{N}_{c_1}$. Since this holds for all $c_1\in C_Q$, the function in \eqref{fun1} is continuous on $C_Q$.\\
(5) Note that $\theta^{c_1,c(t)}=\theta^{c_1,c_3}+\theta^{c_3,c(t)}$. We can therefore assume $c_3$ to be equal to $c_1$.
We have $Q(c(t))=Q(c_2+tc_1)=tB(c_1,c_2)+t^2Q(c_1)<0$ and $B(c_1,c(t))=B(c_1,c_2)+2tQ(c_1)<0$ for all $t\in(0,\infty)$, since $B(c_1,c_2)<0$ and $Q(c_1)<0$. Hence $c(t)\in C_Q$ for all $t\in(0,\infty)$.

Using $\theta^{c_1,c(t)}=\theta^{c_1,c_2}+\theta^{c_2,c(t)}$ and the decomposition of $\rho$ as the sum of \eqref{al1}, \eqref{al2} and \eqref{al3} we see that it's sufficient to prove that
\begin{equation}\label{eerste}
\lim_{t\downarrow 0} \sum_{\nu\in a+\Z^r} \Bigl\{\sign\Bigl(B(c_2,\nu)\Bigr) - \sign\Bigl(B(c(t),\nu)\Bigr)\Bigr\}\ e^{2\pi iQ(\nu)\tau+2\pi iB(\nu,b)}=0,
\end{equation}
and
\begin{equation}\label{232}
\lim_{t\downarrow 0} \sum_{\nu\in a+\Z^r} \sign\Bigl(B(c(t),\nu)\Bigr) \beta\left( -\frac{B(c(t),\nu)^2}{Q(c(t))}y\right) e^{2\pi iQ(\nu)\tau+2\pi iB(\nu,b)}=0.
\end{equation}
It is easy to see that
\begin{equation*}
\Bigl|\sign\Bigl(B(c_2,\nu)\Bigr) - \sign\Bigl(B(c(t),\nu)\Bigr)\Bigr|\leq \Bigl|\sign\Bigl(B(c_1,\nu)\Bigr) - \sign\Bigl(B(c_2,\nu)\Bigr)\Bigr|
\end{equation*}
for all $\nu\in a+\Z^r$ and $t\in(0,\infty)$ (Both sides can take on the values 0,1 and 2. If the right hand side is 0, then $\sign(B(c_1,\nu))=\sign(B(c_2,\nu))$, so $\sign(B(c_2,\nu))=\sign(B(c(t),\nu))$. Hence the left hand side is also 0, and the equation holds. If the right hand side is 1, then either $B(c_1,\nu)$ or $B(c_2,\nu)$ is zero. If $B(c_2,\nu)=0$ we get that the left hand side equals the right hand side. If $B(c_1,\nu)=0$ the left hand side equals 0). Hence
\begin{equation*}
\begin{split}
&\left| \Bigl\{\sign\Bigl(B(c_2,\nu)\Bigr) - \sign\Bigl(B(c(t),\nu)\Bigr)\Bigr\}\ e^{2\pi iQ(\nu)\tau+2\pi iB(\nu,b)} \right| \\
&\qquad \leq \left| \Bigl\{\sign\Bigl(B(c_1,\nu)\Bigr) - \sign\Bigl(B(c_2,\nu)\Bigr)\Bigr\}\ e^{2\pi iQ(\nu)\tau+2\pi iB(\nu,b)} \right|,
\end{split}
\end{equation*}
for all $\nu\in a+\Z^r$ and $t\in(0,\infty)$. In the proof of Proposition \ref{lempje} (Case 2) we have seen that \eqref{222} converges absolutely, i.e. 
\begin{equation*}
\sum_{\nu\in a+\Z^r} \left| \Bigl\{\sign\Bigl(B(c_1,\nu)\Bigr) - \sign\Bigl(B(c_2,\nu)\Bigr)\Bigr\}\ e^{2\pi iQ(\nu)\tau+2\pi iB(\nu,b)} \right|
\end{equation*}
converges. Hence 
\begin{equation*}
\sum_{\nu\in a+\Z^r} \Bigl\{\sign\Bigl(B(c_2,\nu)\Bigr) - \sign\Bigl(B(c(t),\nu)\Bigr)\Bigr\}\ e^{2\pi iQ(\nu)\tau+2\pi iB(\nu,b)}
\end{equation*}
converges uniformly for $t\in (0,\infty)$. Using this we find 
\begin{equation*}
\begin{split}
\lim_{t\downarrow 0} &\sum_{\nu\in a+\Z^r} \Bigl\{\sign\Bigl(B(c_2,\nu)\Bigr) - \sign\Bigl(B(c(t),\nu)\Bigr)\Bigr\}\ e^{2\pi iQ(\nu)\tau+2\pi iB(\nu,b)}\\
&= \sum_{\nu\in a+\Z^r} \lim_{t\downarrow 0} \Bigl\{\sign\Bigl(B(c_2,\nu)\Bigr) - \sign\Bigl(B(c(t),\nu)\Bigr)\Bigr\}\ e^{2\pi iQ(\nu)\tau+2\pi iB(\nu,b)}=0.
\end{split}
\end{equation*}
This proves \eqref{eerste}.

We will now prove \eqref{232}: Using \eqref{223}, we see that 
\begin{equation*}
\left| \sign\Bigl(B(c(t),\nu)\Bigr) \beta\left( -\frac{B(c(t),\nu)^2}{Q(c(t))}y\right) e^{2\pi iQ(\nu)\tau+2\pi iB(\nu,b)} \right| \leq e^{-2\pi \left( Q(\nu) -\frac{B(c(t),\nu)^2}{2Q(c(t))}\right)y}.
\end{equation*}
We write $a+\Z^r$ as the union of $P_1$, $P_2$ and $P_3$, with
\begin{equation*}
\begin{split}
P_1:= \{ \nu\in a+\Z^r \mid & \sign (B(c_2,\nu)) =-\sign(B(c_1,\nu)) \}\\
P_2:= \{ \nu\in a+\Z^r \mid & B(c_1,\nu)(B(c_1,c_2)B(c_1,\nu)-2Q(c_1)B(c_2,\nu))\geq 0 \}\\
P_3:= \{ \nu\in a+\Z^r \mid & \sign(B(c_2,\nu)) =-\sign(B(c_1,c_2)B(c_1,\nu)-2Q(c_1)B(c_2,\nu))\}
\end{split}
\end{equation*}
Note that $B(c_2,\nu)\not=0$ for all $\nu\in a+\Z^r$, which is guaranteed by the fact that $(z,\tau)\in D(c_2)$.

On $P_1$ we use
\begin{equation*}
e^{-2\pi \left( Q(\nu) -\frac{B(c(t),\nu)^2}{2Q(c(t))}\right)y}\leq e^{-2\pi Q(\nu)y},
\end{equation*}
for all $t\in(0,\infty)$. We have seen in the proof of Proposition \ref{lempje} (Case 2) that the series in \eqref{222} converges absolutely. Hence the series
\begin{equation*}
\sum_{\nu\in P_1} e^{-2\pi Q(\nu)y}
\end{equation*}
converges.

On $P_2$ we have
\begin{equation*}
\frac{B(c(t),\nu)^2}{2Q(c(t))}\leq \frac{B(c_1,\nu)^2}{2Q(c_1)}
\end{equation*}
for all $t\in(0,\infty)$, which we get from
\begin{equation*}
B(c_2,\nu)^2+\left(2B(c_2,\nu)B(c_1,\nu)-\frac{B(c_1,c_2)B(c_1,\nu)^2}{Q(c_1)} \right)t \geq 0
\end{equation*}
for all $t\in(0,\infty)$. Hence we find
\begin{equation*}
e^{-2\pi \left( Q(\nu) -\frac{B(c(t),\nu)^2}{2Q(c(t))}\right)y}\leq e^{-2\pi \left( Q(\nu) -\frac{B(c_1,\nu)^2}{2Q(c_1)}\right)y},
\end{equation*}
for all $t\in(0,\infty)$. Using Lemma \ref{sublem1} we see that the series
\begin{equation*}
\sum_{\nu\in P_2}e^{-2\pi \left( Q(\nu) -\frac{B(c_1,\nu)^2}{2Q(c_1)}\right)y}
\end{equation*}
converges.

On $P_3$ we use
\begin{equation}\label{vergl}
\frac{B(c(t),\nu)^2}{2Q(c(t))}\leq \frac{2B(c_2,\nu)}{B(c_1,c_2)^2} \Bigl(B(c_1,c_2)B(c_1,\nu)-Q(c_1)B(c_2,\nu)\Bigr),
\end{equation}
for all $t\in(0,\infty)$, which we get from the inequality
\begin{equation*}
\left( B(c_2,\nu)+\left( -B(c_1,\nu)+ \frac{2Q(c_1)}{B(c_1,c_2)} B(c_2,\nu)\right) t\right)^2 \geq 0.
\end{equation*}
Note that \eqref{vergl} holds also on $P_1$ and $P_2$, but we use it only on $P_3$. Using it we find
\begin{equation*}
e^{-2\pi \left( Q(\nu) -\frac{B(c(t),\nu)^2}{2Q(c(t))}\right)y} \leq e^{-2\pi \widetilde{Q}(\nu)},
\end{equation*}
for all $t\in(0,\infty)$, with
\begin{equation}\label{defq}
\widetilde{Q}(\nu):= Q(\nu)-\frac{2B(c_2,\nu)}{B(c_1,c_2)^2}\Bigl(B(c_1,c_2)B(c_1,\nu)-Q(c_1)B(c_2,\nu)\Bigr).
\end{equation}

Write $\nu = \nu_{c_1} c_1 +\nu_{c_2} c_2 +\nu^{\perp}$, with $\nu^{\perp}$ such that $B(c_1,\nu^\perp)=B(c_2,\nu^\perp)=0$. We see
\begin{equation*}
\begin{split}
B(c_1,\nu)&= 2Q(c_1)\nu_{c_1} + B(c_1,c_2)\nu_{c_2}\\
B(c_2,\nu)&= B(c_1,c_2)\nu_{c_1},
\end{split}
\end{equation*}
so 
\begin{equation*}
\begin{split}
\nu_{c_1} &= \frac{1}{B(c_1,c_2)} B(c_2,\nu)\\
\nu_{c_2} &= \frac{1}{B(c_1,c_2)} B(c_1,\nu)-\frac{2Q(c_1)}{B(c_1,c_2)^2}B(c_2,\nu).
\end{split}
\end{equation*}
Hence
\begin{align*}
Q(\nu)&= Q(\nu^\perp)+Q(c_1)\nu_{c_1}^2 +B(c_1,c_2)\nu_{c_1} \nu_{c_2}\\
&= Q(\nu^\perp)+\frac{B(c_2,\nu)}{B(c_1,c_2)^2}\Bigl(B(c_1,c_2)B(c_1,\nu)-Q(c_1)B(c_2,\nu)\Bigr)\\
\intertext{and}
\widetilde{Q}(\nu)&= Q(\nu^\perp)-\frac{B(c_2,\nu)}{B(c_1,c_2)^2}\Bigl(B(c_1,c_2)B(c_1,\nu)-Q(c_1)B(c_2,\nu)\Bigr).
\end{align*}
The quadratic form $\widetilde{Q}$ has type $(r-1,1)$: $Q$ has type $(1,1)$ on $\left<c_1,c_2\right>_\R$ and type $(r-2,0)$ on $\left<c_1,c_2\right>_\R^\perp$. On $\left<c_1,c_2\right>_\R^\perp$ we have $\widetilde{Q}=Q$ and on $\left<c_1,c_2\right>_\R$ we have $\widetilde{Q}=-Q$. Hence $\widetilde{Q}$ has type $(1,1)$ on $\left<c_1,c_2\right>_\R$ and type $(r-2,0)$ on $\left<c_1,c_2\right>_\R^\perp$.

Set $\tilde{c}_1=\frac{B(c_1,c_2)}{2Q(c_1)}c_1-c_2$ and $\tilde{c}_2=-c_2$ then 
\begin{equation*}
\begin{split}
\widetilde{Q}(\tilde{c}_1)&= \widetilde{Q}\left(\frac{B(c_1,c_2)}{2Q(c_1)}c_1-c_2\right) =-Q\left(\frac{B(c_1,c_2)}{2Q(c_1)}c_1-c_2\right)\\
&=- \frac{B(c_1,c_2)^2}{4Q(c_1)} +\frac{B(c_1,c_2)^2}{2Q(c_1)}=\frac{B(c_1,c_2)^2}{4Q(c_1)}<0\\
\widetilde{Q}(\tilde{c}_2)&= \widetilde{Q}(-c_2)=-Q(c_2)=0\\
\widetilde{B}(\tilde{c}_1,\tilde{c}_2)&= -\widetilde{Q}(\tilde{c}_1-\tilde{c}_2)+\widetilde{Q}(\tilde{c}_1)+\widetilde{Q}(\tilde{c}_2)\\
&=Q\left(\frac{B(c_1,c_2)}{2Q(c_1)}c_1\right)+\frac{B(c_1,c_2)^2}{4Q(c_1)}=\frac{B(c_1,c_2)^2}{2Q(c_1)}<0.
\end{split}
\end{equation*}
If we choose $C_{\widetilde{Q}}$ such that $\tilde{c}_1\in C_{\widetilde{Q}}$ then we see that $\tilde{c}_2\in S_{\widetilde{Q}}$. 

Using \eqref{defq} we see
\begin{equation*}
\begin{split}
\widetilde{B}(x,y)&=\widetilde{Q}(x+y)-\widetilde{Q}(x)-\widetilde{Q}(y)\\
&= B(x,y) -\frac{2}{B(c_1,c_2)} \bigl(B(c_2,x)B(c_1,y)+B(c_1,x)B(c_2,y)\Bigr)\\
&\qquad \qquad +4\frac{Q(c_1)}{B(c_1,c_2)^2}B(c_2,x)B(c_2,y),
\end{split}
\end{equation*}
so
\begin{equation*}
\begin{split}
\widetilde{B}(c_1,\nu)&=-B(c_1,\nu)\\
\widetilde{B}(c_2,\nu)&=-B(c_2,\nu).
\end{split}
\end{equation*}
Since $\tilde{c}_1$ and $\tilde{c}_2$ are linear combinations of $c_1$ and $c_2$, we have
\begin{equation*}
\begin{split}
\widetilde{B}(\tilde{c}_1,\nu)&=-B(\tilde{c}_1,\nu)\\
\widetilde{B}(\tilde{c}_2,\nu)&=-B(\tilde{c}_2,\nu)
\end{split}
\end{equation*}
for all $\nu\in\R^r$. 

We rewrite the set $P_3$:
\begin{equation*}
\begin{split}
P_3&= \{ \nu\in a+\Z^r \mid \sign(B(c_2,\nu)) =-\sign(B(c_1,c_2)B(c_1,\nu)-2Q(c_1)B(c_2,\nu))\}\\
&= \{ \nu\in a+\Z^r \mid  \sign(B(-c_2,\nu)) =-\sign\left(\frac{B(c_1,c_2)}{2Q(c_1)}B(c_1,\nu)-B(c_2,\nu)\right)\}\\
&= \{ \nu\in a+\Z^r \mid  \sign(B(\tilde{c}_2,\nu)) =-\sign(B(\tilde{c}_1,\nu))\}\\
&= \{ \nu\in a+\Z^r \mid  \sign (\widetilde{B}(\tilde{c}_2,\nu)) =-\sign(\widetilde{B}(\tilde{c}_1,\nu)) \}.
\end{split}
\end{equation*}
In the proof of Proposition \ref{lempje} (Case 2, \eqref{222}), we have seen that
\begin{equation*}
\sum_{\nu\in a+\Z^r} \Bigl\{\sign\Bigl(\widetilde{B}(\tilde{c}_1,\nu)\Bigr) - \sign\Bigl(\widetilde{B}(\tilde{c}_2,\nu)\Bigr)\Bigr\}\ e^{2\pi i\widetilde{Q}(\nu)\tau+2\pi i\widetilde{B}(\nu,b)}
\end{equation*}
converges absolutely, i.e.\ the series 
\begin{equation*}
\sum_{\nu\in P_3} e^{-2\pi \widetilde{Q}(\nu)y}
\end{equation*}
converges. 

On all three sets $P_1$, $P_2$ and $P_3$, we have found a suitable majorant, independent of $t\in (0,\infty)$. Combining these results, we see that the series 
\begin{equation*}
\sum_{\nu\in a+\Z^r} \sign\Bigl(B(c(t),\nu)\Bigr) \beta\left( -\frac{B(c(t),\nu)^2}{Q(c(t))}y\right) e^{2\pi iQ(\nu)\tau+2\pi iB(\nu,b)}
\end{equation*}
converges uniformly for $t\in (0,\infty)$. Hence
\begin{equation*}
\begin{split}
\lim_{t\downarrow 0} &\sum_{\nu\in a+\Z^r} \sign\Bigl(B(c(t),\nu)\Bigr) \beta\left( -\frac{B(c(t),\nu)^2}{Q(c(t))}y\right) e^{2\pi iQ(\nu)\tau+2\pi iB(\nu,b)}\\
&=\sum_{\nu\in a+\Z^r} \lim_{t\downarrow 0} \sign\Bigl(B(c(t),\nu)\Bigr) \beta\left( -\frac{B(c(t),\nu)^2}{Q(c(t))}y\right) e^{2\pi iQ(\nu)\tau+2\pi iB(\nu,b)}.
\end{split}
\end{equation*}
We have
\begin{equation*}
\lim_{t\downarrow 0} -\frac{B(c(t),\nu)^2}{Q(c(t))}y =\infty
\end{equation*}
and 
\begin{equation*}
\lim_{x\rightarrow\infty} \beta(x) =0.
\end{equation*}
Hence we get equation \eqref{232}.\\
(6) Since $\rho(a;\tau)$ depends only on $\im(\tau)$, we have $\rho(a;\tau+1)=\rho(a;\tau)$. Hence
\begin{equation}\label{odd}
\theta(z;\tau+1) = \sum_{n\in\Z^r} \rho(n+a;\tau) e^{2\pi iQ(n)} e^{2\pi iQ(n)\tau+2\pi iB(n,z)},
\end{equation}
but since $A$ has integer coefficients we find
\begin{equation*}
e^{2\pi iQ(n)}= e^{\pi i \sum_{l=1}^{r} A_{ll} n_l^2} =e^{\pi i \sum_{l=1}^{r} A_{ll} n_l}= e^{2\pi iB(\frac{1}{2}A^{-1} A^*,n)}.
\end{equation*}
If we put this into \eqref{odd} we get (6).

We first prove (7) for the case $c_1,c_2\in C_Q$. We do this using the Poisson summation formula. The main point -- and the reason for the definition of the function $\rho$ -- is that $a \mapsto \rho(a;\tau) e^{2\pi iQ(a)\tau}$ is more or less its own Fourier transform:

\begin{lemma}\label{lemfour}
We have for all $\alpha\in\R^r$ and $\tau\in\Ha$
\begin{equation*}
\int_{\R^r} \rho(a;\tau)\ e^{2\pi iQ(a)\tau +2\pi iB(a,\alpha)} da= \frac{1}{\sqrt{-\det A}} \frac{i}{(-i\tau)^{r/2}} \rho\Bigl(\alpha;-\frac{1}{\tau}\Bigr)\ e^{-2\pi iQ(\alpha)/\tau}.
\end{equation*}
\end{lemma}
\proof 
The integral converges. This is analogous to the convergence of $\theta$ for case 1: We write $\rho(a;\tau)$ as the sum of the three expressions
\begin{align*}
-&\sign\Bigl(B(c_1,a)\Bigr)\ \beta\left( -\frac{B(c_1,a)^2}{Q(c_1)}y\right),\\
&\sign\Bigl(B(c_2,a)\Bigr)\ \beta\left( -\frac{B(c_2,a)^2}{Q(c_2)}y\right)\\
\intertext{and}
&\sign\Bigl(B(c_1,a)\Bigr) - \sign\Bigl(B(c_2,a)\Bigr).
\end{align*}
We have
\begin{equation*}
\begin{split}
&\left| \sign\Bigl(B(c,a)\Bigr)\ \beta\left( -\frac{B(c,a)^2}{Q(c)}y\right)e^{2\pi iQ(\nu)\tau+2\pi iB(a,b)}\right|\\
&\leq e^{-2\pi \left( Q(a) -\frac{B(c,a)^2}{2Q(c)}\right)y},
\end{split}
\end{equation*}
(see \eqref{223}), with $a \mapsto Q(a) -\frac{B(c,a)^2}{2Q(c)}$ positive definite (see Lemma \ref{sublem1}). We also have 
\begin{equation*}
\left|\Bigl\{\sign\Bigl(B(c_1,a)\Bigr) - \sign\Bigl(B(c_2,a)\Bigr)\Bigr\}\ e^{2\pi iQ(a)\tau+2\pi iB(a,b)}\right| \leq 2 e^{-2\pi Q^+(a)y},
\end{equation*}
(see \eqref{nouja}), with $Q^+$ positive definite (see Lemma \ref{sublem2}).

Using 
\begin{equation*}
\frac{\partial}{\partial \alpha_l} e^{2\pi iQ(a\tau+\alpha)/\tau}=\frac{1}{\tau} \frac{\partial}{\partial a_l} e^{2\pi iQ(a\tau+\alpha)/\tau}
\end{equation*} 
we see that 
\begin{equation}\label{23form1}
\begin{split}
\frac{\partial}{\partial \alpha_l} &\left\{ e^{2\pi iQ(\alpha)/\tau}\int_{\R^r} \rho(a;\tau)\ e^{2\pi iQ(a)\tau +2\pi iB(a,\alpha)} da\right\}\\
&=  \frac{\partial}{\partial \alpha_l} \int_{\R^r} \rho(a;\tau)\ e^{2\pi iQ(a\tau+\alpha)/\tau} da= \int_{\R^r} \rho(a;\tau) \frac{\partial}{\partial \alpha_l} e^{2\pi iQ(a\tau+\alpha)/\tau} da\\
&= \int_{\R^r} \rho(a;\tau) \frac{1}{\tau} \frac{\partial}{\partial a_l} e^{2\pi iQ(a\tau+\alpha)/\tau} da= -\frac{1}{\tau} \int_{\R^r} \frac{\partial\rho}{\partial a_l} (a;\tau)\ e^{2\pi iQ(a\tau+\alpha)/\tau} da,
\end{split}
\end{equation}
where we have used partial integration in the last step.
From the definition of $\rho$ it follows that
\begin{equation}\label{23form2}
\begin{split}
\frac{\partial\rho}{\partial a_l} &(a;\tau) \\
&=\frac{(Ac_1)_l}{\sqrt{-Q(c_1)}}\ y^{1/2}\ E'\Bigl(\frac{B(c_1,a)}{\sqrt{-Q(c_1)}}\ y^{1/2}\Bigr)- \frac{(Ac_2)_l}{\sqrt{-Q(c_2)}}\ y^{1/2}\ E'\Bigl(\frac{B(c_2,a)}{\sqrt{-Q(c_2)}}\ y^{1/2}\Bigr).
\end{split}
\end{equation}
We have (we will use this result for $c=c_1$ and $c=c_2$)
\begin{equation*}
\begin{split}
\int_{\R^r} E'\Bigl(\frac{B(c,a)}{\sqrt{-Q(c)}} y^{1/2}\Bigr)&\ e^{2\pi iQ(a\tau+\alpha)/\tau} da\\
&=2e^{2\pi iQ(\alpha)/\tau} \int_{\R^r} e^{\pi \frac{B(c,a)^2}{Q(c)} y} e^{2\pi iQ(a)\tau + 2\pi i B(a,\alpha)} da.
\end{split}
\end{equation*}

We substitute $a=\left(\begin{smallmatrix}c &C\end{smallmatrix}\right) \left( \begin{smallmatrix} a_c\\a'  \end{smallmatrix} \right)$, with $a_c\in\R$, $a'\in\R^{r-1}$ and $C$ a $r\times (r-1)$-matrix whose columns form a basis for 
\begin{equation*}
\left<c\right>_\R^\perp := \{a\in\R^r  \mid B(c,a)=0\}.
\end{equation*} 
In that way we can split the integral over $\R^r$ in an integral over $\R$ and an integral over $\R^{r-1}$ (Note that $B(c,Ca')=0$, hence $Q(a) = Q(c)a_c^2+\frac{1}{2} \left<a',C^TACa'\right>$ and $B(c,a)=2Q(c)a_c$):
\begin{equation*}
\begin{split}
y^{1/2}& \int_{\R^r}E'\Bigl(\frac{B(c,a)}{\sqrt{-Q(c)}} y^{1/2}\Bigr)\ e^{2\pi iQ(a\tau+\alpha)/\tau} da\\
&=2y^{1/2} e^{2\pi iQ(\alpha)/\tau} \int_{\R\times \R^{r-1}} e^{4\pi Q(c) a_c^2 y+2\pi iQ(c)a_c^2\tau +\pi i \left<a',C^TACa'\right>\tau}\cdot\\
& ~\hspace{5cm} e^{4\pi iQ(c)a_c\alpha_c +2\pi i\left<a',C^TAC\alpha'\right>}\left|\det\begin{pmatrix} c&C\end{pmatrix}\right| da' da_c\\
&= 2\left|\det\begin{pmatrix} c&C\end{pmatrix}\right|y^{1/2} e^{2\pi iQ(\alpha)/\tau}\cdot\\
& \qquad \int_{\R} e^{2\pi iQ(c)a_c^2 \overline{\tau} +4\pi iQ(c)a_c\alpha_c}da_c \cdot \int_{\R^{r-1}} e^{\pi i\left<a',C^TACa'\right>\tau+2\pi i\left<a',C^TAC\alpha'\right>}da',
\end{split}
\end{equation*}
with $\alpha=\left(\begin{smallmatrix}c &C\end{smallmatrix}\right) \left( \begin{smallmatrix} \alpha_c\\\alpha'  \end{smallmatrix} \right)$.

If $\tau\in\Ha$ and $M$ is a positive definite symmetric $n\times n$-matrix, we have the well known result
\begin{equation*}
\int_{\R^n} e^{\pi i\left<a,Ma\right>\tau+2\pi i\left<a,M\alpha\right>}da = \frac{1}{(-i\tau)^{n/2}} \frac{1}{\sqrt{\det M}} e^{-\pi i\left<\alpha,M\alpha\right>/\tau}.
\end{equation*}
(By a change of basis in $\R^n$ one can reduce to the case when $M$ is diagonal).

$Q$ is positive definite on $\left<c\right>_\R^\perp$, so $C^T AC$ is positive definite. Hence we find using the result twice:
\begin{equation*}
\begin{split}
y^{1/2}& \int_{\R^r}E'\Bigl(\frac{B(c,a)}{\sqrt{-Q(c)}} y^{1/2}\Bigr)\ e^{2\pi iQ(a\tau+\alpha)/\tau} da\\
&= 2\left|\det\begin{pmatrix} c&C\end{pmatrix}\right|y^{1/2} e^{2\pi iQ(\alpha)/\tau} \frac{1}{\sqrt{-2iQ(c)\overline{\tau}}}\ e^{-2\pi i Q(c)\alpha_c^2/ \overline{\tau}} \cdot \\
& \hspace{5cm} \frac{1}{(-i\tau)^{(r-1)/2}} \frac{1}{\sqrt{\det C^TAC}}\ e^{-\pi i\left<\alpha',C^TAC\alpha'\right>/\tau}\\
&= \frac{2y^{1/2}}{\sqrt{i\overline{\tau}}} \frac{1}{(-i\tau)^{(r-1)/2}} \frac{\left|\det\begin{pmatrix} c&C\end{pmatrix}\right|}{\sqrt{-2Q(c) \det C^TAC}}\ e^{2\pi i Q(c)\alpha_c^2/ \tau-2\pi i Q(c)\alpha_c^2/ \overline{\tau}}\\
&= \frac{2\sqrt{y'}}{(-i\tau)^{r/2-1}} \frac{\left|\det\begin{pmatrix} c&C\end{pmatrix}\right|}{\sqrt{-2Q(c) \det C^TAC}}\ e^{\pi \frac{B(c,\alpha)^2}{Q(c)}y'}\\
&= \frac{\sqrt{y'}}{(-i\tau)^{r/2-1}} \frac{1}{\sqrt{-\det A}}\ E'\Bigl(\frac{B(c,\alpha)}{\sqrt{-Q(c)}} \sqrt{y'}\Bigr),
\end{split}
\end{equation*}
with $y'=\im\left(-\frac{1}{\tau}\right)$. In the last step we have used 
\begin{equation*}
\left( \det \left(\begin{smallmatrix}c&C\end{smallmatrix}\right)\right)^2 \det A= 2Q(c) \det C^TAC,
\end{equation*}
which follows from
\begin{equation*}
\begin{pmatrix} c&C\end{pmatrix}^T A \begin{pmatrix} c&C\end{pmatrix} = \begin{pmatrix} 2Q(c)&0\\0&C^TAC\end{pmatrix}
\end{equation*}
by taking the determinant. 

We see
\begin{equation*}
\begin{split}
-\frac{1}{\tau} & \int_{\R^r} \frac{(Ac)_l}{\sqrt{-Q(c)}}\ y^{1/2} E'\Bigl(\frac{B(c,a)}{\sqrt{-Q(c_1)}} y^{1/2}\Bigr)\ e^{2\pi iQ(a\tau+\alpha)/\tau} da\\
&= -\frac{1}{\tau} \frac{(Ac)_l}{\sqrt{-Q(c)}} \frac{\sqrt{y'}}{(-i\tau)^{r/2-1}} \frac{1}{\sqrt{-\det A}}\ E'\Bigl(\frac{B(c,\alpha)}{\sqrt{-Q(c)}} \sqrt{y'}\Bigr)\\
&= \frac{\partial}{\partial \alpha_l} \frac{1}{\sqrt{-\det A}} \frac{i}{(-i\tau)^{r/2}}\ E\Bigl(\frac{B(c,\alpha)}{\sqrt{-Q(c)}} \sqrt{y'}\Bigr).
\end{split}
\end{equation*}
Combining this with \eqref{23form1} and \eqref{23form2} we find
\begin{align*}
\frac{\partial}{\partial \alpha_l} &\left\{ e^{2\pi iQ(\alpha)/\tau}\int_{\R^r} \rho(a;\tau) e^{2\pi iQ(a)\tau +2\pi iB(a,\alpha)} da\right\}\\ &=\frac{\partial}{\partial \alpha_l} \frac{1}{\sqrt{-\det A}} \frac{i}{(-i\tau)^{r/2}}\ \rho\Bigl(\alpha;-\frac{1}{\tau}\Bigr).
\end{align*}
So 
\begin{equation*}
e^{2\pi iQ(\alpha)/\tau}\int_{\R^r} \rho(a;\tau)\ e^{2\pi iQ(a)\tau +2\pi iB(a,\alpha)} da - \frac{1}{\sqrt{-\det A}} \frac{i}{(-i\tau)^{r/2}}\ \rho\Bigl(\alpha;-\frac{1}{\tau}\Bigr)
\end{equation*}
is constant as a function of $\alpha$. Since both terms are odd as a function of $\alpha$, that constant is zero. This proves the lemma. \qed
\proof[ of (7)] \emph{Case 1:} $c_1,c_2\in C_Q$.

Using the Poisson summation formula
\begin{equation*}
\sum_{\nu\in\Z^r} f(\nu) = \sum_{\nu\in A^{-1} \Z^r} \tilde{f}(\nu),
\end{equation*}
with $\tilde{f}(\nu)= \int_{\R^r} f(a)\ e^{2\pi iB(\nu,a)} da$, and Lemma \ref{lemfour}, we find that $\theta_{a,b}$  satisfies
\begin{equation}\label{th}
\theta_{a,b}\Bigl(-\frac{1}{\tau}\Bigr) = \frac{i}{\sqrt{-\det A}}\ (-i\tau)^{r/2} e^{2\pi iB(a,b)} \sum_{p\in A^{-1} \Z^r \smod \Z^r} \theta_{b+p,-a}(\tau).
\end{equation}
If we put
\begin{equation*}
\theta_{a,b}(\tau)= e^{2\pi iQ(a)\tau+2\pi iB(a,b)}\ \theta(a\tau+b;\tau)
\end{equation*}
into \eqref{th} (on the left replace $(a,b,\tau)$ by $(a,b,-1/\tau)$, on the right by $(b+p,-a,\tau)$) and multiply both sides by $e^{2\pi iQ(a)/\tau -2\pi iB(a,b)}$, then we find
\begin{equation*}
\theta\Bigl(\frac{b\tau -a}{\tau};-\frac{1}{\tau}\Bigr)= \frac{i}{\sqrt{-\det A}}\ (-i\tau)^{r/2} \sum_{p\in A^{-1}\Z^r /\Z^r} e^{2\pi i Q(b\tau-a+p\tau)/\tau}\theta(b\tau-a+p\tau;\tau),
\end{equation*}
which is the desired result for $z=b\tau-a$.

\noindent{\it Case $2$}: $c_1\in C_Q$ and $c_2\in S_Q$.

We use (5): We have proven the identity for $\theta^{c_1,c(t)}$; if we take $\lim_{t\downarrow 0}$ on both sides we get the desired result.

The other two cases follow using the cocycle conditions given in (1). \qed

\begin{corollary}\label{cor1} The function $\theta_{a,b}$ has the following elliptic and modular transformation properties:
\begin{itemize}
\item[(1)] $\theta_{a+\lambda,b} =\theta_{a,b}$ for all $\lambda\in \Z^r$.
\item[(2)] $\theta_{a,b+\mu} = e^{2\pi iB(a,\mu)}\ \theta_{a,b}$ for all $\mu\in A^{-1}\Z^r$.
\item[(3)] $\theta_{-a,-b} =-\theta_{a,b}$.
\item[(4)] $\theta_{a,b}(\tau+1) = e^{-2\pi iQ(a)-\pi iB(A^{-1}A^*,a)}\ \theta_{a,a+b+\frac{1}{2}A^{-1}A^*}(\tau)$ with $A^*$ the vector of diagonal elements of $A$.
\item[(5)] If $a,b\in R(c_1)\cap R(c_2)$ then
\begin{equation*}
\theta_{a,b}\Bigl(-\frac{1}{\tau}\Bigr) = \frac{i}{\sqrt{-\det A}}\  (-i\tau)^{r/2} e^{2\pi iB(a,b)} \sum_{p\in A^{-1} \Z^r \smod \Z^r} \theta_{b+p,-a}(\tau).
\end{equation*}
\end{itemize}
\end{corollary}

\section{Transformation properties of $\theta$ with respect to $O_A^+(\Z)$}

We consider the group
\begin{equation*}
O_A(\R) := \{C\in \operatorname{GL}_r(\R) \mid C^t A C=A \}.
\end{equation*}
If $C\in O_A(\R)$ and $c\in C_Q \subset \R^r$ then $Q(Cc)=Q(c)$, so $C\cdot C_Q$ is either $C_Q$ or $-C_Q$. We consider only matrices $C$ that leave $C_Q$ invariant, i.e.\ $B(Cc,c)<0$, for all $c\in C_Q$. Set
\begin{equation*}
O_A^+(\R):= \{C\in \operatorname{GL}_r(\R) \mid C^t A C=A,\ B(Cc,c)<0\ \forall c\in C_Q \}.
\end{equation*}
This is a subgroup of $O_A(\R)$ of index 2. 
\begin{definition} Let
\begin{equation*}
O_A^+(\Z):=O_A^+(\R) \cap \operatorname{GL}_r(\Z).
\end{equation*}
\end{definition}
\begin{remark}
From $C^t AC=A$, we find $\det(C)=\pm 1$, so $O_A^+(\Z)$ is the group of elements of $O_A^+(\R)$ that have integer coefficients.
\end{remark}
\begin{remark}
In some cases $O_A^+(\Z)$ is very small. For example if $A= \left(\begin{smallmatrix}1&0\\0&-1\end{smallmatrix}\right)$ then $O_A^+(\Z)$ has only two elements: $\left(\begin{smallmatrix}1&0\\0&1\end{smallmatrix}\right)$ and $\left(\begin{smallmatrix}-1&0\\0&1\end{smallmatrix}\right)$. However, in general $O_A^+(\Z)$ is an infinite group.
\end{remark}
If we consider the theta functions in Definition \ref{deftheta} not only as a function of $z$ and $\tau$, but also as a function of $c_1$ and $c_2$, we get transformation properties with respect to $O_A^+(\Z)$:
\begin{proposition}
Let $C\in O_A^+(\Z)$, $c_1,c_2\in \overline{C}_Q$ and let $(z,\tau)\in D(c_1)\cap D(c_2)$. Let $\theta_A^{c_1,c_2} (z;\tau)$ be as in Definition \ref{deftheta}. Then we have $C\cdot C_Q=C_Q$, $C\cdot S_Q=S_Q$, $(Cz,\tau)\in D(Cc_1)\cap D(Cc_2)$, and 
\begin{equation*}
\theta^{Cc_1,Cc_2} (Cz;\tau)=\theta^{c_1,c_2} (z;\tau).
\end{equation*}
\end{proposition}
\proof $C\cdot C_Q=C_Q$ holds by definition. If $c\in\Z^r$ is primitive, then $Cc$ is also primitive. Hence we find $C\cdot S_Q=S_Q$. We have $Q(Cx)=Q(x)$ and $B(Cx,Cy)=B(x,y)$, for all $x,y\in\R^r$ and $C\Z^r=\Z^r$. We see
\begin{equation*}
B\Bigl(c,\frac{\im(z)}{\im(\tau)}\Bigr)=B\Bigl(Cc,C\frac{\im(z)}{\im(\tau)}\Bigr)=B\Bigl(Cc,\frac{\im(Cz)}{\im(\tau)}\Bigr).
\end{equation*}
Hence if $(z,\tau)\in D(c)$ then $(Cz,\tau)\in D(Cc)$.

If we replace $(c_1,c_2,z,n)$ by $(Cc_1,Cc_2,Cz,Cn)$ in the definition of $\theta$ we get the desired transformation property. \qed

\begin{remark} The $C$ acts on both $c_1$ and $c_2$ at the same time.
\end{remark}

\begin{corollary}\label{cor2}
Let $C\in O_A^+(\Z)$, $c_1,c_2\in \overline{C}_Q$ and let $a\in R(c_1)\cap R(c_2)$. Let $\theta_{a,b}^{c_1,c_2} (\tau)$ be as in Definition \ref{deftheta}. Then
\begin{equation*}
\theta_{Ca,Cb}^{Cc_1,Cc_2} (\tau)=\theta_{a,b}^{c_1,c_2} (\tau).
\end{equation*}
\end{corollary}

\section{Some examples}

\begin{example}
Let $A=\left(\begin{smallmatrix} 1&2\\2&1 \end{smallmatrix}\right)$, $c_1=\left(\begin{smallmatrix} -1\\2 \end{smallmatrix}\right)$, $c_2=\left(\begin{smallmatrix} -2\\1 \end{smallmatrix}\right)$, $e:=\left(\begin{smallmatrix} 1\\1 \end{smallmatrix}\right)$, and $a=b=\frac{1}{6}e$.
Then $B(c_1,c_2)=-6$ and $Q(c_1)=Q(c_2)=-\frac{3}{2}$. If we choose $C_Q$ such that $c_1\in C_Q$ then also $c_2\in C_Q$. Using (4) and (2) of Corollary \ref{cor1}, we see
\begin{equation}\label{transT}
\theta_{\frac{1}{6}e,\frac{1}{6}e}(\tau+1)= e^{-\frac{\pi i}{2}} \theta_{\frac{1}{6}e,\frac{1}{2}e}(\tau)= e^{\frac{\pi i}{6}} \theta_{\frac{1}{6}e,\frac{1}{6}e}(\tau).
\end{equation}
Using (5) of Corollary \ref{cor1}, we see
\begin{equation*}
\theta_{\frac{1}{6}e,\frac{1}{6}e}\Bigl(-\frac{1}{\tau}\Bigr) = \frac{\tau}{\sqrt{3}}\ e^{\frac{\pi i}{3}} \left(\theta_{-\frac{1}{6}e,-\frac{1}{6}e} (\tau)+\theta_{\frac{1}{6}e,-\frac{1}{6}e} (\tau)+\theta_{\frac{1}{2}e,-\frac{1}{6}e} (\tau)\right).
\end{equation*}
Using (3), (2), (3), (1) and (2) of Corollary \ref{cor1}, we see
\begin{equation*}
\begin{split}
\theta_{-\frac{1}{6}e,-\frac{1}{6}e} (\tau) &= -\theta_{\frac{1}{6}e,\frac{1}{6}e} (\tau)\\
\theta_{\frac{1}{6}e,-\frac{1}{6}e} (\tau) &= e^{-\frac{2\pi i}{3}} \theta_{\frac{1}{6}e,\frac{1}{6}e} (\tau)\\
\theta_{\frac{1}{2}e,-\frac{1}{6}e} (\tau) &= -\theta_{-\frac{1}{2}e,\frac{1}{6}e} (\tau)= -\theta_{\frac{1}{2}e,\frac{1}{6}e} (\tau) =-\theta_{\frac{1}{2}e,-\frac{1}{6}e} (\tau).
\end{split}
\end{equation*}
Hence
\begin{equation}\label{transS}
\theta_{\frac{1}{6}e,\frac{1}{6}e}\Bigl(-\frac{1}{\tau}\Bigr) = -i\tau\  \theta_{\frac{1}{6}e,\frac{1}{6}e}(\tau).
\end{equation}
We write $\rho(\nu;\tau)$ as the sum of the three expressions \eqref{al1}, \eqref{al2} and \eqref{al3}. 
We will see that
\begin{align}
\sum_{\nu \in a+\Z^2} \sign\Bigl(B(c_1,\nu)\Bigr) \beta\left( -\frac{B(c_1,\nu)^2}{Q(c_1)}y\right) e^{2\pi i Q(\nu)\tau +2\pi iB(\nu,b)} =0 \label{deel1}\\
\intertext{and}
\sum_{\nu \in a+\Z^2} \sign\Bigl(B(c_2,\nu)\Bigr) \beta\left( -\frac{B(c_2,\nu)^2}{Q(c_2)}y\right) e^{2\pi i Q(\nu)\tau +2\pi iB(\nu,b)} =0. \label{deel2}
\end{align}
To show that \eqref{deel1} holds, consider $C=\left(\begin{smallmatrix}1&0\\-4&-1 \end{smallmatrix}\right)\in O_A^+(\Z)$. If we replace $\nu$ by $C\nu$ in the left hand side of \eqref{deel1} and use $Cc_1=c_1$, we see
\begin{equation*}
\begin{split}
\sum_{\nu\in a+\Z^2} &\sign\Bigl(B(c_1,\nu)\Bigr) \beta\left( -\frac{B(c_1,\nu)^2}{Q(c_1)}y\right) e^{2\pi i Q(\nu)\tau +2\pi iB(\nu,b)}\\
&=\sum_{\nu\in C^{-1}a+\Z^2} \sign\Bigl(B(c_1,\nu)\Bigr) \beta\left( -\frac{B(c_1,\nu)^2}{Q(c_1)}y\right) e^{2\pi i Q(\nu)\tau +2\pi iB(\nu,C^{-1}b)}.
\end{split}
\end{equation*}
Using $C^{-1} a= a-\left(\begin{smallmatrix} 0\\1\end{smallmatrix}\right)$, $C^{-1} b= b-\left(\begin{smallmatrix} 0\\1\end{smallmatrix}\right)$ and $B\left(\nu,\left(\begin{smallmatrix} 0\\1\end{smallmatrix}\right)\right)=2\nu_1+\nu_2\equiv \frac{1}{2} \smod{1}$ for $\nu\in a+\Z^2$, we see
\begin{equation*}
\begin{split}
\sum_{\nu\in a+\Z^2} &\sign\Bigl(B(c_1,\nu)\Bigr) \beta\left( -\frac{B(c_1,\nu)^2}{Q(c_1)}y\right) e^{2\pi i Q(\nu)\tau +2\pi iB(\nu,b)}\\
&=-\sum_{\nu\in a+\Z^2} \sign\Bigl(B(c_1,\nu)\Bigr) \beta\left( -\frac{B(c_1,\nu)^2}{Q(c_1)}y\right) e^{2\pi i Q(\nu)\tau +2\pi iB(\nu,b)}.
\end{split}
\end{equation*}
Hence we get \eqref{deel1}. The proof of \eqref{deel2} is similar. Here we have to use $C=\left(\begin{smallmatrix}-1&-4\\0&1 \end{smallmatrix}\right)\in O_A^+(\Z)$.

Using \eqref{deel1} and \eqref{deel2} we see
\begin{equation}\label{ber}
\begin{split}
\theta_{a,b}(\tau) &= \sum_{\nu\in a+\Z^2} \Bigl\{\sign\Bigl(B(c_1,\nu)\Bigr)-\sign\Bigl(B(c_2,\nu)\Bigr)\Bigr\}\ e^{2\pi i Q(\nu)\tau +2\pi iB(\nu,b)}\\
&= \sum_{\nu\in\left(\frac{1}{6}+\Z\right)^2} \Bigl\{\sign(\nu_1) +\sign(\nu_2)\Bigr\}\ e^{2\pi i(\frac{1}{2}\nu_1^2 +2\nu_1 \nu_2 +\frac{1}{2} \nu_2^2)\tau +\pi i(\nu_1+\nu_2)}\\
&= 2e^{\frac{\pi i}{3}} q^{\frac{1}{12}} \Bigl( \sum_{n,m\geq 0} - \sum_{n,m<0} \Bigr) (-1)^{n+m} q^{\frac{1}{2}n^2+2nm+\frac{1}{2}m^2 +\frac{1}{2}n+\frac{1}{2}m},
\end{split}
\end{equation}
where we have substituted $\nu_1=\frac{1}{6}+n$ and $\nu_2=\frac{1}{6}+m$ in the last step. 

From \eqref{ber} together with \eqref{transT} and \eqref{transS} we see that $\theta_{a,b}$ is a holomorphic modular form of weight 1, with the same transformation properties as $\eta^2$. Hence $\theta_{a,b}$ is a multiple of $\eta^2$ (their quotient is a holomorphic function on the compact Riemann surface $\Ha/\operatorname{SL}_2(\Z) \cup \{\infty \}$). By comparing the first Fourier coefficients we find
\begin{equation*}
\theta_{\frac{1}{6}e,\frac{1}{6}e} =2e^{\frac{\pi i}{3}} \eta^2,
\end{equation*}
or equivalently
\begin{equation*}
\Bigl( \sum_{n,m\geq 0} - \sum_{n,m<0} \Bigr) (-1)^{n+m} q^{\frac{1}{2}n^2+2nm+\frac{1}{2}m^2 +\frac{1}{2}n+\frac{1}{2}m} = (q)_\infty^2.
\end{equation*}
\end{example}
\begin{example} This example is similar to the previous one, so some of the details are omitted.

Let $A=\left(\begin{smallmatrix} 1&0\\0&-3 \end{smallmatrix}\right)$, $c_1=\left(\begin{smallmatrix} -3\\2 \end{smallmatrix}\right)$, $c_2=\left(\begin{smallmatrix} 3\\2 \end{smallmatrix}\right)$ and $a=b=\frac{1}{6}\left(\begin{smallmatrix} 3\\-1 \end{smallmatrix}\right)$.
Then $B(c_1,c_2)=-21$ and $Q(c_1)=Q(c_2)=-\frac{3}{2}$. Using Corollary \ref{cor1}, we see
\begin{align*}
\theta_{a,b}(\tau+1)&= e^{\frac{\pi i}{6}} \theta_{a,b}(\tau)\\
\intertext{and}
\theta_{a,b}\Bigl(-\frac{1}{\tau}\Bigr) &= -i\tau\ \theta_{a,b}(\tau).
\end{align*}
We write $\rho(\nu;\tau)$ as the sum of the three expressions \eqref{al1}, \eqref{al2} and \eqref{al3}.
We have
\begin{align*}
\sum_{\nu \in a+\Z^2} \sign\Bigl(B(c_1,\nu)\Bigr) \beta\left( -\frac{B(c_1,\nu)^2}{Q(c_1)}y\right) e^{2\pi i Q(\nu)\tau +2\pi iB(\nu,b)} =0 \\
\intertext{and}
\sum_{\nu \in a+\Z^2} \sign\Bigl(B(c_2,\nu)\Bigr) \beta\left( -\frac{B(c_2,\nu)^2}{Q(c_2)}y\right) e^{2\pi i Q(\nu)\tau +2\pi iB(\nu,b)} =0. 
\end{align*}
To get the first equation we use $C=\left(\begin{smallmatrix}-7&-12\\4&7 \end{smallmatrix}\right)\in O_A^+(\Z)$ ($Cc_1=c_1$, $C^{-1}a=a+\left(\begin{smallmatrix} -2\\1\end{smallmatrix}\right)$ and $C^{-1}b=b+\left(\begin{smallmatrix} -2\\1\end{smallmatrix}\right)$). To get the second equation we use $C=\left(\begin{smallmatrix}-7&12\\-4&7 \end{smallmatrix}\right)\in O_A^+(\Z)$ ($Cc_2=c_2$, $C^{-1}a=a-\left(\begin{smallmatrix} 4\\3\end{smallmatrix}\right)$ and $C^{-1}b=b-\left(\begin{smallmatrix} 4\\3\end{smallmatrix}\right)$). 
Hence we see
\begin{equation*}
\begin{split}
\theta_{a,b}(\tau) &= \sum_{\nu\in a+\Z^2} \Bigl\{\sign\Bigl(B(c_1,\nu)\Bigr)-\sign\Bigl(B(c_2,\nu)\Bigr)\Bigr\}\ e^{2\pi i Q(\nu)\tau +2\pi iB(\nu,b)}\\
&= -\sum_{\nu\in\frac{1}{6}\left(\begin{smallmatrix} 3\\-1\end{smallmatrix}\right)+\Z^2} \Bigl\{\sign(\nu_1+2\nu_2) +\sign(\nu_1-2\nu_2)\Bigr\}\ e^{2\pi i(\frac{1}{2}\nu_1^2 -\frac{3}{2} \nu_2^2)\tau +\pi i(\nu_1+\nu_2)}\\
&= -2e^{\frac{\pi i}{3}} q^{\frac{1}{12}} \Bigl( \sum_{n+2m,n-2m\geq 0} - \sum_{n+2m,n-2m<0} \Bigr)\ (-1)^{n+m} q^{\frac{1}{2}n^2-\frac{3}{2}m^2 +\frac{1}{2}n+\frac{1}{2}m},
\end{split}
\end{equation*}
where we have substituted $\nu_1=\frac{1}{2}+n$ and $\nu_2=-\frac{1}{6}+m$ in the last step. Replacing $n$ by $-n-1$, we see
\begin{equation*}
\begin{split}
\sum_{n+2m,n-2m<0}&(-1)^{n+m} q^{\frac{1}{2}n^2-\frac{3}{2}m^2 +\frac{1}{2}n+\frac{1}{2}m}\\
&= \sum_{n+2m,n-2m\geq 0} (-1)^{n+m} q^{\frac{1}{2}n^2-\frac{3}{2}m^2 +\frac{1}{2}n+\frac{1}{2}m},
\end{split}
\end{equation*}
so
\begin{equation*}
\theta_{a,b}(\tau)= -4e^{\frac{\pi i}{3}} q^{\frac{1}{12}} \sum_{n\geq 2|m|} (-1)^{n+m} q^{\frac{1}{2}n^2-\frac{3}{2}m^2 +\frac{1}{2}n+\frac{1}{2}m}.
\end{equation*}
We see that $\theta_{a,b}$ is a holomorphic modular form of weight 1, with the same transformation properties as $\eta^2$. Hence $\theta_{a,b}$ is a multiple of $\eta^2$. By comparing the first Fourier coefficients we find
\begin{equation*}
\theta_{a,b} =-4e^{\frac{\pi i}{3}} \eta^2,
\end{equation*}
or equivalently
\begin{equation*}
\sum_{n\geq 2|m|} (-1)^{n+m} q^{\frac{1}{2}n^2-\frac{3}{2}m^2 +\frac{1}{2}n+\frac{1}{2}m} = (q)_\infty^2.
\end{equation*}
This last equation is proven in \cite[pp.\ 451]{andrews}, using different techniques. In that article several similar results are proven. The modular transformation properties of the functions involved can be found using the same method as in the examples presented here. 
\end{example}

These examples are very special: in general, $\theta_{a,b}$ is not a holomorphic function. However, for the special values of $c_1$, $c_2$, $a$ and $b$ given here, $\theta_{a,b}$ is holomorphic.

In \cite{polishchuk} a theorem about the modularity of a certain family of $q$-series associated with indefinite binary quadratic forms is given. This result may also be found using the same method as in the examples presented here.

In the next two chapters, we will see some other examples. In these examples, the $\theta$-functions are not holomorphic.

\clearemptydoublepage
\chapter{Fourier Coefficients of Meromorphic Jacobi Forms}
\section{Introduction}

In this chapter we consider functions $\phi:\C\times\Ha \longrightarrow \C$ that satisfy
\begin{equation*}
\phi(z+\lambda \tau+\mu;\tau)= e^{-2\pi im(\lambda^2 \tau+2\lambda z)} \phi(z;\tau) \qquad \forall \lambda,\mu\in\Z \tag{E}
\end{equation*}
and
\begin{equation*}
\phi \left( \frac{z}{c\tau+d};\frac{a\tau+b}{c\tau+d}\right) = (c\tau+d)^k e^{2\pi imcz^2/(c\tau+d)} \phi(z;\tau) \qquad \forall \left(\begin{smallmatrix}a&b\\c&d\end{smallmatrix}\right) \in\operatorname{SL}_2(\Z), \tag{M}
\end{equation*}
with $k\in\Z$ and $m\in\Z_{>0}$.
The first equation gives the transformation law with respect to $z\mapsto z+\lambda\tau+\mu$ and will be denoted by (E), for elliptic. The second equation gives the transformation law with respect to $\operatorname{SL}_2(\Z)$ and will be denoted by (M), for modular. Jacobi forms of weight $k$ and index $m$ satisfy both (E) and (M).

It is a classical result, see \cite[pp.\ 57--59]{eichler}, that the space of Jacobi forms of weight $k$ and index $m$ is isomorphic to a certain space of (vector-valued) modular forms of weight $k-\frac{1}{2}$ in one variable:

\begin{theorem}\label{them311} If $\phi$ is holomorphic as a function of $z$ and satisfies (E), we have
\begin{equation}\label{them311form0}
\phi(z;\tau) = \sum_{l\smod{2m}} h_l(\tau) \theta_{m,l} (z;\tau),
\end{equation}
with Fourier coefficients
\begin{equation*}
h_l(\tau) = e^{-\pi i l^2\tau/2m} \int_p^{p+1} \phi(z;\tau) e^{-2\pi ilz} dz \qquad p\in\C
\end{equation*}
and
\begin{equation}\label{them311form1}
\theta_{m,l} (z;\tau) = \underset{\lambda\equiv l\smod{2m}}{\sum_{\lambda\in \Z}} e^{\pi i\lambda^2 \tau/2m +2\pi i\lambda z}. 
\end{equation}
If $\phi$ also satisfies the transformation (M), then we have for each $l$:
\begin{equation}\label{them311form2}
h_l(\tau+1) = e^{-\pi il^2/2m} h_l(\tau)
\end{equation}
and 
\begin{equation}\label{them311form3}
h_l \Bigl(-\frac{1}{\tau}\Bigr) = \frac{\tau^k}{\sqrt{-i\tau}} \frac{1}{\sqrt{2m}} \sum_{\nu\smod{2m}} e^{\pi i l\nu/m}h_\nu(\tau).
\end{equation}
\end{theorem}
\proof{}If we take $\lambda=0$ and $\mu=1$ in (E) we see that $\phi$ is 1-periodic. Hence we can write
\begin{equation}\label{theta1}
\phi(z;\tau) = \sum_{r\in\Z} h_r(\tau) e^{\pi i r^2\tau/2m +2\pi irz},
\end{equation}
with
\begin{equation*}
h_r(\tau) = e^{-\pi i r^2\tau/2m} \int_p^{p+1} \phi(z;\tau) e^{-2\pi irz} dz \qquad p\in\C.
\end{equation*}
(The extra factor $e^{\pi i r^2\tau/2m}$ in the Fourier coefficients is for convenience).

If we use (E) with $\lambda=1$ and $\mu=0$ we see that $h_{r+2m}=h_r$. Hence $h_r$ depends only on $r\smod{2m}$. Putting this into \eqref{theta1} gives
\begin{equation*}
\begin{split}
\phi(z;\tau) &= \sum_{r\in\Z} h_r(\tau) e^{\pi i r^2\tau/2m +2\pi irz} = \sum_{l\smod{2m}} \underset{\lambda\equiv l\smod{2m}}{\sum_{\lambda\in \Z}} h_{\lambda} (\tau) e^{\pi i \lambda ^2\tau/2m +2\pi i\lambda z}\\
&= \sum_{l\smod{2m}} h_l(\tau) \theta_{m,l}(z;\tau),
\end{split}
\end{equation*}
with $\theta_{m,l}$ as in \eqref{them311form1}. 

If $\phi$ also satisfies the transformation (M), then we get the transformation properties \eqref{them311form2} and \eqref{them311form3} of $h_l$ from the transformation properties of $\theta_{m,l}$ and the decomposition given in \eqref{them311form0}; see \cite[pp. 58--59]{eichler} for details. \qed

The $h_l$ are more or less the Fourier coefficients of $\phi$, if we consider $\phi$ as a function of $z$. These Fourier coefficients form a vector-valued modular form. 

In  \cite{coeff} Andrews gives most of the fifth order mock theta functions as Fourier coefficients of meromorphic Jacobi forms (i.e.\ meromorphic as a function of $z$), namely certain quotients of ordinary Jacobi theta-series. In this chapter (see Theorem \ref{themv}), we generalize Theorem \ref{them311} to include meromorphic Jacobi forms. We give the result only for Jacobi forms on the full Jacobi group (i.e.\ satisfying (E) and (M) without any congruence restrictions on $(\lambda,\mu)$ or $\left(\begin{smallmatrix}a&b\\c&d\end{smallmatrix}\right)$), but it could certainly be generalized to congruence subgroups (and vector-valued Jacobi forms) and could then be combined with Andrews's identities to obtain information about the modular properties of the fifth order mock theta functions. We will not carry this out, since the same results will be obtained in Chapter 4 using instead the results on indefinite $\theta$-functions from Chapter 2. 

\section{A building block}

In this section we define functions $f_u:\C\times\Ha \longrightarrow \C$ and $\tilde{f}_u:\C\times\Ha \longrightarrow \C$, which will be used in the next section as building blocks for meromorphic Jacobi forms. 
\begin{definition}
Let $u\in\C$ and $m\in\Z_{>0}$. Define $f_u:\C\times\Ha \longrightarrow \C$ by
\begin{equation*}
f_u(z;\tau) = f_u^{(m)} (z;\tau):= \sum_{\lambda\in\Z} \frac{e^{2\pi im\lambda^2 \tau +4\pi im\lambda z}}{1-e^{2\pi i\lambda \tau +2\pi i(z-u)}}.
\end{equation*}
\end{definition}
Note the similarity of this sum with the Lerch sums studied in Chapter 1. The function $f^{(1/2)}$ is the sum studied in Section 1.2. The following result is the analogue of Proposition \ref{prop3} and Proposition \ref{prop4}.
\begin{proposition}\label{prop121}
We have
\begin{description}
\item[(1)] $f_u$ satisfies (E),
\item[(2)] $z\mapsto f_u(z;\tau)$ is a meromorphic function, with simple poles in $u+\Z\tau+\Z$, and residue $-\frac{1}{2\pi i}$ in $z=u$,
\item[(3)] $f_{u+1}(z;\tau)=f_u(z;\tau)$ 
\item[(4)] $\displaystyle{f_u(z;\tau)-e^{-2\pi im\tau-4\pi imu} f_{u+\tau} (z;\tau)= \sum_{l=0}^{2m-1} e^{-\pi il^2\tau/2m -2\pi ilu} \theta_{m,l}(z;\tau)}$,
\item[(5)] $f_u(z;\tau +1)=f_u(z;\tau)$,
\item[(6)] $\displaystyle{f_u(z;\tau) -\frac{1}{\tau} e^{2\pi im (u^2-z^2)/\tau} f_{\frac{u}{\tau}} \Bigl( \frac{z}{\tau}; -\frac{1}{\tau} \Bigr) = \sum_{l=0}^{2m-1} h_l (u;\tau) \theta_{m,l}(z;\tau)}$, with 
\begin{equation*}
h_l(u;\tau) = i e^{-\pi i l^2 \tau/2m -2\pi ilu} \int_L \frac{e^{2\pi im\tau x^2 -2\pi(2mu+l\tau)x}}{1-e^{2\pi x}} dx,
\end{equation*} 
where $L=\R-it$ with $0<t<1$. This path can be deformed into the real axis indented by the lower half of a small circle with the origin as its centre.
\end{description}
\end{proposition}
\proof{}We see immediately that the series converges absolutely, unless $z=u-\lambda\tau+\mu$ for some $\lambda,\mu\in\Z$, in which case one term in the sum becomes infinite. Hence $z\mapsto f_u(z;\tau)$ is meromorphic, with simple poles only in the points $z=u-\lambda\tau+\mu\ (\lambda,\mu\in\Z)$.
(1) It is easy to see that $f_u(z+\mu)=f_u(z)$ for all $\mu\in\Z$. Also
\begin{equation*} 
\begin{split}
f_u(z+\mu\tau;\tau)&=  \sum_{\lambda\in\Z} \frac{e^{2\pi im\lambda^2 \tau +4\pi im\lambda z +4\pi im\mu\lambda \tau}}{1-e^{2\pi i\lambda \tau +2\pi i(z-u)+2\pi i\mu\tau}}\\
&= e^{-2\pi im\mu^2\tau -4\pi im\mu z} \sum_{\lambda\in\Z} \frac{e^{2\pi im(\lambda+\mu)^2 \tau +4\pi im(\lambda+\mu) z}}{1-e^{2\pi i(\lambda+\mu) \tau +2\pi i(z-u)}}\\
&= e^{-2\pi im\mu^2\tau -4\pi im\mu z} f_u(z;\tau).
\end{split}
\end{equation*}
(2) We have already seen that $z\mapsto f_u(z;\tau)$ is a meromorphic function, with simple poles in $u+\Z\tau+\Z$. The pole in $z=u$ comes from the term $\lambda=0$. We see
\begin{equation*}
\lim_{z\rightarrow u} (z-u) f_u(z;\tau) = \lim_{z\rightarrow u} \frac{(z-u)}{1-e^{2\pi i(z-u)}} = -\frac{1}{2\pi i}.
\end{equation*}
(3) Trivial.\\
(4) If we replace $\lambda$ by $\lambda+1$ in the definition we find
\begin{equation*}
e^{-2\pi im\tau-4\pi imu}f_{u+\tau}(z;\tau)= \sum_{\lambda\in\Z} \frac{e^{2\pi im\lambda^2 \tau +4\pi im\lambda z+4\pi im\lambda\tau+4\pi im(z-u)}}{1-e^{2\pi i\lambda \tau +2\pi i(z-u)}}.
\end{equation*}
Hence
\begin{equation*}
\begin{split}
f_u(z;\tau)&-e^{-2\pi im\tau-4\pi imu} f_{u+\tau} (z;\tau) = \sum_{\lambda\in\Z} e^{2\pi im\lambda^2 \tau +4\pi im\lambda z}\frac{1-e^{4\pi im\lambda\tau+4\pi im(z-u)}}{1-e^{2\pi i\lambda \tau +2\pi i(z-u)}}\\
&= \sum_{\lambda\in\Z} e^{2\pi im\lambda^2 \tau +4\pi im\lambda z} \sum_{l=0}^{2m-1} e^{2\pi il\lambda\tau +2\pi il(z-u)}\\
&= \sum_{l=0}^{2m-1} e^{-\pi il^2\tau/2m -2\pi ilu} \theta_{m,l}(z;\tau).
\end{split}
\end{equation*}
In the last step we have changed the order of summation and substituted $\mu=2m\lambda+l$.\\
(5) Trivial.\\
(6) If $z\mapsto f(z;\tau)$ satisfies (E), then so does $z\mapsto e^{-2\pi im z^2/\tau}f\left( \frac{z}{\tau};-\frac{1}{\tau}\right)$. The function $z\mapsto e^{-2\pi im z^2/\tau} f_{\frac{u}{\tau}}\left( \frac{z}{\tau};-\frac{1}{\tau}\right)$ is meromorphic, with simple poles in $u+\Z\tau+\Z$, and residue $-\frac{\tau}{2\pi i} e^{-2\pi im u^2/\tau}$ in $z=u$. So
\begin{equation*}
z\mapsto f_u(z;\tau) -\frac{1}{\tau} e^{2\pi im (u^2-z^2)/\tau} f_{\frac{u}{\tau}}\Bigl( \frac{z}{\tau};-\frac{1}{\tau}\Bigr)
\end{equation*}
is a holomorphic function, which satisfies (E). Theorem \ref{them311} shows that there are $h_l$ such that
\begin{equation*}
f_u(z;\tau) -\frac{1}{\tau} e^{2\pi im (u^2-z^2)/\tau} f_{\frac{u}{\tau}}\Bigl( \frac{z}{\tau};-\frac{1}{\tau}\Bigr) = \sum_{l\smod{2m}} h_l(\tau) \theta_{m,l} (z;\tau).
\end{equation*}
If we restrict $z$ by $0<\im(u-z)<\im(\tau)$ and expand $\Bigl(1-e^{2\pi i\lambda \tau +2\pi i(z-u)}\Bigr)^{-1}$ into a geometric series, we see
\begin{equation*}
\begin{split}
f_u(z;\tau) &= \sum_{\lambda\in\Z} e^{2\pi im\lambda^2 \tau +4\pi im\lambda z} \sign \Bigl(\lambda -\frac{1}{2}\Bigr) \underset{\sign(\lambda-\frac{1}{2})=\sign(\mu+\frac{1}{2})}{\sum_{\mu\in\Z}} e^{2\pi i\lambda \mu\tau +2\pi i(z-u)\mu}\\
&= \underset{\sign(\lambda-\frac{1}{2})=\sign(\mu+\frac{1}{2})}{\sum_{\lambda,\mu\in\Z}} \sign \Bigl(\lambda -\frac{1}{2}\Bigr)\ e^{2\pi im\lambda^2 \tau +2\pi i\lambda \mu\tau -2\pi i\mu u+2\pi i(2m\lambda+\mu)z}.
\end{split}
\end{equation*}
But if $\sign(\lambda-\frac{1}{2})=\sign(\mu+\frac{1}{2})$ then $2m\lambda+\mu$ is either $\geq 2m$ or $<0$, so
\begin{equation*}
\int_p^{p+1} f_u(z;\tau) e^{-2\pi ilz} dz =0, 
\end{equation*}
for $0\leq l\leq 2m-1$, for any $p$ satisfying $0<\im(u-p)<\im(\tau)$. Hence for $0\leq l\leq 2m-1$ we find
\begin{equation*}
\begin{split}
h_l(\tau) &= -\frac{1}{\tau} e^{2\pi imu^2/\tau -\pi i l^2\tau/2m} \int_p^{p+1} e^{-2\pi imz^2/\tau -2\pi ilz} f_{\frac{u}{\tau}} \left( \frac{z}{\tau};-\frac{1}{\tau}\right) dz\\
&= -\frac{1}{\tau} e^{2\pi imu^2/\tau -\pi i l^2\tau/2m} \int_p^{p+1} \sum_{\lambda\in\Z} \frac{e^{-2\pi im(z-\lambda)^2/\tau -2\pi il(z-\lambda)}}{1-e^{2\pi i (z-\lambda -u)/\tau}} dz\\
&= -\frac{1}{\tau} e^{2\pi imu^2/\tau -\pi i l^2\tau/2m} \int_{p+\R} \frac{e^{-2\pi imz^2/\tau -2\pi ilz}}{1-e^{2\pi i (z-u)/\tau}} dz.
\end{split}
\end{equation*}
If we substitute $x= i(z-u)/\tau$ in this last integral and use Cauchy's theorem, we get the desired result. \qed

We see that the transformation law with respect to $\operatorname{SL}_2(\Z)$ for $f_u$ is rather complicated. However, if we modify the definition a little (in Theorem \ref{them1} we did something similar), we get a function which is no longer holomorphic as a function of $\tau$, but has simpler transformation properties with respect to $\operatorname{SL}_2(\Z)$.

\begin{definition}\label{defft}
Let $u\in\C$. Define $\tilde{f}_u:\C\times\Ha \longrightarrow \C$ by
\begin{equation*}
\tilde{f}_u (z;\tau) =f_u (z;\tau) - \frac{1}{2} \sum_{l\smod{2m}} R_{m,l}(u;\tau) \theta_{m,l} (z;\tau),
\end{equation*}
with
\begin{equation*}
R_{m,l}(u;\tau) = \underset{\lambda\equiv l\smod{2m}}{\sum_{\lambda\in \Z}} \Bigl\{\sign\Bigl(\lambda+\frac{1}{2}\Bigr) -E\Bigl((\lambda+2m\im (u) /y)\sqrt{y/m}\Bigr)\Bigr\}\ e^{-\pi i\lambda^2 \tau/2m-2\pi i\lambda u},
\end{equation*}
$y=\im(\tau)$ and $E$ as in Definition \ref{defE}.
\end{definition}
The function $R_{m,l}$ is the analogue of the function $R$ defined in Lemma \ref{defR}. We will not show that the series defining $R_{m,l}$ converges, since this is similar to the convergence of $R$, proven in Lemma \ref{defR}.
\begin{proposition}\label{prop222}
We have
\begin{description}
\item[(1)] $\tilde{f}_u(z;\tau)$ transforms like a 2-variable Jacobi form of weight 1 and index $\left(\begin{smallmatrix}2m&0\\0&-2m\end{smallmatrix}\right)$ with respect to $(z,u,\tau)\in\C^2\times\Ha$, i.e.
\begin{description}
\item[(a)] $\tilde{f}_u$ satisfies (E),
\item[(b)] $\tilde{f}_{u+\lambda\tau+\mu} (z;\tau)= e^{2\pi im (\lambda^2 \tau+2\lambda u)} \tilde{f}_u (z;\tau)$ for all $\lambda,\mu\in\Z$,
\item[(c)] 
\begin{equation*}
\tilde{f}_{\frac{u}{c\tau+d}} \left( \frac{z}{c\tau+d};\frac{a\tau+b}{c\tau+d} \right) = (c\tau+d) e^{2\pi imc (z^2-u^2)/(c\tau+d)} \tilde{f}_u(z;\tau)
\end{equation*}
for all $\left(\begin{smallmatrix} a&b\\c&d\end{smallmatrix}\right)\in\operatorname{SL}_2(\Z)$.
\end{description}
\item[(2)] for fixed $u$, $z\mapsto \tilde{f}_u(z;\tau)$ is a meromorphic function, with simple poles in $u+\Z\tau+\Z$, and no other poles, and residue $-\frac{1}{2\pi i}$ in $z=u$,
\item[(3)] $\tilde{f}$ can be seen as a indefinite $\theta$-series (see Definition \ref{deftheta}), namely
\begin{equation*}
\tilde{f}_u (z;\tau) = \frac{1}{2} \theta_A^{c_1,c_2}\left(\left(\begin{smallmatrix}z-u\\2mu\end{smallmatrix}\right);\tau\right),
\end{equation*}
with $A=\left(\begin{smallmatrix}2m&1\\1&0\end{smallmatrix}\right)$,  $c_1=\left(\begin{smallmatrix}0\\1\end{smallmatrix}\right)\in S_Q$ and $c_2=\left(\begin{smallmatrix}-1\\2m\end{smallmatrix}\right)\in C_Q$.
\end{description}
\end{proposition}
\proof{} (2) This follows directly from (2) of Proposition \ref{prop121} and the fact that $\theta_{m,l}$ is holomorphic.\\
(3) Using the geometric series expansion we see
\begin{equation*}
\frac{1}{1-e^{2\pi i\lambda \tau+2\pi i(z-u)}}=\frac{1}{2} \sum_{\mu\in\Z} \Bigl\{ \sign \Bigl( \lambda+\im(z-u)/y\Bigr) +\sign\Bigl(\mu+\frac{1}{2}\Bigr) \Bigr\}\ e^{2\pi i\lambda \mu \tau +2\pi i(z-u)\mu}.
\end{equation*}
Hence
\begin{equation*}
\begin{split}
f_u&(z;\tau)\\
&= \sum_{\lambda,\mu\in\Z} \frac{1}{2} \Bigl\{\sign \Bigl(\lambda +\im(z-u)/y\Bigr)+\sign \Bigl(\mu+\frac{1}{2}\Bigr)\Bigr\}\ e^{2\pi i(m\lambda^2 +\lambda \mu)\tau +2\pi i(2mz\lambda +(z-u)\mu)}\\
&= \sum_{n\in\Z^2} \frac{1}{2} \Bigl\{\sign \Bigl(n_1 +\im(z-u)/y\Bigr)+\sign \Bigl(n_2+\frac{1}{2}\Bigr)\Bigr\}\ e^{2\pi iQ(n)\tau +2\pi iB\left(n,\left(\begin{smallmatrix}z-u\\2mu\end{smallmatrix}\right)\right)}.
\end{split}
\end{equation*}
We also have 
\begin{equation*}
\begin{split}
&\sum_{l\smod{2m}} R_{m,l}(u;\tau) \theta_{m,l} (z;\tau)\\
&= \underset{\lambda\equiv \mu \smod{2m}}{\sum_{\lambda,\mu\in\Z}} \Bigl\{\sign\Bigl(\lambda+\frac{1}{2}\Bigr) -E\Bigl((\lambda+2m\im u /y)\sqrt{y/m}\Bigr)\Bigr\}\ e^{\pi i(\mu^2-\lambda^2) \tau/2m+2\pi i(\mu z-\lambda u)}.
\end{split}
\end{equation*}
If we substitute $\lambda=n_2$ and $\mu=2mn_1+n_2$, we find
\begin{equation*}
\begin{split}
&\sum_{l\smod{2m}} R_{m,l}(u;\tau) \theta_{m,l} (z;\tau)\\
&= \sum_{n\in\Z^2} \Bigl\{\sign\Bigl(n_2+\frac{1}{2}\Bigr) -E\Bigl((n_2+2m\im u /y)\sqrt{y/m}\Bigr)\Bigr\}\ e^{2\pi iQ(n) \tau+2\pi iB\left(n,\left(\begin{smallmatrix}z-u\\2mu\end{smallmatrix}\right)\right)}.
\end{split}
\end{equation*}
Hence
\begin{equation*}
\begin{split}
\tilde{f}_u (z;\tau)&=f_u (z;\tau) - \frac{1}{2} \sum_{l\smod{2m}} R_{m,l}(u;\tau) \theta_{m,l} (z;\tau)\\
&=\sum_{n\in\Z^2} \frac{1}{2}\left\{\sign\Bigl(n_1+\im(z-u)/y\Bigr) +E\Bigl((n_2+2m\im u /y\Bigr)\sqrt{y/m})\right\} \cdot \\
&\hspace{5cm}\cdot e^{2\pi iQ(n) \tau+2\pi iB\left(n,\left(\begin{smallmatrix}z-u\\2mu\end{smallmatrix}\right)\right)}\\
&= \frac{1}{2} \theta_A^{c_1,c_2}\Bigl(\left(\begin{smallmatrix}z-u\\2mu\end{smallmatrix}\right);\tau\Bigr),
\end{split}
\end{equation*}
which proves (3).\\
(1a) and (1b) follow directly from (3) and the transformation properties of $\theta_A^{c_1,c_2}$ given in (2) of Proposition \ref{thet}.\\
(1c) Using (3) and the transformation properties of $\theta_A^{c_1,c_2}$ given in (6) and (7) of Proposition \ref{thet} we see
\begin{equation}\label{jacft}
\tilde{f}_u (z;\tau+1) = \tilde{f}_u (z;\tau)
\end{equation}
and
\begin{equation}\label{jacfs}
\tilde{f}_{\frac{u}{\tau}} \Bigl( \frac{z}{\tau};-\frac{1}{\tau} \Bigr) = \tau e^{2\pi im (z^2-u^2)/\tau} \tilde{f}_u(z;\tau).
\end{equation}
Combining these results we get (1c). \qed

\begin{remark}
We do not need (3) to prove (1): We could also prove \eqref{jacfs} using (6) of Proposition \ref{prop121} and an analogue of (2) of Proposition \ref{prop6} for $R_{m,l}$:
\begin{equation*}
R_{m,l}(u;\tau) +\frac{i}{\sqrt{-i\tau}} e^{2\pi imu^2/\tau} \frac{1}{\sqrt{2m}} \sum_{\nu\smod{2m}} e^{-\pi il\nu/m} R_{m,\nu}\Bigl( \frac{u}{\tau};-\frac{1}{\tau}\Bigr) =2 h_l(u;\tau).
\end{equation*} 
We will not prove this equation. Part (1b) and equation \eqref{jacft} may also be proved by using properties of $f_u$ given in Proposition \ref{prop121} and properties of $R_{m,l}$, which we will not give here. Part (1a) follows directly from the fact that both $f_u$ and $\theta_{m,l}$ satisfy (E).
\end{remark}

\begin{proposition} \label{prop3v}
Let $R_{m,l}$ be as in Definition \ref{defft}. Then
\begin{description}
\item[(1)] if $\alpha,\beta\in\R$ then
\begin{equation*}
\frac{\partial}{\partial \overline{\tau}} e^{-2\pi im\alpha^2\tau} R_{m,l}(\alpha\tau+\beta;\tau) = -i \sqrt{\frac{m}{y}} e^{4\pi im\alpha\beta} \sum_{\lambda \in \alpha+\frac{l}{2m}+\Z} \lambda e^{-2\pi im\lambda^2\overline{\tau} -4\pi im\lambda \beta},
\end{equation*}
\item[(2)] $\tau\mapsto e^{-2\pi im\alpha^2\tau} R_{m,l}(\alpha\tau+\beta;\tau)$ is an eigenfunction of the weight 1/2 Casimir operator $\Omega_{\frac{1}{2}}=-4y^2 \frac{\partial^2}{\partial \tau \partial \overline{\tau}} +iy\frac{\partial}{\partial\overline{\tau}} +\frac{3}{16}$ with eigenvalue $\frac{3}{16}$.
\end{description}
\end{proposition}
\proof{} (1) We have 
\begin{equation*}
\begin{split}
&\frac{\partial}{\partial \overline{\tau}} e^{-2\pi im\alpha^2\tau} R_{m,l}(\alpha\tau+\beta;\tau)\\
&= \frac{1}{2}\Bigl(\frac{\partial}{\partial x}+ i\frac{\partial}{\partial y}\Bigr) \cdot \\
& \qquad \underset{\lambda\equiv l\smod{2m}}{\sum_{\lambda\in \Z}} \Bigl\{\sign\Bigl(\lambda+\frac{1}{2}\Bigr) -E\Bigl((\lambda+2m\alpha)\sqrt{y/m}\Bigr)\Bigr\}\ e^{-2\pi im(\lambda/2m+\alpha)^2 \tau-2\pi i\lambda \beta}\\
&= -\frac{i}{2} \sqrt{\frac{m}{y}} \underset{\lambda\equiv l\smod{2m}}{\sum_{\lambda\in \Z}} \Bigl(\frac{\lambda}{2m}+\alpha\Bigr)\ E'\Bigl((\lambda+2m\alpha)\sqrt{y/m}\Bigr)\ e^{-2\pi im(\lambda/2m+\alpha)^2 \tau-2\pi i\lambda \beta}\\
&= -i \sqrt{\frac{m}{y}} \underset{\lambda\equiv l\smod{2m}}{\sum_{\lambda\in \Z}} \Bigl(\frac{\lambda}{2m}+\alpha\Bigr)\  e^{-2\pi im(\lambda/2m+\alpha)^2 \overline{\tau}-2\pi i\lambda \beta}.
\end{split}
\end{equation*}
If we now substitute $\lambda'=\frac{\lambda}{2m}+\alpha$ we get the desired result.\\
(2) From (1) we see that $\tau\mapsto\sqrt{y}\frac{\partial}{\partial \overline{\tau}} e^{-2\pi im\alpha^2\tau} R_{m,l}(\alpha\tau+\beta;\tau)$ is anti-holomorphic, so
\begin{equation*}
\frac{\partial}{\partial \tau} \sqrt{y}\frac{\partial}{\partial \overline{\tau}} e^{-2\pi im\alpha^2\tau} R_{m,l}(\alpha\tau+\beta;\tau) =0
\end{equation*}
We can write the operator $\Omega_{\frac{1}{2}}=-4y^2 \frac{\partial^2}{\partial \tau \partial \overline{\tau}} +iy\frac{\partial}{\partial\overline{\tau}} +\frac{3}{16}$ as
\begin{equation*}
\Omega_{\frac{1}{2}}= \frac{3}{16}-4y^{3/2} \frac{\partial}{\partial \tau} \sqrt{y}\frac{\partial}{\partial \overline{\tau}}.
\end{equation*}
Hence 
\begin{equation*}
\Omega_{\frac{1}{2}} e^{-2\pi im\alpha^2\tau} R_{m,l}(\alpha\tau+\beta;\tau) =\frac{3}{16} e^{-2\pi im\alpha^2\tau} R_{m,l}(\alpha\tau+\beta;\tau),
\end{equation*}
which proves (2). \qed

\section{Transformation properties}

Before we can state the main result we need the following
\begin{definition}
Let $u\in\C$ and let $f$ be a real-analytic function in a neighbourhood of $u$. If $g$ is a meromorphic function with a pole of order $s$ in $u$, then $f\cdot g$ has, in a neighbourhood of $u$, an expansion 
\begin{equation*}
\sum_{n\geq -s} \sum_{m\geq 0} a_{nm} (v-u)^n (\overline{v}-\overline{u})^m.
\end{equation*}
We define
\begin{equation*} 
\Res_{v=u} \Bigl[f(v) g(v)\Bigr] = a_{-1,0}= \frac{1}{(s-1)!} \left. \frac{\partial^{s-1}}{\partial v^{s-1}}\right|_{v=u} \Bigl(f(v)\cdot (v-u)^s g(v)\Bigr) .
\end{equation*}
If $f$ is holomorphic the definition coincides with the usual definition of the residue.
\end{definition}

Now the main result:
\begin{theorem}\label{themv}
Let $\phi$ be such that $z\mapsto \phi(z;\tau)$, for fixed $\tau\in\Ha$, is a meromorphic function. If $\phi$ satisfies (E), then $\phi$ has a development of the form
\begin{equation*}
\phi(z;\tau) = \sum_{l\smod{2m}} h_l(\tau) \theta_{m,l} (z;\tau) -2\pi i \sum_{u\in\operatorname{Sing} \phi(\cdot;\tau)\smod{\Lambda_\tau}}\Res_{v=u}\ \Bigl[\tilde{f}_v(z;\tau)\phi(v;\tau)\Bigr],
\end{equation*}
for $\tau\in\Ha$ and $z\not\in\operatorname{Sing}\phi(\cdot;\tau)$, with 
\begin{equation*}
\operatorname{Sing}\phi(\cdot;\tau) := \{ u\in \C \mid z\mapsto \phi(z;\tau)\ \text{has a pole in}\ u\},
\end{equation*}
\begin{equation*}
\Lambda_\tau=\Z\tau+\Z,
\end{equation*}
and, for $0\leq l\leq 2m-1$ and any $p\in\C$ such that there are no poles on the boundary $\partial P_p$ of $P_p:=p+(0,1)\tau+(0,1)$,
\begin{equation*}
h_l(\tau) = e^{-\pi i l^2\tau/2m} \int_p^{p+1} \phi(z;\tau) e^{-2\pi ilz} dz -\pi i \sum_{u\in\operatorname{sing}_p\phi(\cdot;\tau)} \Res_{v=u}\ \Bigl[R_{m,l}(v;\tau)\phi(v;\tau)\Bigr],
\end{equation*}
with $\operatorname{sing}_p \phi(\cdot;\tau)=\operatorname{Sing}\phi(\cdot;\tau)\cap P_p$.

If $\phi$ also satisfies (M), then the vector $\bigl(h_l\bigr)_{l\smod{2m}}$ transforms under the action of $\operatorname{SL}_2(\Z)$ as in equations \eqref{them311form2} and \eqref{them311form3}.
\end{theorem}
\proof{} Let $p\in\C$ be such that $z\mapsto \phi(z;\tau)$ has no poles on $\partial P_p$. Let $z\in P_p$, $z\not\in\operatorname{Sing}\phi(\cdot;\tau)$. Now consider
\begin{equation*}
\int_{\partial P_p} f_v(z;\tau) \phi(v;\tau) dv.
\end{equation*}
We compute this integral in two different ways. On the one hand, the function we are integrating is 1-periodic. Hence
\begin{equation*}
\begin{split}
\int_{\partial P_p} &f_v(z;\tau) \phi(v;\tau) dv\\
&=\int_p^{p+1} f_v(z;\tau) \phi(v;\tau) dv - \int_p^{p+1} f_{v+\tau}(z;\tau) \phi(v+\tau;\tau) dv\\
&=\int_p^{p+1} \Bigl( f_v(z;\tau) - e^{-2\pi im\tau-4\pi imv}f_{v+\tau}(z;\tau)\Bigr) \phi(v;\tau) dv\\
&= \sum_{l=0}^{2m-1} e^{-\pi il^2\tau/2m} \theta_{m,l}(z;\tau) \int_p^{p+1} \phi(v;\tau) e^{-2\pi ilv}dv,
\end{split} 
\end{equation*}
by (4) of Proposition \ref{prop121}. On the other hand, we can compute the integral using the residue theorem. The poles of $v\mapsto f_v(z;\tau) \phi(v;\tau)$ inside $P_p$ are the poles of $\phi$ inside $P_p$, together with $z$, and the residue in $v=z$ is $\frac{1}{2\pi i} \phi(z;\tau)$. Hence
\begin{equation*}
\begin{split}
\int_{\partial P_p} &f_v(z;\tau) \phi(v;\tau) dv\\
&=2\pi i \sum_{u\in\operatorname{sing}_p \phi(\cdot;\tau)} \Res_{v=u}\ \Bigl[f_v(z;\tau) \phi(v;\tau)\Bigr] +\phi(z;\tau)\\
&= 2\pi i \sum_{u\in\operatorname{sing}_p \phi(\cdot;\tau)} \Res_{v=u}\ \Bigl[\tilde{f}_v(z;\tau) \phi(v;\tau)\Bigr] +\phi(z;\tau)\\
&\qquad +\pi i \sum_{l\smod{2m}} \theta_{m,l}(z;\tau) \sum_{u\in\operatorname{sing}_p  \phi(\cdot;\tau)} \Res_{v=u}\ \Bigl[R_{m,l}(v;\tau)\phi(v;\tau)\Bigr]
\end{split}
\end{equation*}
If we compare the two evaluations of the integral, we obtain 
\begin{equation*}
\phi(z;\tau) = \sum_{l\smod{2m}} h_l(\tau) \theta_{m,l} (z;\tau) -2\pi i \sum_{u\in\operatorname{sing}_p \phi(\cdot;\tau)} \Res_{v=u}\ \Bigl[\tilde{f}_v(z;\tau) \phi(v;\tau)\Bigr],
\end{equation*}
with $h_l$ as in the theorem. Since $v\mapsto \tilde{f}_v(z;\tau) \phi(v;\tau)$ is invariant under translation by a lattice point, so is $\displaystyle{u\mapsto \Res_{v=u}\ \Bigl[\tilde{f}_v(z;\tau)} \phi(v;\tau)\Bigr]$. Hence we can replace $\displaystyle{\sum_{u\in\operatorname{sing}_p  \phi(\cdot;\tau)}}$ by $\displaystyle{\sum_{u\in\operatorname{Sing} \phi(\cdot;\tau)\smod{\Lambda_\tau}}}$. So far we have only proven the identity for $z\in P_p$. However, both sides satisfy (E), so the identity holds for all $z\in\C$, $z\not\in\operatorname{Sing}\phi(\cdot;\tau)$.

In the rest of the proof we assume that $\phi$ satisfies (M). Let 
\begin{equation*}
\tilde{\phi}(z;\tau) := 2\pi i \sum_{u\in\operatorname{Sing} \phi(\cdot;\tau)\smod{\Lambda_\tau}}\Res_{v=u}\ \Bigl[\tilde{f}_v(z;\tau) \phi(v;\tau)\Bigr].
\end{equation*}
From the first part of the theorem we see that $\phi +\tilde{\phi}$ is a holomorphic function, which satisfies (E), and
\begin{equation*}
\phi(z;\tau) +\tilde{\phi}(z;\tau) = \sum_{l\smod{2m}} h_l(\tau) \theta_{m,l} (z;\tau)
\end{equation*}
If we can show that $\tilde{\phi}$ also satisfies (M), then the second part of the theorem follows from the second part of Theorem \ref{them311} applied to $\phi +\tilde{\phi}$. 

Let $\gamma=\left(\begin{smallmatrix} a&b\\c&d \end{smallmatrix}\right)\in\operatorname{SL}_2(\Z)$. If $u$ is a pole of $\phi \left(\cdot;\gamma\tau\right)$, then $u'=(c\tau+d)u$ is a pole of $\phi (\cdot;\tau)$, and $\Lambda_\tau =(c\tau+d) \Lambda_{\gamma\tau}$. Hence 
\begin{equation*}
\operatorname{Sing} \phi\left(\cdot;\gamma\tau\right)\smod{\Lambda_{\gamma\tau}} = (c\tau+d) \operatorname{Sing} \phi(\cdot;\tau)\smod{\Lambda_\tau}.
\end{equation*}

Using (5) of Proposition \ref{prop222} we find
\begin{equation*}
\begin{split}
\tilde{f}_{\frac{v}{c\tau +d}} \left( \frac{z}{c\tau+d};\frac{a\tau+b}{c\tau+d}\right) &\left(\frac{v}{c\tau+d}-\frac{u}{c\tau+d}\right)^s \phi \left( \frac{v}{c\tau+d};\frac{a\tau+b}{c\tau+d}\right)\\
&= (c\tau+d)^{k+1-s} e^{2\pi imcz^2/(c\tau+d)} \tilde{f}_v (z;\tau) (v-u)^s \phi(v;\tau) .
\end{split}
\end{equation*}
If we apply $\frac{1}{(s-1)!} \left.\frac{\partial^{s-1}}{\partial v^{s-1}} \right|_{v=u}$ to both sides and multiply by $(c\tau+d)^{s-1}$ we find
\begin{equation*}
\begin{split}
\Res_{v=\frac{u}{c\tau+d}} & \left[\tilde{f}_{\frac{v}{c\tau +d}} \left( \frac{z}{c\tau+d};\frac{a\tau+d}{c\tau+d}\right) \phi \left( \frac{v}{c\tau+d};\frac{a\tau+d}{c\tau+d}\right) \right]\\
&=(c\tau+d)^k e^{2\pi imcz^2/(c\tau+d)}\ \Res_{v=u}\ \Bigl[\tilde{f}_v (z;\tau) \phi (v;\tau)\Bigr].
\end{split}
\end{equation*}
Hence
\begin{equation*}
\begin{split}
\tilde{\phi}&\left(\frac{z}{c\tau+d};\frac{a\tau+b}{c\tau+d}\right) \\
&= 2\pi i\sum_{u\in\operatorname{Sing} \phi\left(\cdot;\gamma\tau\right)\smod{\Lambda_{\gamma\tau}}}  \Res_{v=u}\ \left[\tilde{f}_u \left(\frac{z}{c\tau+d};\frac{a\tau+b}{c\tau+d}\right) \phi \left( \frac{v}{c\tau+d};\frac{a\tau+d}{c\tau+d}\right)\right]\\
&= 2\pi i\sum_{u'\in\operatorname{Sing}\phi(\cdot;\tau)\smod{\Lambda_\tau}} \Res_{v=\frac{u'}{c\tau+d}} \left[\tilde{f}_u \left(\frac{z}{c\tau+d};\frac{a\tau+b}{c\tau+d}\right) \phi \left( \frac{v}{c\tau+d};\frac{a\tau+d}{c\tau+d}\right)\right]\\
&= 2\pi i\sum_{u'\in\operatorname{Sing}\phi(\cdot;\tau)\smod{\Lambda_\tau}} (c\tau+d)^k e^{2\pi imcz^2/(c\tau+d)}\ \Res_{v=u'}\ \Bigl[\tilde{f}_v (z;\tau) \phi(v;\tau)\Bigr]\\
&= (c\tau+d)^k e^{2\pi imcz^2/(c\tau+d)}\tilde{\phi}(z;\tau),
\end{split}
\end{equation*}
where in the second step we have substituted $u=\frac{u'}{c\tau+d}$. \qed

\section{Simple poles}

If all the poles of $\phi$ are simple, the theorem from the previous section reduces to

\begin{corollary}\label{cor}
Let $\phi$ be such that $z\mapsto \phi(z;\tau)$, for fixed $\tau\in\Ha$, is a meromorphic function having only simple poles. If $\phi$ satisfies (E), then
\begin{equation*}
\phi(z;\tau) = \sum_{l\smod{2m}} h_l(\tau) \theta_{m,l} (z;\tau) + \sum_{u\in\operatorname{Sing} \phi(\cdot;\tau)\smod{\Lambda_\tau}} d_u(\tau) \tilde{f}_u(z;\tau),
\end{equation*}
with 
\begin{equation*}
\operatorname{Sing}\phi(\cdot;\tau) := \{ u\in \C \mid z\mapsto \phi(z;\tau)\ \text{has a pole in}\ u\},
\end{equation*}
\begin{equation*}
\Lambda_\tau=\Z\tau+\Z,
\end{equation*}
\begin{equation*}
d_u(\tau) =-2\pi i \Res_{z=u} \phi(z;\tau),
\end{equation*}
and, for $0\leq l\leq 2m-1$ and any $p\in\C$ such that there are no poles on the boundary $\partial P_p$ of $P_p:= p+(0,1)\tau+(0,1)$,
\begin{equation*}
h_l(\tau) = e^{-\pi i l^2\tau/2m} \int_p^{p+1} \phi(z;\tau) e^{-2\pi ilz} dz +\frac{1}{2} \sum_{u\in\operatorname{sing}_p\phi(\cdot;\tau)} d_u(\tau) R_{m,l}(u;\tau),
\end{equation*}
with $\operatorname{sing}_p \phi(\cdot;\tau)=\operatorname{Sing}\phi(\cdot;\tau)\cap P_p$.

If $\phi$ also satisfies (M), then the vector $\bigl(h_l\bigr)_{l\smod{2m}}$ transforms under the action of $\operatorname{SL}_2(\Z)$ as in equations \eqref{them311form2} and \eqref{them311form3}.
\end{corollary}

In a special case, the residue function $\tau\mapsto d_u(\tau)$ has modular transformation properties:
\begin{proposition}
Let $\phi$ and $d_u(\tau)$ be as in Corollary \ref{cor} and suppose that the pole $u$ of $\phi(\cdot;\tau)$ is of the form $u=\alpha\tau+\beta$, with $\alpha,\beta\in\Q$ independent of $\tau$. Then $\tau\mapsto e^{2\pi im\alpha^2\tau} d_{u}(\tau)$ transforms as a modular form of weight $k-1$ on some subgroup $\Gamma_{\alpha,\beta}$ of $\operatorname{SL}_2(\Z)$.
\end{proposition}
\proof{}
From the definition of $d_u$ we can easily verify that 
\begin{equation*}
d_{u+\lambda\tau+\mu} (\tau) = e^{-2\pi im(\lambda^2\tau+2\lambda u)}d_u(\tau),
\end{equation*}
and 
\begin{equation*}
d_{\frac{u}{c\tau+d}} \left(\frac{a\tau+b}{c\tau+d} \right) = (c\tau+d)^{k-1} e^{2\pi im cu^2/(c\tau+d)} d_u(\tau). 
\end{equation*}
Hence $(u,\tau)\mapsto d_u(\tau)$ transforms as a Jacobi form of weight $k-1$ and index $m$. By Theorem 1.3 of \cite[pp.\ 10]{eichler} we get the desired result. Actually, that theorem also assumes a growth condition, but one can check in the proof that the growth condition is not needed to prove the modular transformation properties.
\qed  

\section{An example}

Define $\phi$ by:
\begin{equation*}
\begin{split}
\phi(z;\tau) &:= \frac{\left( \theta_{0,0}(z;\tau) \theta_{0,\frac{1}{2}}(z;\tau) \theta_{\frac{1}{2},0}(z;\tau)\right)^9}{\Delta(\tau) \theta_{\frac{1}{2},\frac{1}{2}}(z;\tau)}\\
&= -i\Bigl(\zeta^{\frac{1}{2}}+\zeta^{-\frac{1}{2}}\Bigr)^9 \left\{ \frac{1}{\zeta^{\frac{1}{2}}-\zeta^{-\frac{1}{2}}}-\Bigl( 9 \zeta^{\frac{3}{2}}-\zeta^{\frac{1}{2}}+\zeta^{-\frac{1}{2}}-9\zeta^{-\frac{3}{2}}\Bigr) q+ \ldots\right\},
\end{split}
\end{equation*}
with $\theta_{a,b}(z;\tau) := \sum_{\lambda\in a+\Z}e^{\pi i \lambda^2 \tau +2\pi i\lambda(z+b)}$, $\zeta=e^{2\pi iz}$ and $q=e^{2\pi i\tau}$. 

Using Table V on page 36 of \cite{mumford} we see that $\phi$ transforms like a Jacobi form of weight $k=1$ and index $m=13$ on the full modular group. (Note that Mumford uses the notation $\theta_{01}$, $\theta_{10}$ and $\theta_{11}$ for the functions denoted here by $\theta_{0,\frac{1}{2}}$, $\theta_{\frac{1}{2},0}$ and $\theta_{\frac{1}{2},\frac{1}{2}}$. Also there's a mistake in the 4th formula on the right: it should read $\theta_{11}(z/\tau,-1/\tau)=-i(-i\tau)^{\frac{1}{2}} \exp(\pi iz^2/\tau) \theta_{11}(z,\tau)$.)

The function $\phi$ is meromorphic in $z$ with simple poles in $\Z\tau+\Z$. If we take $p=-\frac{1}{2}\tau-\frac{1}{2}$ then $\operatorname{sing}_p \phi(\cdot;\tau) =\{0\}$. Further
\begin{equation*}
\begin{split}
\operatornamewithlimits{Res}_{z=0}\ \phi(z;\tau)&=\frac{\left( \theta_{0,0}(0;\tau) \theta_{0,\frac{1}{2}}(0;\tau) \theta_{\frac{1}{2},0}(0;\tau)\right)^9}{\Delta(\tau) \theta_{\frac{1}{2},\frac{1}{2}}'(0;\tau)}= -\frac{1}{\pi^9} \frac{\theta_{\frac{1}{2},\frac{1}{2}}'(0;\tau)^8}{\Delta (\tau)}\\
&= -\frac{1}{\pi^9} \frac{(-2\pi \eta(\tau)^3)^8}{\Delta (\tau)}= -\frac{128}{\pi},
\end{split}
\end{equation*}
where we have used Jacobi's derivative formula (see \cite[pp.\ 64]{mumford}), and (10) of Proposition \ref{prop2}.

Corollary \ref{cor} gives:
\begin{equation}\label{deve}
\phi(z;\tau) = \sum_{l\smod{26}} h_l(\tau) \theta_{13,l} (z;\tau) +512i\tilde{f}_0(z;\tau)
\end{equation}
with
\begin{equation}\label{ha}
h_l(\tau) = e^{-\pi i l^2\tau/26} \int_p^{p+1} \phi(z;\tau) e^{-2\pi ilz} dz +256i R_{13,l}(0;\tau),
\end{equation}
for $0\leq l\leq 25$.
According to the corollary the $h_l$ transform as a vector-valued modular form of weight $\frac{1}{2}$. 

Since $\phi$ is holomorphic as a function of $\tau$, the Fourier coefficients
\begin{equation*}
e^{-\pi i l^2\tau/26} \int_p^{p+1} \phi(z;\tau) e^{-2\pi ilz} dz
\end{equation*}
are holomorphic as a function of $\tau$. In particular they are eigenfunctions of $\Omega_{\frac{1}{2}}$ with eigenvalue $\frac{3}{16}$. Using (2) of Proposition \ref{prop3v} we see that $\tau\mapsto R_{13,l}(0;\tau)$ is also a eigenfunction of $\Omega_{\frac{1}{2}}$ with eigenvalue $\frac{3}{16}$. Hence 
\begin{equation*}
\Omega_{\frac{1}{2}} h_l =\frac{3}{16} h_l.
\end{equation*}
So the $h_l$ form a vector-valued real-analytic modular form of weight $\frac{1}{2}$. The transformations are:
\begin{equation*}
h_l(\tau+1) = e^{-\pi il^2/26} h_l(\tau)
\end{equation*}
and 
\begin{equation*}
h_l \Bigl(-\frac{1}{\tau}\Bigr) = i\sqrt{-i\tau} \frac{1}{\sqrt{26}} \sum_{\nu\smod{26}} e^{\pi i l\nu/13}h_\nu(\tau)
\end{equation*}

Using (1) of Proposition \ref{prop3v} we see
\begin{equation*}
\frac{\partial}{\partial \overline{\tau}} h_l(\tau) = 256 \sqrt{13}\ y^{-\frac{1}{2}} \sum_{\lambda \in \frac{l}{26}+\Z} \lambda e^{-26\pi i\lambda^2\overline{\tau} }.
\end{equation*}
Note that $\sum_{\lambda \in \frac{l}{26}+\Z} \lambda e^{-26\pi i\lambda^2\overline{\tau} }$ is the complex conjugate of $\sum_{\lambda \in \frac{l}{26}+\Z} \lambda e^{26\pi i\lambda^2\tau }$, which is a theta function of weight $3/2$.

Summarizing, we have proved:

\begin{proposition}
Let $\phi$ be the function given by
\begin{equation*}
\phi(z;\tau) = \frac{\left( \theta_{0,0}(z;\tau) \theta_{0,\frac{1}{2}}(z;\tau) \theta_{\frac{1}{2},0}(z;\tau)\right)^9}{\Delta(\tau) \theta_{\frac{1}{2},\frac{1}{2}}(z;\tau)}.
\end{equation*}
Then $\phi$ can be decomposed as in \eqref{deve}, with $h_l$ as in \eqref{ha}, and $\bigl(h_l\bigr)_{l\smod{26}}$ is a vector-valued real-analytic modular form of weight 1/2, with eigenvalue 3/16 for the weight 1/2 Casimir operator.
\end{proposition}

This is a very special example: It has been constructed is such a way that the $d_u$ are constant (as a function of $\tau$). As a result the $h_l$ are eigenfunctions of a Casimir operator. However, in general the $d_u$ will not be constant and the $h_l$ will not be eigenfunctions of a Casimir operator. So in general we do not end up with a real-analytic modular form.

\clearemptydoublepage
\chapter{Mock $\theta$-functions}
\section{Introduction}

Mock $\theta$-functions were introduced by S. Ramanujan  in the last letter he wrote to G.H. Hardy, dated January, 1920. For a photocopy of the mathematical part of this letter see \cite[pp. 127--131]{lost} (also reproduced in \cite{andrews1}). In this letter, Ramanujan provided a list of 17 mock $\theta$-functions, together with identities they satisfy. Ramanujan divided his list of functions into ``third order'', ``fifth order'' and ``seventh order'' functions, but did not say what he meant. There's still no formal definition of ``order'', but known identities for these mock $\theta$-functions make it clear that they are related to the numbers 3, 5 and 7. Therefore we regard the order of a mock $\theta$-function merely as a convenient label, which may or may not have a deeper significance. 

In his letter, Ramanujan explained what he meant by a mock $\theta$-function.
In \cite{andrews2} we find a formal definition. Slightly rephrased it reads: a mock $\theta$-function is a function $f$ of the complex variable $q$, defined by a $q$-series of a particular type (Ramanujan calls this the Eulerian form), which converges for $|q|<1$ and satisfies the following conditions:
\begin{itemize}
\item[(1)] infinitely many roots of unity are exponential singularities,
\item[(2)] for every root of unity $\xi$ there is a $\theta$-function $\theta_{\xi} (q)$ such that the difference $f(q)-\theta_{\xi} (q)$ is bounded as $q\rightarrow \xi$ radially (presumably with only finitely many of the $\theta_{\xi}$ being different),
\item[(3)] there is no $\theta$-function that works for all $\xi$, i.e.\ $f$ is not the sum of two functions, one of which is a $\theta$-function and the other a function which is bounded in all roots of unity.
\end{itemize}
(When Ramanujan refers to $\theta$-functions, he means sums, products, and quotients of series of the form $\sum_{n\in\Z} \epsilon^n q^{an^2 +bn}$ with $a,b\in\Q$ and $\epsilon=-1,1$).

The 17 functions given by Ramanujan indeed satisfy conditions (1) and (2) (see \cite{watson1}, \cite{watson2} and \cite{selberg}). However no proof has ever been given that they also satisfy condition (3). Watson (see \cite{watson1}) proved a very weak form of condition (3) for the ``third order'' mock $\theta$-functions, namely, that they are not equal to $\theta$-functions.

In this chapter we will see that condition (3) is not satisfied if we strengthen it slightly. Indeed, we shall discuss vector-valued mock $\theta$-functions $F$ for which there is a vector-valued real-analytic modular form $H$ such that $F-H$ is bounded in all roots of unity.

There are several ways to get these results: For example let us consider the ``fifth order'' mock $\theta$-function (using Watson's notation)
\begin{equation*}
f_0(q) = \sum_{n=0}^\infty \frac{q^{n^2}}{(-q;q)_n},
\end{equation*}
with $(a)_n=(a;q)_n:= (1-a)(1-aq)\cdots (1-aq^{n-1})$. Andrews (see \cite{andrews3}) showed that
\begin{equation}\label{in5}
f_0(q) = \frac{1}{(q)_\infty} \sum_{n\geq0} \sum_{|j|\leq n} (-1)^j q^{\frac{5}{2}n^2+\frac{1}{2}n -j^2}(1-q^{4n+2}),
\end{equation} 
with $(q)_\infty =\prod_{n=1}^\infty (1-q^n)$. Using this identity, Andrews (see \cite{coeff}) showed that $f_0(q)$ can be seen as a Fourier coefficient of a certain quotient of Jacobi theta functions. Next Hickerson (see \cite{hick}) showed that $f_0(q)$ can be related to a sum similar to the Lerch sum discussed in Chapter I. Similar results have been found for most other mock $\theta$-functions. Hence we can extend the results of Chapter I to include these ``Lerch-like'' sums and thereby get the transformation properties of $f_0$. We may also use the techniques from Chapter 3 and the representation of the mock $\theta$-functions as a Fourier coefficient of a meromorphic Jacobi form, to find the transformation properties. However, we will use \eqref{in5} and similar identities for the other mock $\theta$-functions and apply the results from Chapter 2.

\section{General results}

Before we start with the mock $\theta$-functions, we derive some general results, which we will use repeatedly.

\begin{definition}
For $a,b\in\R$ and $\tau\in\Ha$ we define
\begin{equation*}
R_{a,b}(\tau)=\sum_{\nu\in a+\Z} \sign(\nu) \beta (2\nu^2y)\ e^{-\pi i\nu^2 \tau -2\pi i\nu b},
\end{equation*}
with $y=\im(\tau)$ and $\beta$ as in Lemma \ref{lem131}.
\end{definition}\
This function, a kind of non-holomorphic unary theta-series, is a slight modification of the function $R$ studied in Chapter 1 (cf.\ (1) of the following proposition) and is also similar to the function $R_{m,l}(z;\tau)$ defined in Definition \ref{defft}.

\begin{proposition}\label{propRab}
Let $a\in(0,1)$, $b\in\R$ and $\tau\in\Ha$. We have
\begin{itemize}
\item[(1)] $R_{a,b} (\tau) =i e^{-\pi i(a-\frac{1}{2})^2\tau-2\pi i(a-\frac{1}{2})b}\ R\Bigl((a-\frac{1}{2})\tau +b+\frac{1}{2};\tau\Bigr)$, with $R$ as defined in Lemma \ref{defR}.
\item[(2)] $\displaystyle{R_{a,b} (\tau) =-i \int_{-\overline{\tau}}^{i\infty} \frac{g_{a,-b}(z)}{\sqrt{-i(z+\tau)}}\, dz}$, with $g_{a,b}$ as in Definition \ref{defg}.
\item[(3)] If $\xi\in\Q$, then $R_{a,b}(\tau)$ is bounded as $\tau \downarrow \xi$.
\item[(4)] $\frac{\partial}{\partial \overline{\tau}} R_{a,b}(\tau)=-i\frac{1}{\sqrt{2y}}\ g_{a,-b}(-\overline{\tau})$.
\item[(5)] $\tau\mapsto R_{a,b}(\tau)$ is an eigenfunction of the weight 1/2 Casimir operator $\Omega_{\frac{1}{2}}=-4y^2 \frac{\partial^2}{\partial \tau \partial \overline{\tau}} +iy\frac{\partial}{\partial\overline{\tau}} +\frac{3}{16}$, with eigenvalue $\frac{3}{16}$.
\end{itemize}
\end{proposition}
\proof (1) By definition we have
\begin{equation*}
\begin{split}
R&\Bigl(\Bigl(a-\frac{1}{2}\Bigr)\tau+b+\frac{1}{2};\tau\Bigr) \\
&= \sum_{\nu\in \frac{1}{2}+\Z} \Bigl\{ \sign(\nu) - E\Bigl(\Bigl(\nu+a-\frac{1}{2}\Bigr)\sqrt{2y}\Bigr)\Bigr\}\  (-1)^{\nu-\frac{1}{2}} e^{-\pi i \nu^2\tau -2\pi i\nu ((a-\frac{1}{2})\tau+b+\frac{1}{2})}.
\end{split}
\end{equation*}
Using Lemma \ref{lem131}, we write $\sign(\nu) - E\left(\left(\nu+a-\frac{1}{2}\right)\sqrt{2y}\right)$ as the sum of $\sign(\nu)-\sign(\nu+a-\frac{1}{2})$ and $\sign(\nu+a-\frac{1}{2})\beta(2(\nu+a-\frac{1}{2})^2 y)$. We see that $\sign(\nu)-\sign(\nu+a-\frac{1}{2})=0$ for all $\nu\in\frac{1}{2}+\Z$, since $a\in(0,1)$. Hence
\begin{equation*}
\begin{split}
R&\Bigl(\Bigl(a-\frac{1}{2}\Bigr)\tau+b+\frac{1}{2};\tau\Bigr)=\\ &\sum_{\nu\in\frac{1}{2}+\Z}\sign\Bigl(\nu+a-\frac{1}{2}\Bigr)\beta\Bigl(2\Bigl(\nu+a-\frac{1}{2}\Bigr)^2 y\Bigr)\ (-1)^{\nu-\frac{1}{2}} e^{-\pi i \nu^2\tau -2\pi i\nu ((a-\frac{1}{2})\tau+b+\frac{1}{2})}\\
&= e^{\pi i(a-\frac{1}{2})^2 \tau +2\pi i(a-\frac{1}{2})(b+\frac{1}{2})}\sum_{\nu\in a+\Z}\sign(\nu)\beta(2\nu^2 y)\ (-1)^{\nu-a} e^{-\pi i \nu^2\tau -2\pi i\nu (b+\frac{1}{2})}\\
&= -i e^{\pi i(a-\frac{1}{2})^2 \tau +2\pi i(a-\frac{1}{2})b}\ R_{a,b}(\tau),
\end{split}
\end{equation*}
where we have substituted $\nu \rightarrow \nu -a +\frac{1}{2}$ in the second step.\\
(2) Use (1) of this proposition and (1) of Theorem \ref{period} with $a$ replaced by $a-\frac{1}{2}$ and $b$ replaced by $-b-\frac{1}{2}$.\\
(3) This follows directly from (2) and the fact that $\lim_{z \downarrow \xi} g_{a,-b}(z)=0$ for all $\xi\in\Q$.\\
(4) This follows directly from (2) by taking $\frac{\partial}{\partial\overline{\tau}}$ on both sides.\\
(5) From (4) we see that $\tau\mapsto\sqrt{y}\frac{\partial}{\partial \overline{\tau}} R_{a,b}(\tau)$ is anti-holomorphic, so
\begin{equation*}
\frac{\partial}{\partial \tau} \sqrt{y}\frac{\partial}{\partial \overline{\tau}} R_{a,b}(\tau) =0.
\end{equation*}
We can write the operator $\Omega_{\frac{1}{2}}=-4y^2 \frac{\partial^2}{\partial \tau \partial \overline{\tau}} +iy\frac{\partial}{\partial\overline{\tau}} +\frac{3}{16}$ as
\begin{equation*}
\Omega_{\frac{1}{2}}= \frac{3}{16}-4y^{3/2} \frac{\partial}{\partial \tau} \sqrt{y}\frac{\partial}{\partial \overline{\tau}}.
\end{equation*}
Hence $\Omega_{\frac{1}{2}} R_{a,b}=\frac{3}{16} R_{a,b}$. \qed

We now return to the general setup of Chapter 2, i.e.\ indefinite $\theta$-series for a quadratic form $Q$ of type $(r-1,1)$. In the next proposition, we will rewrite
\begin{equation*}
\sum_{\nu\in a+\Z^r} \sign\Bigl(B(c,\nu)\Bigr) \beta \left(-\frac{B(c,\nu)^2}{Q(c)}y\right) e^{2\pi iQ(\nu)\tau +2\pi iB(\nu,b)},
\end{equation*}
(this is the same series as in \eqref{221}) for $c\in C_Q \cap \Z^r$.
In order to do so, we write $\nu=\mu+nc$ with $\mu\in a+\Z^r$, $n\in\Z$, such that $\frac{B(c,\mu)}{2Q(c)}\in [0,1)$. Since $c\in\Z^r$ we can write
\begin{equation*}
\left\{ \mu\in a+\Z^r \left| \frac{B(c,\mu)}{2Q(c)} \in [0,1)\right.\right\}=\bigsqcup_{\mu_0\in P_0} \left( \mu_0 + \left<c\right>_\Z^{\perp} \right)
\end{equation*}
(disjoint union), for a suitable finite set $P_0$, with $\left<c\right>_\Z^{\perp}:=\{ \xi\in\Z^r \mid B(c,\xi)=0\}$. We can now state the result:
\begin{proposition}\label{split}
Let $c\in\C_Q\cap \Z^r$ be primitive. Then there is a finite set $P_0$ (see above), such that 
\begin{equation*}
\begin{split}
\sum_{\nu\in a+\Z^r} &\sign\Bigl(B(c,\nu)\Bigr) \beta \left(-\frac{B(c,\nu)^2}{Q(c)}y\right) e^{2\pi iQ(\nu)\tau +2\pi iB(\nu,b)}\\
&= -\sum_{\mu_0 \in P_0} R_{\frac{B(c,\mu_0)}{2Q(c)},-B(c,b)} (-2Q(c)\tau) \cdot \sum_{\xi\in\mu_0^\perp + \left< c\right>_\Z^\perp} e^{2\pi iQ(\xi)\tau+2\pi iB(\xi,b^\perp)},
\end{split}
\end{equation*}
with $\mu_0^\perp= \mu_0 - \frac{B(c,\mu_0)}{2Q(c)}c$ and $b^\perp= b-\frac{B(c,b)}{2Q(c)}c$.
\end{proposition}
\begin{remark}
Since $\mu_0^\perp + \left< c\right>_\Z^\perp$ is a shifted $(r-1)$-dimensional lattice, on which $Q$ is positive definite, the inner sum
\begin{equation*}
\sum_{\xi\in\mu_0^\perp + \left< c\right>_\Z^\perp} e^{2\pi iQ(\xi)\tau+2\pi iB(\xi,b^\perp)}
\end{equation*}
is a classical (positive definite) theta function, and is in particular modular of weight $(r-1)/2$.
\end{remark}
\proof We write $\nu=\mu_0 +\xi +nc$, with $\mu_0\in P_0$, $\xi\in \left<c\right>_\Z^\perp$ and $n\in\Z$. Set $\mu_0^\perp =\mu_0 - \frac{B(c,\mu_0)}{2Q(c)}c$. Then $B(c,\mu_0^\perp)=0$ and
\begin{equation*}
\nu=\Bigl(n+\frac{B(c,\mu_0)}{2Q(c)}\Bigr)c+\mu_0^\perp +\xi,
\end{equation*}
so
\begin{equation*}
\begin{split}
\sum_{\nu\in a+\Z^r} &\sign \Bigl(B(c,\nu)\Bigr) \beta \left(-\frac{B(c,\nu)^2}{Q(c)}y\right) e^{2\pi iQ(\nu)\tau +2\pi iB(\nu,b)}\\
&= - \sum_{\mu_0 \in P_0} \sum_{\xi\in\left< c\right>_\Z^\perp} \sum_{n\in\Z} \sign\Bigl( n+\frac{B(c,\mu_0)}{2Q(c)}\Bigr) \beta\left(-4Q(c)\Bigl( n+\frac{B(c,\mu_0)}{2Q(c)}\Bigr)^2y\right)\cdot\\
&\qquad \cdot e^{2\pi iQ(c)(n+\frac{B(c,\mu_0)}{2Q(c)})^2\tau +2\pi iQ(\mu_0^\perp +\xi)\tau+2\pi iB(c,b)(n+\frac{B(c,\mu_0)}{2Q(c)})+2\pi iB(\mu_0^\perp +\xi,b)}\\
&= -\sum_{\mu_0 \in P_0} R_{\frac{B(c,\mu_0)}{2Q(c)},-B(c,b)} (-2Q(c)\tau) \cdot \sum_{\alpha\in\mu_0^\perp + \left< c\right>_\Z^\perp} e^{2\pi iQ(\alpha)\tau+2\pi iB(\alpha,b)}.
\end{split}
\end{equation*}
If we use $B(\alpha,b)=B(\alpha,b^\perp)$ for all $\alpha\in\mu_0^\perp + \left< c\right>_\Z^\perp$, we get the desired result. \qed

\section{The seventh order mock $\theta$-functions}

In this section we deal with the ``seventh order'' mock $\theta$-functions from Ramanujan's letter. 
In \cite[pp.\ 666]{hick2} we find the following (slightly rewritten) identities:
\begin{equation*}
\begin{split}
(q)_\infty \ef_0(q) &= \Bigl( \sum_{r,s\geq0} -\sum_{r,s<0}\Bigr) (-1)^{r+s} q^{\frac{3}{2} r^2 +4rs+\frac{3}{2}s^2+\frac{1}{2}r+\frac{1}{2}s}\\
(q)_\infty \ef_1(q) &= \Bigl( \sum_{r,s\geq0} -\sum_{r,s<0}\Bigr) (-1)^{r+s} q^{\frac{3}{2} r^2 +4rs+\frac{3}{2}s^2+\frac{5}{2}r+\frac{5}{2}s+1}\\
(q)_\infty \ef_2(q) &= \Bigl( \sum_{r,s\geq0} -\sum_{r,s<0}\Bigr) (-1)^{r+s} q^{\frac{3}{2} r^2 +4rs+\frac{3}{2}s^2+\frac{3}{2}r+\frac{3}{2}s}.
\end{split}
\end{equation*}
We will use these as the definitions of the mock $\theta$-functions. We rewrite these identities ($\zeta_n := e^{2\pi i/n}$):
\begin{equation*}
\begin{split}
2\eta(\tau) \zeta_{14}\ q^{-\frac{1}{168}} \ef_0(q) &= \sum_{\nu\in \frac{1}{14}e +\Z^r} \Bigl\{\sign(B(\nu,c_1))-\sign(B(\nu,c_2))\Bigr\}\ e^{2\pi iQ(\nu)\tau+2\pi i B(\nu,\frac{1}{14}e)}\\
2\eta(\tau) \zeta_{14}\ q^{\frac{47}{168}} \ef_2(q) &= \sum_{\nu\in \frac{3}{14}e +\Z^r} \Bigl\{\sign(B(\nu,c_1))-\sign(B(\nu,c_2))\Bigr\}\ e^{2\pi iQ(\nu)\tau+2\pi i B(\nu,\frac{5}{14}e)}\\
2\eta(\tau) \zeta_{14}\ q^{-\frac{25}{168}} \ef_1(q) &= \sum_{\nu\in \frac{5}{14}e +\Z^r} \Bigl\{\sign(B(\nu,c_1))-\sign(B(\nu,c_2))\Bigr\}\ e^{2\pi iQ(\nu)\tau+2\pi i B(\nu,\frac{3}{14}e)},
\end{split}
\end{equation*}
with $A=\left(\begin{smallmatrix} 3&4\\4&3 \end{smallmatrix}\right)$, $c_1=\left(\begin{smallmatrix} -3\\4 \end{smallmatrix}\right)$, $c_2=\left(\begin{smallmatrix} -4\\3 \end{smallmatrix}\right)$ and $e:=\left(\begin{smallmatrix} 1\\1 \end{smallmatrix}\right)$.
We have $B(c_1,c_2)=-28$ and $Q(c_1)=Q(c_2)=-\frac{21}{2}$. If we choose $C_Q$ such that $c_1\in C_Q$ then also $c_2\in C_Q$. 

We collect these three mock $\theta$-functions into a single vector-valued mock $\theta$-function
\begin{align*}
F_7(\tau) &:=\begin{pmatrix} q^{-\frac{1}{168}}\ef_0(q)\\q^{\frac{47}{168}}\ef_2(q)\\q^{-\frac{25}{168}} \ef_1(q) \end{pmatrix}.\\
\intertext{To express its modular transformation behaviour, we also introduce}
H_7(\tau) &:= \frac{\zeta_{14}^{-1}}{2\eta(\tau)} \begin{pmatrix} \theta_{\frac{1}{14}e,\frac{1}{14}e}\\
\theta_{\frac{3}{14}e,\frac{5}{14}e}\\
\theta_{\frac{5}{14}e,\frac{3}{14}e}\end{pmatrix}(\tau),\\
\intertext{with $\theta_{a,b}$ as in Definition \ref{deftheta}, and}
G_7(\tau) &:=-\begin{pmatrix} \zeta_{84}^{-13} R_{\frac{13}{42},-\frac{1}{2}} + \zeta_{84} R_{\frac{41}{42},-\frac{1}{2}}\\  \zeta_{84}^{29} R_{\frac{11}{42},-\frac{5}{2}} + \zeta_{84}^{-41} R_{\frac{25}{42},-\frac{5}{2}}\\ 
 \zeta_{28}^5 R_{\frac{23}{42},-\frac{3}{2}} + \zeta_{28}^{-9} R_{\frac{37}{42},-\frac{3}{2}} \end{pmatrix}(21\tau).
\end{align*}
Note that the components of $H_7$ are the quotient of a (real-analytic) binary theta series by $\eta$, while the components of $G_7$ are (real-analytic) unary theta series.

\begin{proposition}\label{propm7}
We have
\begin{equation*}
F_7=H_7+G_7,
\end{equation*}
where
\begin{description}
\item[(1)] The function $H_7$ is a (vector-valued) real-analytic modular form of weight 1/2, satisfying
\begin{align*}
H_7(\tau+1) &= \begin{pmatrix} \zeta_{168}^{-1}&0&0\\ 0&\zeta_{168}^{47}&0\\0&0&\zeta_{168}^{-25}\end{pmatrix} H_7(\tau),\\
\intertext{and}
H_7\Bigl(-\frac{1}{\tau}\Bigr) &=\sqrt{-i\tau} \frac{2}{\sqrt{7}} \begin{pmatrix} \sin \frac{\pi}{7}&\sin \frac{3\pi}{7}&\sin \frac{2\pi}{7}\\ \sin \frac{3\pi}{7}&-\sin \frac{2\pi}{7}&\sin \frac{\pi}{7}\\ \sin \frac{2\pi}{7}&\sin \frac{\pi}{7}&-\sin \frac{3\pi}{7}\end{pmatrix} H_7(\tau),
\end{align*}
and is an eigenfunction of the Casimir operator $\Omega_{\frac{1}{2}}$, with eigenvalue $\frac{3}{16}$.
\item[(2)] The function $G_7$ is bounded if $\tau\downarrow \xi$, with $\xi\in\Q$.
\end{description}
\end{proposition}
\proof 
We consider the functions $\theta_{\frac{1}{14}e,\frac{1}{14}e}$, $\theta_{\frac{3}{14}e,\frac{5}{14}e}$ and $\theta_{\frac{5}{14}e,\frac{3}{14}e}$. Using (4) and (2) of Corollary \ref{cor1} we see
\begin{equation*}
\begin{pmatrix}\theta_{\frac{1}{14}e,\frac{1}{14}e}\\
\theta_{\frac{3}{14}e,\frac{5}{14}e}\\
\theta_{\frac{5}{14}e,\frac{3}{14}e}\end{pmatrix}(\tau+1)= \begin{pmatrix} \zeta_{28} &0&0\\0&\zeta_{28}^9&0\\0&0&\zeta_{28}^{25} \end{pmatrix} \begin{pmatrix}\theta_{\frac{1}{14}e,\frac{1}{14}e}\\
\theta_{\frac{3}{14}e,\frac{5}{14}e}\\
\theta_{\frac{5}{14}e,\frac{3}{14}e}\end{pmatrix}(\tau).
\end{equation*}
Using Corollary \ref{cor1} and $\theta_{s,\frac{1}{14}e}=-e^{2\pi i(s_1+s_2)} \theta_{e-s,\frac{1}{14}e}$ for all $s=\left(\begin{smallmatrix} s_1\\s_2 \end{smallmatrix} \right)\in\R^2$, which we get from (1), (2) and (3) of Corollary \ref{cor1}, we obtain
\begin{equation*}
\begin{pmatrix}\theta_{\frac{1}{14}e,\frac{1}{14}e}\\
\theta_{\frac{3}{14}e,\frac{5}{14}e}\\
\theta_{\frac{5}{14}e,\frac{3}{14}e}\end{pmatrix}\Bigl(-\frac{1}{\tau}\Bigr)=
-i\tau\frac{2}{\sqrt{7}} \begin{pmatrix} \sin \frac{\pi}{7}&\sin \frac{3\pi}{7}&\sin \frac{2\pi}{7}\\ \sin \frac{3\pi}{7}&-\sin \frac{2\pi}{7}&\sin \frac{\pi}{7}\\ \sin \frac{2\pi}{7}&\sin \frac{\pi}{7}&-\sin \frac{3\pi}{7}\end{pmatrix} \begin{pmatrix}\theta_{\frac{1}{14}e,\frac{1}{14}e}\\
\theta_{\frac{3}{14}e,\frac{5}{14}e}\\
\theta_{\frac{5}{14}e,\frac{3}{14}e}\end{pmatrix}(\tau).
\end{equation*}
If we use $\eta(\tau+1) =\zeta_{24} \eta(\tau)$ and $\eta\left(-\frac{1}{\tau}\right)= \sqrt{-i\tau}\ \eta(\tau)$, we get the transformations for $H_7$.

Using Proposition \ref{split} we see ($P_0=\left\{ -\frac{13}{14}e, -\frac{27}{14}e, -\frac{41}{14}e\right\}$ and $\left<c\right>_\Z^\perp = \{\xi\in\Z^2 \mid \xi_1=0\} = \left\{\left.\left(\begin{smallmatrix}0\\\xi_2\end{smallmatrix}\right)\right|\xi_2\in\Z\right\}$)
\begin{equation*}
\begin{split}
&\sum_{\nu\in \frac{1}{14}e+\Z^r} \sign (B(c_1,\nu)) \beta \left(-\frac{B(c_1,\nu)^2}{Q(c_1)}y\right) e^{2\pi iQ(\nu)\tau +2\pi iB(\nu,\frac{1}{14}e)}\\
&= -R_{\frac{13}{42},-\frac{1}{2}}(21\tau) \sum_{\mu_2\in -\frac{13}{6}+\Z} e^{3\pi i\mu_2^2\tau +\pi i\mu_2} -R_{\frac{27}{42},-\frac{1}{2}}(21\tau) \sum_{\mu_2\in -\frac{27}{6}+\Z} e^{3\pi i\mu_2^2\tau +\pi i\mu_2}\\
&\quad -R_{\frac{41}{42},-\frac{1}{2}}(21\tau) \sum_{\mu_2\in -\frac{41}{6}+\Z} e^{3\pi i\mu_2^2\tau +\pi i\mu_2}\\
&= -\eta(\tau)\left( \zeta_{12}^{-1} R_{\frac{13}{42},-\frac{1}{2}}(21\tau) +\zeta_{12} R_{\frac{41}{42},-\frac{1}{2}}(21\tau)\right).
\end{split}
\end{equation*}
Similarly we find
\begin{equation*}
\begin{split}
\sum_{\nu\in \frac{3}{14}e+\Z^r} \sign (B(c_1,\nu)) &\beta \left(-\frac{B(c_1,\nu)^2}{Q(c_1)}y\right) e^{2\pi iQ(\nu)\tau +2\pi iB(\nu,\frac{5}{14}e)}\\
&= -\eta(\tau)\left( \zeta_{12}^5 R_{\frac{11}{42},-\frac{5}{2}}(21\tau) +\zeta_{12}^{-5} R_{\frac{25}{42},-\frac{5}{2}}(21\tau)\right)\\
\sum_{\nu\in \frac{5}{14}e+\Z^r} \sign (B(c_1,\nu)) &\beta \left(-\frac{B(c_1,\nu)^2}{Q(c_1)}y\right) e^{2\pi iQ(\nu)\tau +2\pi iB(\nu,\frac{3}{14}e)}\\
&= -\eta(\tau)\left( \zeta_4 R_{\frac{23}{42},-\frac{3}{2}}(21\tau) +\zeta_4^{-1} R_{\frac{37}{42},-\frac{3}{2}}(21\tau)\right),
\end{split}
\end{equation*}
and
\begin{align*}
\sum_{\nu\in \frac{1}{14}e+\Z^r} \sign (B(c_2,\nu)) &\beta \left(-\frac{B(c_2,\nu)^2}{Q(c_2)}y\right) e^{2\pi iQ(\nu)\tau +2\pi iB(\nu,\frac{1}{14}e)}\\
&= \eta(\tau)\left( \zeta_{12}^{-1} R_{\frac{13}{42},-\frac{1}{2}}(21\tau) +\zeta_{12} R_{\frac{41}{42},-\frac{1}{2}}(21\tau)\right)\displaybreak[0]\\
\sum_{\nu\in \frac{3}{14}e+\Z^r} \sign (B(c_2,\nu)) &\beta \left(-\frac{B(c_2,\nu)^2}{Q(c_2)}y\right) e^{2\pi iQ(\nu)\tau +2\pi iB(\nu,\frac{5}{14}e)}\\
&= \eta(\tau)\left( \zeta_{12}^5 R_{\frac{11}{42},-\frac{5}{2}}(21\tau) +\zeta_{12}^{-5} R_{\frac{25}{42},-\frac{5}{2}}(21\tau)\right)\displaybreak[0]\\
\sum_{\nu\in \frac{5}{14}e+\Z^r} \sign (B(c_2,\nu)) &\beta \left(-\frac{B(c_2,\nu)^2}{Q(c_2)}y\right) e^{2\pi iQ(\nu)\tau +2\pi iB(\nu,\frac{3}{14}e)}\\
&= \eta(\tau)\left( \zeta_4 R_{\frac{23}{42},-\frac{3}{2}}(21\tau) +\zeta_4^{-1} R_{\frac{37}{42},-\frac{3}{2}}(21\tau)\right).
\end{align*}
Hence, if we write $\rho(\nu;\tau)$ as the sum of the three expressions \eqref{al1}, \eqref{al2} and \eqref{al3}, we find
\begin{equation*}
\begin{split}
\theta_{\frac{1}{14}e,\frac{1}{14}e}(\tau) &= \sum_{\nu\in \frac{1}{14}e +\Z^r} \bigl\{\sign(B(\nu,c_1))-\sign(B(\nu,c_2))\Bigr\}\ e^{2\pi iQ(\nu)\tau+2\pi i B(\nu,\frac{1}{14}e)}\\
&\quad +2\eta(\tau)\left( \zeta_{12}^{-1} R_{\frac{13}{42},-\frac{1}{2}}(21\tau) +\zeta_{12} R_{\frac{41}{42},-\frac{1}{2}}(21\tau)\right)\\
&= 2\zeta_{14}\ q^{-\frac{1}{168}} \eta(\tau) \ef_0(q)+2\eta(\tau)\left( \zeta_{12}^{-1} R_{\frac{13}{42},-\frac{1}{2}}(21\tau) +\zeta_{12} R_{\frac{41}{42},-\frac{1}{2}}(21\tau)\right),
\end{split}
\end{equation*}
so
\begin{equation*}
\frac{\zeta_{14}^{-1}}{2\eta(\tau)}\ \theta_{\frac{1}{14}e,\frac{1}{14}e}(\tau) = q^{-\frac{1}{168}} \ef_0(q)+\zeta_{84}^{-13} R_{\frac{13}{42},-\frac{1}{2}}(21\tau) +\zeta_{84} R_{\frac{41}{42},-\frac{1}{2}}(21\tau).
\end{equation*}
Similarly, we find
\begin{equation*}
\begin{split}
\frac{\zeta_{14}^{-1}}{2\eta(\tau)}\ \theta_{\frac{3}{14}e,\frac{5}{14}e}(\tau) &= q^{\frac{47}{168}} \ef_2(q)+\zeta_{84}^{29} R_{\frac{11}{42},-\frac{5}{2}}(21\tau) +\zeta_{84}^{-41} R_{\frac{25}{42},-\frac{5}{2}}(21\tau),\\
\frac{\zeta_{14}^{-1}}{2\eta(\tau)}\ \theta_{\frac{5}{14}e,\frac{3}{14}e}(\tau) &= q^{-\frac{25}{168}} \ef_1(q)+\zeta_{28}^5 R_{\frac{23}{42},-\frac{3}{2}}(21\tau) +\zeta_{28}^{-9} R_{\frac{37}{42},-\frac{3}{2}}(21\tau).
\end{split}
\end{equation*}
So
\begin{equation*}
H_7 = F_7 -G_7.
\end{equation*}
Using (5) of Proposition \ref{propRab} and the fact that $F_7$ is holomorphic, we find
\begin{align*}
\Omega_{\frac{1}{2}} G_7&=\frac{3}{16} G_7\\
\Omega_{\frac{1}{2}} F_7&=\frac{3}{16} F_7,\\
\intertext{and so}
\Omega_{\frac{1}{2}} H_7&=\frac{3}{16} H_7.
\end{align*}
Hence $H_7$ is a vector valued real-analytic modular form of weight 1/2. From (3) of Proposition \ref{propRab} we obtain that $G_7$ is bounded if $\tau\downarrow \xi$, with $\xi\in\Q$. \qed

As a corollary we get the description of the non-modularity of the ``seventh order'' mock $\theta$-function $F_7$. To state the corollary we need the following vector of theta functions of weight $3/2$:
\begin{equation*}
g_7 (\tau) := \begin{pmatrix} \zeta_{84}^{-13}\ g_{\frac{13}{42},\frac{1}{2}} +\zeta_{84}\ g_{\frac{41}{42},\frac{1}{2}}\\ 
\zeta_{84}^{73}\ g_{\frac{11}{42},\frac{1}{2}} +\zeta_{84}^{59}\ g_{\frac{25}{42},\frac{1}{2}}\\
\zeta_{84}^{61}\ g_{\frac{23}{42},\frac{1}{2}} +\zeta_{84}^{47}\ g_{\frac{37}{42},\frac{1}{2}} \end{pmatrix}(21\tau).
\end{equation*}
This function has the following modular transformation property, which can be verified using standard methods:
\begin{equation*}
g_7\Bigl( -\frac{1}{\tau} \Bigr) = -M_7 (-i\tau)^{3/2} g_7 (\tau),
\end{equation*}
with 
\begin{equation*}
M_7 := \frac{2}{\sqrt{7}} \begin{pmatrix} \sin \frac{\pi}{7}&\sin \frac{3\pi}{7}&\sin \frac{2\pi}{7}\\ \sin \frac{3\pi}{7}&-\sin \frac{2\pi}{7}&\sin \frac{\pi}{7}\\ \sin \frac{2\pi}{7}&\sin \frac{\pi}{7}&-\sin \frac{3\pi}{7}\end{pmatrix}.
\end{equation*}

\begin{corollary}\label{cor7}
We have
\begin{equation*}
F_7 (\tau) - \frac{1}{\sqrt{-i\tau}} M_7 F_7 \Bigl( -\frac{1}{\tau} \Bigr) = i\sqrt{21} \int_0^{i\infty} \frac{g_7(z)}{\sqrt{-i(z+\tau)}}\ dz,
\end{equation*}
where we have to integrate each component of the vector,
as well as the obvious equation
\begin{equation*}
F_7(\tau+1) = \begin{pmatrix} \zeta_{168}^{-1}&0&0\\ 0&\zeta_{168}^{47}&0\\0&0&\zeta_{168}^{-25}\end{pmatrix} F_7(\tau).
\end{equation*}
\end{corollary}
\proof According to Proposition \ref{propm7} we have
\begin{equation*}
H_7\Bigl( -\frac{1}{\tau}\Bigr) = \sqrt{-i\tau} M_7 H_7 (\tau).
\end{equation*}
If we replace $\tau$ by $-1/\tau$ in the equation (or multiply both sides by $M_7/\sqrt{-i\tau}$ and use $M_7^2=I$) we find
\begin{equation*}
H_7 (\tau) = \frac{1}{\sqrt{-i\tau}} M_7 H_7 \Bigl( -\frac{1}{\tau} \Bigr),
\end{equation*}
so
\begin{equation}\label{transjeF}
F_7 (\tau) - \frac{1}{\sqrt{-i\tau}} M_7 F_7 \Bigl( -\frac{1}{\tau} \Bigr)=G_7 (\tau) - \frac{1}{\sqrt{-i\tau}} M_7 G_7 \Bigl( -\frac{1}{\tau} \Bigr).
\end{equation}
Using (2) of Proposition \ref{propRab} and (2) of Proposition \ref{prop141} we see:
\begin{equation}\label{g71}
G_7(\tau)= i\sqrt{21} \int_{-\overline{\tau}}^{i\infty} \frac{g_7(z)}{\sqrt{-i(z+\tau)}}\ dz.
\end{equation}
Hence
\begin{equation*}
\begin{split}
\frac{1}{\sqrt{-i\tau}} G_7\Bigl( -\frac{1}{\tau} \Bigr) &= \frac{i\sqrt{21}}{\sqrt{-i\tau}} \int_{1/\overline{\tau}}^{i\infty} \frac{g_7(z)}{\sqrt{-i(z-1/\tau)}}\ dz\\
&= i\sqrt{21} \int_{0}^{-\overline{\tau}} \frac{g_7(-1/z)}{\sqrt{1+\tau/z}}\frac{dz}{(-iz)^2},
\end{split}
\end{equation*}
where we have replaced $z$ by $-1/z$ in the integral. Using the transformation property of $g_7$, we find
\begin{equation}\label{g72}
\frac{1}{\sqrt{-i\tau}} G_7\Bigl( -\frac{1}{\tau} \Bigr)= -i\sqrt{21} M_7 \int_0^{-\overline{\tau}} \frac{g_7(z)}{\sqrt{-i(z+\tau)}}\ dz.
\end{equation}
Putting \eqref{g71} and \eqref{g72} in \eqref{transjeF} we get the desired result. \qed

\begin{remark} Using (2) of Theorem \ref{period} we could give the non-modularity of $F_7$ in terms of the function $h$ from Chapter 1. The result is similar to results found by Watson in \cite{watson1} for the ``third order'' mock $\theta$-functions.
\end{remark}

\section{The fifth order mock $\theta$-functions}

In this section we deal with eight of the ten ``fifth order'' mock $\theta$-functions from Ramanujan's letter. The remaining two will be discussed in the next section. In \cite{andrews3} we find the following identities: 
{\allowdisplaybreaks
\begin{align*} 
f_0(q) &=\frac{1}{(q)_\infty} \underset{|j|\leq n}{\sum_{n=0}^\infty} (-1)^j q^{\frac{5}{2}n^2 +\frac{1}{2}n -j^2}(1-q^{4n+2}),\\
F_0(q)&= \frac{1}{(q^2;q^2)_\infty} \sum_{n=0}^\infty \sum_{j=0}^{2n} (-1)^n q^{5n^2+2n-\frac{1}{2}j^2-\frac{1}{2}j}(1+q^{6n+3}),\\
1+2\psi_0(q)&=\frac{(-q)_\infty}{(q)_\infty} \Bigl( 1+2\sum_{n=1}^\infty (-1)^n q^{n^2+n}-2\underset{|j|<n}{\sum_{n=1}^\infty} (-1)^j q^{\frac{5}{2}n^2-\frac{1}{2}n-\frac{3}{2}j^2-\frac{1}{2}j}(1-q^n)\Bigr),\\
\phi_0(q)&= \frac{(-q;q^2)_\infty}{(q^2;q^2)_\infty} \underset{|j|\leq n}{\sum_{n=0}^\infty} (-1)^j q^{5n^2+2n-3j^2-j}(1-q^{6n+3}),\\
f_1(q)&= \frac{1}{(q)_\infty} \underset{|j|\leq n}{\sum_{n=0}^\infty} (-1)^j q^{\frac{5}{2}n^2 +\frac{3}{2}n -j^2}(1-q^{2n+1}),\\
F_1(q)&= \frac{1}{(q^2;q^2)_\infty} \sum_{n=0}^\infty \sum_{j=0}^{2n} (-1)^n q^{5n^2+4n-\frac{1}{2}j^2-\frac{1}{2}j}(1+q^{2n+1}),\\
\psi_1(q)&=\frac{(-q)_\infty}{(q)_\infty} \underset{|j|\leq n}{\sum_{n=0}^\infty} (-1)^j q^{\frac{5}{2}n^2+\frac{3}{2}n-\frac{3}{2}j^2-\frac{1}{2}j}(1-q^{2n+1}),\\
\phi_1(q)&= q \frac{(-q;q^2)_\infty}{(q^2;q^2)_\infty} \underset{|j|\leq n}{\sum_{n=0}^\infty} (-1)^j q^{5n^2+4n-3j^2-j}(1-q^{2n+1}).
\end{align*}}
Note that there are mistakes in the 3rd and 8th formula in \cite{andrews3}.
We will use these identities as the definitions of the mock $\theta$-functions. We write four of these identities in a more suitable form (the other four  will be discussed later):
\begin{lemma}\label{le1}
We have
\begin{align*}
2\eta(\tau) &q^{-\frac{1}{60}} f_0(q) \\
&= \sum_{\nu\in\left(\begin{smallmatrix} 1/10\\0\end{smallmatrix}\right)+\Z^2} \Bigl\{\sign(B(\nu,c_1))-\sign(B(\nu,c_2))\Bigr\}\ e^{2\pi iQ(\nu)\tau+2\pi i B\left(\nu,\left(\begin{smallmatrix} 0\\1/4\end{smallmatrix}\right)\right)}\displaybreak[0]\\ 
2\eta(\tau) &q^{\frac{11}{60}} f_1(q) \\
&= \sum_{\nu\in\left(\begin{smallmatrix} 3/10\\0\end{smallmatrix}\right)+\Z^2} \Bigl\{\sign(B(\nu,c_1))-\sign(B(\nu,c_2))\Bigr\}\ e^{2\pi iQ(\nu)\tau+2\pi i B\left(\nu,\left(\begin{smallmatrix} 0\\1/4\end{smallmatrix}\right)\right)}\displaybreak[0]\\
2\eta(\tau) &q^{-\frac{1}{240}} (-1+F_0(q^{\frac{1}{2}})) \\
&= \sum_{\nu\in\left(\begin{smallmatrix} 1/5\\1/4\end{smallmatrix}\right)+\Z^2} \Bigl\{\sign(B(\nu,c_1))-\sign(B(\nu,c_2))\Bigr\}\ e^{2\pi iQ(\nu)\tau+2\pi i B\left(\nu,\left(\begin{smallmatrix} 1/2\\1\end{smallmatrix}\right)\right)}\displaybreak[0]\\
2\eta(\tau) &q^{\frac{71}{240}} F_1(q^{\frac{1}{2}}) \\
&= \sum_{\nu\in\left(\begin{smallmatrix} 2/5\\1/4\end{smallmatrix}\right)+\Z^2} \Bigl\{\sign(B(\nu,c_1))-\sign(B(\nu,c_2))\Bigr\}\ e^{2\pi iQ(\nu)\tau+2\pi i B\left(\nu,\left(\begin{smallmatrix} 1/2\\2\end{smallmatrix}\right)\right)}\displaybreak[0]\\
2\eta(\tau) \zeta_8^{-1} &q^{-\frac{1}{240}} (-1+F_0(-q^{\frac{1}{2}})) \\
&= \sum_{\nu\in\left(\begin{smallmatrix} 1/5\\1/4\end{smallmatrix}\right)+\Z^2} \Bigl\{\sign(B(\nu,c_1))-\sign(B(\nu,c_2))\Bigr\}\ e^{2\pi iQ(\nu)\tau+2\pi i B\left(\nu,\left(\begin{smallmatrix} 0\\1/4\end{smallmatrix}\right)\right)}\displaybreak[0]\\
2\eta(\tau) \zeta_8^{-1} &q^{\frac{71}{240}} F_1(-q^{\frac{1}{2}}) \\
&= \sum_{\nu\in\left(\begin{smallmatrix} 2/5\\1/4\end{smallmatrix}\right)+\Z^2} \Bigl\{\sign(B(\nu,c_1))-\sign(B(\nu,c_2))\Bigr\}\ e^{2\pi iQ(\nu)\tau+2\pi i B\left(\nu,\left(\begin{smallmatrix} 0\\1/4\end{smallmatrix}\right)\right)},
\end{align*}
with $A=\left(\begin{smallmatrix} 5&0\\0&-2 \end{smallmatrix}\right)$, $c_1=\left(\begin{smallmatrix} 2\\5 \end{smallmatrix}\right)$ and  $c_2=\left(\begin{smallmatrix} -2\\5 \end{smallmatrix}\right)$.
\end{lemma}
\begin{remark}
We have $B(c_1,c_2)=-70$ and $Q(c_1)=Q(c_2)=-15$. If we choose $C_Q$ such that $c_1\in C_Q$ then also $c_2\in C_Q$.
\end{remark}
\proof We have 
{\allowdisplaybreaks
\begin{align*}
\underset{|j|\leq n}{\sum_{n=0}^\infty} (-1)^j q^{\frac{5}{2}n^2 +\frac{1}{2}n -j^2}&(1-q^{4n+2})\\
&= \underset{|j|\leq n}{\sum_{n=0}^\infty} (-1)^j q^{\frac{5}{2}n^2 +\frac{1}{2}n -j^2} - \underset{|j|\leq n}{\sum_{n=0}^\infty} (-1)^j q^{\frac{5}{2}n^2 +\frac{9}{2}n +2 -j^2}\\
&= \Bigl(\sum_{n+j\geq0, n-j\geq0} -\sum_{n+j<0, n-j<0}\Bigr) (-1)^j q^{\frac{5}{2}n^2+\frac{1}{2}n -j^2},
\end{align*}}
where we have replaced $n$ by $-n-1$ in the second sum. From this we get the first identity. The proof of the 2rd identity is similar. 

Using
\begin{equation}\label{fo}
\sum_{j=0}^{2n} q^{-\frac{1}{2}j^2-\frac{1}{2}j} = \sum_{j=-n}^n q^{-2j^2-j},
\end{equation}
we see
\begin{equation*}
\begin{split}
\sum_{n=0}^\infty &\sum_{j=0}^{2n} (-1)^n q^{5n^2+2n-\frac{1}{2}j^2-\frac{1}{2}j}(1+q^{6n+3})\\
&= \sum_{n=0}^\infty \sum_{j=-n}^{n} (-1)^n q^{5n^2+2n-2j^2-j}(1+q^{6n+3})\\
&= \Bigl(\sum_{n+j\geq0, n-j\geq0} -\sum_{n+j<0, n-j<0}\Bigr) (-1)^n q^{5n^2+2n-2j^2-j}\\
&= (q^2;q^2)_\infty + \Bigl(\sum_{n+j\geq0, n-j>0} -\sum_{n+j<0, n-j\leq 0}\Bigr) (-1)^n q^{5n^2+2n-2j^2-j}.
\end{split}
\end{equation*}
From this we get the 3rd and 5th identity. Again using \eqref{fo} we see
\begin{equation*}
\begin{split}
\sum_{n=0}^\infty &\sum_{j=0}^{2n} (-1)^n q^{5n^2+4n-\frac{1}{2}j^2-\frac{1}{2}j}(1+q^{2n+1})\\
&= \sum_{n=0}^\infty \sum_{j=-n}^{n} (-1)^n q^{5n^2+4n-2j^2-j}(1+q^{2n+1})\\
&= \Bigl(\sum_{n+j\geq0, n-j\geq0} -\sum_{n+j<0, n-j<0}\Bigr) (-1)^n q^{5n^2+4n-2j^2-j}.
\end{split}
\end{equation*}
From this we get the 4th and 6th identity. \qed

We collect these six mock $\theta$-functions into a single vector-valued mock $\theta$-function
\begin{align*}
F_{5,1} (\tau) &:= \begin{pmatrix}
q^{-\frac{1}{60}} f_0(q)\\
q^{\frac{11}{60}} f_1(q)\\
q^{-\frac{1}{240}} (-1+F_0 (q^{\frac{1}{2}}))\\
q^{\frac{71}{240}} F_1(q^{\frac{1}{2}})\\
q^{-\frac{1}{240}} (-1+F_0 (-q^{\frac{1}{2}}))\\
q^{\frac{71}{240}} F_1(-q^{\frac{1}{2}}))
\end{pmatrix}.\\
\intertext{To express its modular transformation behaviour, we also introduce}
H_{5,1} (\tau) &:= \frac{1}{2\eta(\tau)} \begin{pmatrix} \theta_{\left(\begin{smallmatrix} 1/10\\0\end{smallmatrix}\right), \left(\begin{smallmatrix} 0\\1/4\end{smallmatrix}\right)}\\
\theta_{\left(\begin{smallmatrix} 3/10\\0\end{smallmatrix}\right), \left(\begin{smallmatrix} 0\\1/4\end{smallmatrix}\right)}\\
\theta_{\left(\begin{smallmatrix} 1/5\\1/4\end{smallmatrix}\right), \left(\begin{smallmatrix} 1/2\\1\end{smallmatrix}\right)}\\
\theta_{\left(\begin{smallmatrix} 2/5\\1/4\end{smallmatrix}\right), \left(\begin{smallmatrix} 1/2\\2\end{smallmatrix}\right)}\\
\zeta_8 \theta_{\left(\begin{smallmatrix} 1/5\\1/4\end{smallmatrix}\right), \left(\begin{smallmatrix} 0\\1/4\end{smallmatrix}\right)}\\
\zeta_8 \theta_{\left(\begin{smallmatrix} 2/5\\1/4\end{smallmatrix}\right), \left(\begin{smallmatrix} 0\\1/4\end{smallmatrix}\right)}
\end{pmatrix} (\tau),\\
\intertext{with $A$, $c_1$ and $c_2$ as in Lemma \ref{le1}, and}
G_{5,1}(\tau)&:= \frac{1}{2}\begin{pmatrix}
2\zeta_{12} R_{\frac{1}{30},\frac{5}{2}}+2\zeta_{12}^{-1} R_{\frac{11}{30},\frac{5}{2}}\\
2\zeta_{12} R_{\frac{13}{30},\frac{5}{2}}+2\zeta_{12}^{-1} R_{\frac{23}{30},\frac{5}{2}}\\
-R_{\frac{19}{60},0}-R_{\frac{29}{60},0} + R_{\frac{49}{60},0}+R_{\frac{59}{60},0}\\
-R_{\frac{13}{60},0}-R_{\frac{23}{60},0} + R_{\frac{43}{60},0} +R_{\frac{53}{60},0}\\
\zeta_{24}^{-5} R_{\frac{19}{60},\frac{5}{2}}+\zeta_{24}^5 R_{\frac{29}{60},\frac{5}{2}} +\zeta_{24} R_{\frac{49}{60},\frac{5}{2}}+\zeta_{24}^{-1} R_{\frac{59}{60},\frac{5}{2}}\\
\zeta_{24} R_{\frac{13}{60},\frac{5}{2}}+\zeta_{24}^{-1} R_{\frac{23}{60},\frac{5}{2}}+\zeta_{24}^{-5} R_{\frac{43}{60},\frac{5}{2}} +\zeta_{24}^5 R_{\frac{53}{60},\frac{5}{2}}
\end{pmatrix}(30\tau).
\end{align*}

\begin{proposition}\label{51}
We have
\begin{equation*}
F_{5,1}=H_{5,1}+G_{5,1},
\end{equation*}
where
\begin{description}
\item[(1)] The function $H_{5,1}$ is a (vector-valued) real-analytic modular form of weight 1/2, satisfying
\begin{align*}
H_{5,1}(\tau+1) &= \begin{pmatrix} \zeta_{60}^{-1}&0&0&0&0&0\\
0&\zeta_{60}^{11}&0&0&0&0\\
0&0&0&0&\zeta_{240}^{-1}&0\\
0&0&0&0&0&\zeta_{240}^{71}\\
0&0&\zeta_{240}^{-1}&0&0&0\\
0&0&0&\zeta_{240}^{71}&0&0
\end{pmatrix} H_{5,1}(\tau),\\
\intertext{and}
H_{5,1}\Bigl(-\frac{1}{\tau}\Bigr) &=\sqrt{-i\tau} \frac{2}{\sqrt{5}}\ M_5\
 H_{5,1}(\tau),
\end{align*}
with
\begin{equation*}
M_5=\begin{pmatrix}
0&0&\sqrt{2}\sin\frac{\pi}{5}&\sqrt{2}\sin\frac{2\pi}{5}&0&0\\
0&0&\sqrt{2}\sin\frac{2\pi}{5}&-\sqrt{2}\sin\frac{\pi}{5}&0&0\\
\frac{1}{\sqrt{2}}\sin\frac{\pi}{5}&\frac{1}{\sqrt{2}}\sin\frac{2\pi}{5}&0&0&0&0\\
\frac{1}{\sqrt{2}}\sin\frac{2\pi}{5}&-\frac{1}{\sqrt{2}}\sin\frac{\pi}{5}&0&0&0&0\\
0&0&0&0&\sin\frac{2\pi}{5}&\sin\frac{\pi}{5}\\
0&0&0&0&\sin\frac{\pi}{5}&-\sin\frac{2\pi}{5}
\end{pmatrix},
\end{equation*}
and is an eigenfunction of the Casimir operator $\Omega_{\frac{1}{2}}$, with eigenvalue $\frac{3}{16}$.
\item[(2)] The function $G_{5,1}$ is bounded if $\tau\downarrow \xi$, with $\xi\in\Q$.
\end{description}
\end{proposition}
\proof
We consider the function
\begin{equation*}
\Theta (\tau) = \begin{pmatrix} \theta_{\left(\begin{smallmatrix} 1/10\\0\end{smallmatrix}\right), \left(\begin{smallmatrix} 0\\1/4\end{smallmatrix}\right)}\\
\theta_{\left(\begin{smallmatrix} 3/10\\0\end{smallmatrix}\right), \left(\begin{smallmatrix} 0\\1/4\end{smallmatrix}\right)}\\
\theta_{\left(\begin{smallmatrix} 1/5\\1/4\end{smallmatrix}\right), \left(\begin{smallmatrix} 1/2\\1\end{smallmatrix}\right)}\\
\theta_{\left(\begin{smallmatrix} 2/5\\1/4\end{smallmatrix}\right), \left(\begin{smallmatrix} 1/2\\2\end{smallmatrix}\right)}\\
\theta_{\left(\begin{smallmatrix} 1/5\\1/4\end{smallmatrix}\right), \left(\begin{smallmatrix} 0\\1/4\end{smallmatrix}\right)}\\
\theta_{\left(\begin{smallmatrix} 2/5\\1/4\end{smallmatrix}\right), \left(\begin{smallmatrix} 0\\1/4\end{smallmatrix}\right)}
\end{pmatrix}(\tau).
\end{equation*}
Using (4) and (2) of Corollary \ref{cor1} we see
\begin{equation*}
\Theta (\tau+1) = \begin{pmatrix} \zeta_{40}&0&0&0&0&0\\
0&\zeta_{40}^9&0&0&0&0\\
0&0&0&0&\zeta_{80}^{13}&0\\
0&0&0&0&0&\zeta_{80}^{37}\\
0&0&\zeta_{80}^{-7}&0&0&0\\
0&0&0&\zeta_{80}^{17}&0&0
\end{pmatrix} \Theta (\tau).
\end{equation*}
Let $C=\left(\begin{smallmatrix} -1&0\\0&1\end{smallmatrix}\right)$, then $C\in O_A^+(\Z)$, $Cc_1=c_2$ and $Cc_2=c_1$. Hence, using Corollary \ref{cor2}, we find
\begin{equation*}
\theta_{a,b}^{c_1,c_2} = \theta_{Ca,Cb}^{c_2,c_1}=-\theta_{Ca,Cb}^{c_1,c_2}
\end{equation*}
Using this and Corollary \ref{cor1}, we obtain
\begin{equation*}
\Theta \Bigl(-\frac{1}{\tau}\Bigr) = -i\tau \frac{2}{\sqrt{5}}\ M_5\ \Theta (\tau).
\end{equation*}
If we use $\eta(\tau+1) =\zeta_{24} \eta(\tau)$ and $\eta\left(-\frac{1}{\tau}\right)= \sqrt{-i\tau}\ \eta(\tau)$, we obtain the transformations for $H_{5,1}$.

If we write $\rho(\nu;\tau)$ as the sum of the three expressions \eqref{al1}, \eqref{al2} and \eqref{al3} and use Proposition \ref{split}, we find
\begin{equation*}
\theta_{\left(\begin{smallmatrix} 1/10\\0\end{smallmatrix}\right), \left(\begin{smallmatrix} 0\\1/4\end{smallmatrix}\right)}(\tau)=2\eta(\tau) q^{-\frac{1}{60}} f_0(q) -2\eta(\tau) \left( \zeta_{12} R_{\frac{1}{30},\frac{5}{2}}(30\tau) +\zeta_{12}^{-1} R_{\frac{11}{30},\frac{5}{2}}(30\tau)\right).
\end{equation*}
Hence
\begin{equation*}
\frac{1}{2\eta(\tau)}\ \theta_{\left(\begin{smallmatrix} 1/10\\0\end{smallmatrix}\right), \left(\begin{smallmatrix} 0\\1/4\end{smallmatrix}\right)}(\tau)= q^{-\frac{1}{60}} f_0(q) -\left( \zeta_{12} R_{\frac{1}{30},\frac{5}{2}}(30\tau) +\zeta_{12}^{-1} R_{\frac{11}{30},\frac{5}{2}}(30\tau)\right).
\end{equation*}
We can find similar identities for the other components of $H_{5,1}$. Combining them gives
\begin{equation*}
H_{5,1} = F_{5,1} -G_{5,1}.
\end{equation*}
Using (5) of Proposition \ref{propRab} and the fact that $F_{5,1}$ is holomorphic, we find
\begin{align*}
\Omega_{\frac{1}{2}} G_{5,1}&=\frac{3}{16} G_{5,1}\\
\Omega_{\frac{1}{2}} F_{5,1}&=\frac{3}{16} F_{5,1},\\
\intertext{and so}
\Omega_{\frac{1}{2}} H_{5,1}&=\frac{3}{16} H_{5,1}.
\end{align*}
Hence $H_{5,1}$ is a vector valued real-analytic modular form of weight 1/2. From (3) of Proposition \ref{propRab} we get that $G_{5,1}$ is bounded as $\tau\downarrow \xi$, with $\xi\in\Q$. \qed

As with Corollary \ref{cor7}, we could also use Proposition \ref{51} to describe the non-modularity of the ``fifth order'' mock $\theta$-function $F_{5,1}$. We omit this.
 
We turn now to the other four ``fifth order'' mock $\theta$-functions.

\begin{lemma}\label{le2}
We have
\begin{align*}
2\zeta_{12}^{-1} &\frac{\eta(\tau)^2}{\eta(2\tau)}\ q^{-\frac{1}{60}} \psi_0(q) \\
&= \sum_{\nu\in\left(\begin{smallmatrix} 1/10\\1/6\end{smallmatrix}\right)+\Z^2} \Bigl\{\sign(B(\nu,c_1))-\sign(B(\nu,c_2))\Bigr\}\ e^{2\pi iQ(\nu)\tau+2\pi i B\left(\nu,\left(\begin{smallmatrix} 0\\1/6\end{smallmatrix}\right)\right)}\displaybreak[0]\\ 
2\zeta_{12}^{-1} &\frac{\eta(\tau)^2}{\eta(2\tau)}\ q^{\frac{11}{60}} \psi_1(q) \\
&= \sum_{\nu\in\left(\begin{smallmatrix} 3/10\\1/6\end{smallmatrix}\right)+\Z^2} \Bigl\{\sign(B(\nu,c_1))-\sign(B(\nu,c_2))\Bigr\}\ e^{2\pi iQ(\nu)\tau+2\pi i B\left(\nu,\left(\begin{smallmatrix} 0\\1/6\end{smallmatrix}\right)\right)}\displaybreak[0]\\
2\zeta_{60} &\frac{\eta(\tau)^2}{\eta(\tau/2)}\ q^{-\frac{1}{240}} \phi_0(-q^{\frac{1}{2}}) \\
&= \sum_{\nu\in\left(\begin{smallmatrix} 1/5\\1/6\end{smallmatrix}\right)+\Z^2} \Bigl\{\sign(B(\nu,c_1))-\sign(B(\nu,c_2))\Bigr\}\ e^{2\pi iQ(\nu)\tau+2\pi i B\left(\nu,\left(\begin{smallmatrix} 1/10\\1/6\end{smallmatrix}\right)\right)}\displaybreak[0]\\
-2\zeta_{60}^7 &\frac{\eta(\tau)^2}{\eta(\tau/2)}\ q^{-\frac{49}{240}} \phi_1(-q^{\frac{1}{2}}) \\
&= \sum_{\nu\in\left(\begin{smallmatrix} 2/5\\1/6\end{smallmatrix}\right)+\Z^2} \Bigl\{\sign(B(\nu,c_1))-\sign(B(\nu,c_2))\Bigr\}\ e^{2\pi iQ(\nu)\tau+2\pi i B\left(\nu,\left(\begin{smallmatrix} 1/10\\1/6\end{smallmatrix}\right)\right)}\displaybreak[0]\\
2\zeta_{16}^{-1} &\frac{\eta(\tau)^2}{\eta((\tau+1)/2)}\ q^{-\frac{1}{240}} \phi_0(q^{\frac{1}{2}}) \\
&= \sum_{\nu\in\left(\begin{smallmatrix} 1/5\\1/6\end{smallmatrix}\right)+\Z^2} \Bigl\{\sign(B(\nu,c_1))-\sign(B(\nu,c_2))\Bigr\}\ e^{2\pi iQ(\nu)\tau+2\pi i B\left(\nu,\left(\begin{smallmatrix} 0\\1/6\end{smallmatrix}\right)\right)}\displaybreak[0]\\
2\zeta_{16}^{-1} &\frac{\eta(\tau)^2}{\eta((\tau+1)/2)}\ q^{-\frac{49}{240}} \phi_1(q^{\frac{1}{2}}) \\
&= \sum_{\nu\in\left(\begin{smallmatrix} 2/5\\1/6\end{smallmatrix}\right)+\Z^2} \Bigl\{\sign(B(\nu,c_1))-\sign(B(\nu,c_2))\Bigr\}\ e^{2\pi iQ(\nu)\tau+2\pi i B\left(\nu,\left(\begin{smallmatrix} 0\\1/6\end{smallmatrix}\right)\right)},
\end{align*}
with $A=\left(\begin{smallmatrix} 5&0\\0&-3 \end{smallmatrix}\right)$, $c_1=\left(\begin{smallmatrix} 3\\5 \end{smallmatrix}\right)$ and  $c_2=\left(\begin{smallmatrix} -3\\5 \end{smallmatrix}\right)$.
\end{lemma}
\begin{remark}
We have $B(c_1,c_2)=-120$ and $Q(c_1)=Q(c_2)=-15$. If we choose $C_Q$ such that $c_1\in C_Q$ then also $c_2\in C_Q$.
\end{remark}
\proof
The proof is similar to the proof of Lemma \ref{le1}. We also have to use
\begin{align*}
\frac{(q^{\frac{1}{2}})_\infty}{(q)_\infty}&= \frac{(q^{\frac{1}{2}};q^{\frac{1}{2}})_\infty}{(q)_\infty^2}= q^{\frac{1}{16}} \frac{\eta(\tau/2)}{\eta(\tau)^2}\\
\frac{(-q^{\frac{1}{2}})_\infty}{(q)_\infty}&=\zeta_{48}^{-1} q^{\frac{1}{16}} \frac{\eta\bigl((\tau+1)/2\bigr)}{\eta(\tau)^2},\\
\intertext{and}
\frac{(-q)_\infty}{(q)_\infty}&= \frac{(q^2;q^2)_\infty}{(q)_\infty^2}=\frac{\eta(2\tau)}{\eta(\tau)^2}.
\end{align*}
Only the first identity is a bit more difficult: We have
\begin{equation*}
\begin{split}
-\underset{|j|< n}{\sum_{n=1}^\infty} (-1)^j &q^{\frac{5}{2}n^2 -\frac{1}{2}n -\frac{3}{2}j^2-\frac{1}{2}j}(1-q^n)\\
&=\underset{|j|< n}{\sum_{n=1}^\infty} (-1)^j q^{\frac{5}{2}n^2 +\frac{1}{2}n -\frac{3}{2}j^2-\frac{1}{2}j} -\underset{|j|< n}{\sum_{n=1}^\infty} (-1)^j q^{\frac{5}{2}n^2 -\frac{1}{2}n -\frac{3}{2}j^2-\frac{1}{2}j}\\
&= \Bigl(\sum_{n+j>0, n-j>0} -\sum_{n+j<0, n-j<0}\Bigr) (-1)^j q^{\frac{5}{2}n^2+\frac{1}{2}n -\frac{3}{2}j^2-\frac{1}{2}j},
\end{split}
\end{equation*}
where we have replaced $n$ by $-n$ in the second sum. Hence \begin{equation}\label{fietje}
\begin{split}
1+&2\sum_{n=1}^\infty (-1)^n q^{n^2+n}-2\underset{|j|<n}{\sum_{n=1}^\infty} (-1)^j q^{\frac{5}{2}n^2-\frac{1}{2}n-\frac{3}{2}j^2-\frac{1}{2}j}(1-q^n)\\
&= 2 \Bigl(\sum_{n+j\geq 0, n-j>0} -\sum_{n+j<0, n-j\leq 0}\Bigr) (-1)^j q^{\frac{5}{2}n^2+\frac{1}{2}n -\frac{3}{2}j^2-\frac{1}{2}j} +2\sum_{n=-\infty}^{-1} (-1)^n q^{n^2}+1.
\end{split}
\end{equation}
We have 
\begin{equation*}
2\sum_{n=-\infty}^{-1} (-1)^n q^{n^2}+1= \sum_{n\in\Z}(-1)^n q^{n^2}=\frac{(q)_\infty}{(-q)_\infty}.
\end{equation*}
Using this and \eqref{fietje} we find 
\begin{equation*}
\psi_0(q) = \frac{(-q)_\infty}{(q)_\infty}\ \Bigl(\sum_{n+j\geq0, n-j>0} -\sum_{n+j<0, n-j\leq 0}\Bigr) (-1)^j q^{\frac{5}{2}n^2+\frac{1}{2}n -\frac{3}{2}j^2-\frac{1}{2}j},
\end{equation*}
from which the result follows. \qed

We collect these six mock $\theta$-functions into a single vector-valued mock $\theta$-function
\begin{align*}
F_{5,2} (\tau) &:= \begin{pmatrix} 2q^{-\frac{1}{60}} \psi_0 (q)\\
2q^{\frac{11}{60}} \psi_1 (q)\\
q^{-\frac{1}{240}} \phi_0(-q^{\frac{1}{2}})\\
-q^{-\frac{49}{240}} \phi_1(-q^{\frac{1}{2}})\\
q^{-\frac{1}{240}} \phi_0(q^{\frac{1}{2}})\\
q^{-\frac{49}{240}} \phi_1(q^{\frac{1}{2}})
\end{pmatrix}\\
\intertext{To express its modular transformation behaviour, we also introduce}
H_{5,2} (\tau) &:= \frac{1}{2\eta(\tau)^2} \begin{pmatrix} 2\zeta_{12}\ \eta(2\tau)\ \theta_{\left(\begin{smallmatrix} 1/10\\1/6\end{smallmatrix}\right), \left(\begin{smallmatrix} 0\\1/6\end{smallmatrix}\right)}(\tau)\\
2\zeta_{12}\ \eta(2\tau)\ \theta_{\left(\begin{smallmatrix} 3/10\\1/6\end{smallmatrix}\right), \left(\begin{smallmatrix} 0\\1/6\end{smallmatrix}\right)}(\tau)\\
\zeta_{60}^{-1} \eta(\tau/2)\ \theta_{\left(\begin{smallmatrix} 1/5\\1/6\end{smallmatrix}\right), \left(\begin{smallmatrix} 1/10\\1/6\end{smallmatrix}\right)}(\tau)\\
\zeta_{60}^{-7} \eta(\tau/2)\ \theta_{\left(\begin{smallmatrix} 2/5\\1/6\end{smallmatrix}\right), \left(\begin{smallmatrix} 1/10\\1/6\end{smallmatrix}\right)}(\tau)\\
\zeta_{16}\ \eta\bigl((\tau+1)/2\bigr)\ \theta_{\left(\begin{smallmatrix} 1/5\\1/6\end{smallmatrix}\right), \left(\begin{smallmatrix} 0\\1/6\end{smallmatrix}\right)}(\tau)\\
\zeta_{16}\ \eta\bigl((\tau+1)/2\bigr)\ \theta_{\left(\begin{smallmatrix} 2/5\\1/6\end{smallmatrix}\right), \left(\begin{smallmatrix} 0\\1/6\end{smallmatrix}\right)}(\tau)
\end{pmatrix},\\
\intertext{with $A$, $c_1$ and $c_2$ as in Lemma \ref{le2}, and}
G_{5,2}(\tau)&:= -\frac{1}{2}\begin{pmatrix}
2\zeta_{12} R_{\frac{1}{30},\frac{5}{2}}+2\zeta_{12}^{-1} R_{\frac{11}{30},\frac{5}{2}}\\
2\zeta_{12} R_{\frac{13}{30},\frac{5}{2}}+2\zeta_{12}^{-1} R_{\frac{23}{30},\frac{5}{2}}\\
-R_{\frac{19}{60},0}-R_{\frac{29}{60},0} + R_{\frac{49}{60},0}+R_{\frac{59}{60},0}\\
-R_{\frac{13}{60},0}-R_{\frac{23}{60},0} + R_{\frac{43}{60},0} +R_{\frac{53}{60},0}\\
\zeta_{24}^{-5} R_{\frac{19}{60},\frac{5}{2}}+\zeta_{24}^5 R_{\frac{29}{60},\frac{5}{2}} +\zeta_{24} R_{\frac{49}{60},\frac{5}{2}}+\zeta_{24}^{-1} R_{\frac{59}{60},\frac{5}{2}}\\
\zeta_{24} R_{\frac{13}{60},\frac{5}{2}}+\zeta_{24}^{-1} R_{\frac{23}{60},\frac{5}{2}}+\zeta_{24}^{-5} R_{\frac{43}{60},\frac{5}{2}} +\zeta_{24}^5 R_{\frac{53}{60},\frac{5}{2}}
\end{pmatrix}(30\tau). 
\end{align*}

\begin{proposition}\label{52}
We have
\begin{equation*}
F_{5,2}=H_{5,2}+G_{5,2},
\end{equation*}
where
\begin{description}
\item[(1)] The function $H_{5,2}$ is a (vector-valued) real-analytic modular form of weight 1/2, satisfying
\begin{align*}
H_{5,2}(\tau+1) &= \begin{pmatrix} \zeta_{60}^{-1}&0&0&0&0&0\\
0&\zeta_{60}^{11}&0&0&0&0\\
0&0&0&0&\zeta_{240}^{-1}&0\\
0&0&0&0&0&\zeta_{240}^{71}\\
0&0&\zeta_{240}^{-1}&0&0&0\\
0&0&0&\zeta_{240}^{71}&0&0
\end{pmatrix} H_{5,2}(\tau),\\
\intertext{and}
H_{5,2}\Bigl(-\frac{1}{\tau}\Bigr) &=\sqrt{-i\tau} \frac{2}{\sqrt{5}}\ M_5\ H_{5,2}(\tau),
\end{align*}
with $M_5$ as in Proposition \ref{51}, and is an eigenfunction of the Casimir operator $\Omega_{\frac{1}{2}}$, with eigenvalue $\frac{3}{16}$.
\item[(2)] The function $G_{5,2}$ is bounded if $\tau\downarrow \xi$, with $\xi\in\Q$.
\end{description}
\end{proposition}
\proof
We consider the function
\begin{equation*}
\Theta (\tau) = \begin{pmatrix} \theta_{\left(\begin{smallmatrix} 1/10\\1/6\end{smallmatrix}\right), \left(\begin{smallmatrix} 0\\1/6\end{smallmatrix}\right)}\\
\theta_{\left(\begin{smallmatrix} 3/10\\1/6\end{smallmatrix}\right), \left(\begin{smallmatrix} 0\\1/6\end{smallmatrix}\right)}\\
\theta_{\left(\begin{smallmatrix} 1/5\\1/6\end{smallmatrix}\right), \left(\begin{smallmatrix} 1/10\\1/6\end{smallmatrix}\right)}\\
\theta_{\left(\begin{smallmatrix} 2/5\\1/6\end{smallmatrix}\right), \left(\begin{smallmatrix} 1/10\\1/6\end{smallmatrix}\right)}\\
\theta_{\left(\begin{smallmatrix} 1/5\\1/6\end{smallmatrix}\right), \left(\begin{smallmatrix} 0\\1/6\end{smallmatrix}\right)}\\
\theta_{\left(\begin{smallmatrix} 2/5\\1/6\end{smallmatrix}\right), \left(\begin{smallmatrix} 0\\1/6\end{smallmatrix}\right)}
\end{pmatrix}(\tau).
\end{equation*}
Using (4) and (2) of Corollary \ref{cor1} we see
\begin{equation*}
\Theta (\tau+1) = \begin{pmatrix} \zeta_{60}^{-1}&0&0&0&0&0\\
0&\zeta_{60}^{11}&0&0&0&0\\
0&0&0&0&\zeta_{120}^{19}&0\\
0&0&0&0&0&\zeta_{120}^{67}\\
0&0&\zeta_{24}^{-1}&0&0&0\\
0&0&0&\zeta_{120}^{19}&0&0
\end{pmatrix} \Theta (\tau).
\end{equation*}
Let $C=\left(\begin{smallmatrix} -1&0\\0&1\end{smallmatrix}\right)$, then $C\in O_A^+(\Z)$, $Cc_1=c_2$ and $Cc_2=c_1$. Hence, using Corollary \ref{cor2}, we find
\begin{equation*}
\theta_{a,b}^{c_1,c_2} = \theta_{Ca,Cb}^{c_2,c_1}=-\theta_{Ca,Cb}^{c_1,c_2}
\end{equation*}
Using this and Corollary \ref{cor1}, we obtain
\begin{align*}
\Theta \Bigl(-\frac{1}{\tau}\Bigr) &= -i\tau \frac{2}{\sqrt{5}} \cdot \\
&\quad \begin{pmatrix}
0&0&\zeta_{10}^{-1}\sin\frac{\pi}{5}&\zeta_5^{-1}\sin\frac{2\pi}{5}&0&0\\
0&0&\zeta_{10}^{-1}\sin\frac{2\pi}{5}&-\zeta_5^{-1}\sin\frac{\pi}{5}&0&0\\
\zeta_{10}\sin\frac{\pi}{5}&\zeta_{10}\sin\frac{2\pi}{5}&0&0&0&0\\
\zeta_5 \sin\frac{2\pi}{5}&-\zeta_5 \sin\frac{\pi}{5}&0&0&0&0\\
0&0&0&0&\sin\frac{2\pi}{5}&\sin\frac{\pi}{5}\\
0&0&0&0&\sin\frac{\pi}{5}&-\sin\frac{2\pi}{5}
\end{pmatrix} \Theta (\tau).
\end{align*}
If we use that
\begin{equation*}
f(\tau):= \begin{pmatrix} \eta(\tau/2)\\\eta\bigl((\tau+1)/2\bigr)\\\eta(2\tau)\end{pmatrix}
\end{equation*}
transforms as
\begin{equation*}
\begin{split}
f(\tau+1)&=\begin{pmatrix} 0&1&0\\ \zeta_{24}&0&0\\0&0&\zeta_{12}\end{pmatrix} f(\tau)\\
f\Bigl(-\frac{1}{\tau}\Bigr)&= \sqrt{-i\tau} \begin{pmatrix} 0&0&\sqrt{2}\\0&1&0\\\frac{1}{\sqrt{2}}&0&0\end{pmatrix} f(\tau),
\end{split}
\end{equation*}
we obtain the transformations for $H_{5,2}$.

If we write $\rho(\nu;\tau)$ as the sum of the three expressions \eqref{al1}, \eqref{al2} and \eqref{al3} and use Proposition \ref{split}, we find
\begin{equation*}
\theta_{\left(\begin{smallmatrix} 1/10\\1/6\end{smallmatrix}\right), \left(\begin{smallmatrix} 0\\1/6\end{smallmatrix}\right)}(\tau)=2\zeta_{12}^{-1} \frac{\eta(\tau)^2}{\eta(2\tau)}\ q^{-\frac{1}{60}} \psi_0(q) +\frac{\eta(\tau)^2}{\eta(2\tau)} \left( R_{\frac{1}{30},\frac{5}{2}}(30\tau) + \zeta_6^{-1} R_{\frac{11}{30},\frac{5}{2}}(30\tau)\right).
\end{equation*}
Hence
\begin{equation*}
\zeta_{12}\frac{\eta(2\tau)}{\eta(\tau)^2}\ \theta_{\left(\begin{smallmatrix} 1/10\\1/6\end{smallmatrix}\right), \left(\begin{smallmatrix} 0\\1/6\end{smallmatrix}\right)}(\tau)= 2q^{-\frac{1}{60}} \psi_0(q) + \zeta_{12} R_{\frac{1}{30},\frac{5}{2}}(30\tau) + \zeta_{12}^{-1} R_{\frac{11}{30},\frac{5}{2}}(30\tau).
\end{equation*}
We can find similar identities for the other components of $H_{5,2}$. Combining them gives
\begin{equation*}
H_{5,2} = F_{5,2} -G_{5,2}.
\end{equation*}
Using (5) of Proposition \ref{propRab} and the fact that $F_{5,2}$ is holomorphic, we find
\begin{align*}
\Omega_{\frac{1}{2}} G_{5,2}&=\frac{3}{16} G_{5,2}\\
\Omega_{\frac{1}{2}} F_{5,2}&=\frac{3}{16} F_{5,2},\\
\intertext{and so}
\Omega_{\frac{1}{2}} H_{5,2}&=\frac{3}{16} H_{5,2}.
\end{align*}
Hence $H_{5,2}$ is a vector valued real-analytic modular form of weight 1/2. From (3) of Proposition \ref{propRab} we get that $G_{5,2}$ is bounded as $\tau\downarrow \xi$, with $\xi\in\Q$. \qed

\begin{proposition}\label{propm5}
The (holomorphic) function $F_{5}:=F_{5,1}+F_{5,2}$ is a (vector-valued) modular form of weight 1/2, satisfying
\begin{align*}
F_{5}(\tau+1) &= \begin{pmatrix} \zeta_{60}^{-1}&0&0&0&0&0\\
0&\zeta_{60}^{11}&0&0&0&0\\
0&0&0&0&\zeta_{240}^{-1}&0\\
0&0&0&0&0&\zeta_{240}^{71}\\
0&0&\zeta_{240}^{-1}&0&0&0\\
0&0&0&\zeta_{240}^{71}&0&0
\end{pmatrix} F_{5}(\tau),\\
\intertext{and}
F_{5}\Bigl(-\frac{1}{\tau}\Bigr) &=\sqrt{-i\tau} \frac{2}{\sqrt{5}}\ M_5\ F_{5}(\tau),
\end{align*}
with $M_5$ as in Proposition \ref{51}.
\end{proposition}
\proof
Since $G_{5,2}=-G_{5,1}$, we get from Proposition \ref{51} and Proposition \ref{52}
\begin{equation*}
F_5=H_{5,1}+H_{5,2}.
\end{equation*}
Using this, the transformation behaviour of $F$ follows directly from the transformation behaviour of $H_{5,1}$ and $H_{5,2}$ given in Proposition \ref{51} and Proposition \ref{52}. \qed

Proposition \ref{propm5} implies that each of the six components of $F_{5}:=F_{5,1}+F_{5,2}$, i.e.\ each of the six functions
\begin{equation*}
\begin{split}
q^{-\frac{1}{60}} f_0(q)&+2q^{-\frac{1}{60}} \psi_0 (q)\\
q^{\frac{11}{60}} f_1(q)&+2q^{\frac{11}{60}} \psi_1 (q)\\
q^{-\frac{1}{240}} (-1+F_0 (q^{\frac{1}{2}}))&+q^{-\frac{1}{240}} \phi_0(-q^{\frac{1}{2}})\\
q^{\frac{71}{240}} F_1(q^{\frac{1}{2}})&-q^{-\frac{49}{240}} \phi_1(-q^{\frac{1}{2}})\\
q^{-\frac{1}{240}} (-1+F_0 (-q^{\frac{1}{2}}))&+q^{-\frac{1}{240}} \phi_0(q^{\frac{1}{2}})\\
q^{\frac{71}{240}} F_1(-q^{\frac{1}{2}}))&+q^{-\frac{49}{240}} \phi_1(q^{\frac{1}{2}})
\end{split}
\end{equation*}
is a modular form on a suitable congruence subgroup of $\operatorname{SL}_2(\Z)$. The first four functions were already given in terms of theta functions in \cite[pp.\ 299]{watson2}. Similar identities for the fifth and sixth function follow directly from the identities for the third and fourth function by replacing $q^{\frac{1}{2}}$ by $-q^{\frac{1}{2}}$.

\section{Other mock $\theta$-functions}
In the previous sections we have dealt with the seventh order and most of the fifth order mock $\theta$-functions from Ramanujan's letter. We were able to do this because identities for these mock $\theta$-functions were available in the literature. Similar identities for the other two fifth order mock $\theta$-functions, $\chi_0$ and $\chi_1$, are not available in the literature, as far as I know. I have found the following identities (which I will not prove here) for $\chi_0$ and $\chi_1$
\begin{equation*}
\begin{split}
\chi_0 (q) &= 2- \frac{1}{(q)_\infty^2} \Bigl( \sum_{k,l,m\geq 0} +\sum_{k,l,m<0}\Bigr) (-1)^{k+l+m} q^{\frac{1}{2}k^2+\frac{1}{2}l^2+\frac{1}{2}m^2+2kl+2km+2lm+\frac{1}{2}(k+l+m)},\\
\chi_1 (q) &=\frac{1}{(q)_\infty^2} \Bigl( \sum_{k,l,m\geq 0} +\sum_{k,l,m<0}\Bigr) (-1)^{k+l+m} q^{\frac{1}{2}k^2+\frac{1}{2}l^2+\frac{1}{2}m^2+2kl+2km+2lm+\frac{3}{2}(k+l+m)}.
\end{split}
\end{equation*}
These series are similar to the ones we used for the seventh and fifth order mock theta functions. However, the quadratic form is of type $(1,2)$. Hence we cannot apply any of the results from Chapter 2. However, in \cite{watson2}  identities for $\chi_0$ and $\chi_1$ are given which give $\chi_0$ and $\chi_1$ as a linear combination of other fifth order mock $\theta$-functions. Hence, we could derive the transformation properties of $\chi_0$ and $\chi_1$.

In Ramanujan's letter four third order mock $\theta$-functions are given. Watson (see \cite{watson1}) defined three more third order mock $\theta$-functions. Later even more exotic mock $\theta$-functions were introduced: of sixth order (see \cite{andrews2}), of eighth order (see \cite{mac}) and of tenth order (see \cite{choi1}, \cite{choi2} and \cite{choi3}). In these articles, identities are given, which relate the mock $\theta$-functions to sums of the same type as the ones we used for the fifth and seventh order mock $\theta$-functions. Hence, using the same techniques, we could derive the transformation properties of these mock $\theta$-functions. In \cite{ikke} I derive the transformation properties of the vector-valued third order mock $\theta$-function
\begin{equation*}
F(\tau) = \begin{pmatrix} q^{-\frac{1}{24}} f(q)\\ 2q^{\frac{1}{3}} \omega (q^{\frac{1}{2}})\\2q^{\frac{1}{3}} \omega (-q^{\frac{1}{2}})\end{pmatrix}.
\end{equation*}
The result is similar to the results found in this chapter.

\clearemptydoublepage
\addcontentsline{toc}{chapter}{\numberline{}Bibliography}
\bibliography{thesis}

\providecommand{\bysame}{\leavevmode\hbox to3em{\hrulefill}\thinspace}
\providecommand{\MR}{\relax\ifhmode\unskip\space\fi MR }
\providecommand{\MRhref}[2]{%
  \href{http://www.ams.org/mathscinet-getitem?mr=#1}{#2}
}
\providecommand{\href}[2]{#2}
\begin{thebibliography}{10}

\bibitem{andrews}
G.E. Andrews, \emph{Hecke modular forms and the {Kac-Peterson} identities},
  Trans. Amer. Math. Soc. \textbf{283} (1984), no.~2, 451--458.

\bibitem{andrews3}
\bysame, \emph{The fifth and seventh order mock theta functions}, Trans. Amer.
  Math. Soc. \textbf{293} (1986), no.~1, 113--134.

\bibitem{coeff}
\bysame, \emph{Ramanujan's fifth order mock theta functions as constant terms},
  Ramanujan Revisited: Proceedings of the Centenary Conference, Univ.\ of
  Illinois at Urbana-Champaign, June 1--5, 1987, San Diego: Acad.\ Press, 1988,
  pp.~47--56.

\bibitem{andrews1}
\bysame, \emph{Mock theta functions}, Theta functions --- Bowdoin 1987, Part 2
  (Brunswick, 1987), Proc. Symp. Pure Math., vol.~49, Amer. Math. Soc., 1989,
  pp.~283--298.

\bibitem{andrews2}
G.E. Andrews and D.~Hickerson, \emph{Ramanujan's ``lost'' notebook: the sixth
  order mock theta functions}, Adv. Math. \textbf{89} (1991), 60--105.

\bibitem{choi1}
Youn-Seo Choi, \emph{Tenth order mock theta functions in {R}amanujan's {L}ost
  {N}otebook}, Invent. Math. \textbf{136} (1999), 497--569.

\bibitem{choi2}
\bysame, \emph{Tenth order mock theta functions in {R}amanujan's {L}ost
  {N}otebook {II}}, Adv. Math. \textbf{156} (2000), 180--285.

\bibitem{choi3}
\bysame, \emph{Tenth order mock theta functions in {R}amanujan's {L}ost
  {N}otebook {IV}}, Trans. Amer. Math. Soc. \textbf{354} (2002), 705--733.

\bibitem{eichler}
M.~Eichler and D.B. Zagier, \emph{{The Theory of Jacobi Forms}}, Progress in
  Mathematics, no.~55, Birkha{\"u}ser, 1985.

\bibitem{mac}
B.~Gordon and R.J. McIntosh, \emph{Some eighth order mock theta functions}, J.
  London Math. Soc. (2) \textbf{62} (2000), 321--335.

\bibitem{gottsche}
L.~{G\"ottsche} and D.~Zagier, \emph{Jacobi forms and the structure of
  {Donaldson} invariants for 4-manifolds with $b_+=1$}, Selecta Math. (N.S.)
  \textbf{4} (1998), no.~1, 69--115.

\bibitem{hick2}
D.R. Hickerson, \emph{On the seventh order mock theta functions}, Invent. Math.
  \textbf{94} (1988), 661--677.

\bibitem{hick}
\bysame, \emph{A proof of the mock theta conjectures}, Invent. Math.
  \textbf{94} (1988), 639--660.

\bibitem{lerch2}
M.~Lerch, \emph{Bemerkungen zur {Theorie} der elliptischen {Functionen}},
  Jahrbuch {\"uber} die {Fortschritte} der {Mathematik} \textbf{24} (1892),
  442--445.

\bibitem{lerch1}
\bysame, \emph{Pozn{\'{a}}mky k theorii funkc{\'{\i}} elliptick{\'{y}}ch},
  Rozpravy {\v{C}}esk{\'{e}} Akademie C{\'{\i}}sa{\v{r}}e Franti{\v{s}}ka
  Josefa pro v{\v{e}}dy, slovesnost a um{\v{e}}n{\'{\i}} v praze (II Cl) I
  \textbf{24} (1892), 465--480.

\bibitem{maass}
H.~Maass, \emph{{Lectures on Modular Function of one Complex Variable}}, Tata
  Institute of Fundamental Research, Bombay, 1964.

\bibitem{mordell1}
L.J. Mordell, \emph{The value of the definite integral
  $\displaystyle{\int_{-\infty}^{\infty} \frac{e^{at^2+bt}}{e^{ct}+d}\,dt}$},
  Quarterly Journal of Math. \textbf{68} (1920), 329--342.

\bibitem{mordell2}
\bysame, \emph{The definite integral $\displaystyle{\int_{-\infty}^{\infty}
  \frac{e^{at^2+bt}}{e^{ct}+d}\,dt}$ and the analytic theory of numbers}, Acta
  Math. \textbf{61} (1933), 323--360.

\bibitem{mumford}
D.~Mumford, \emph{{Tata Lectures on Theta I}}, Progress in Mathematics, no.~28,
  Birkha{\"u}ser, 1983.

\bibitem{polishchuk}
A.~Polishchuk, \emph{A new look at {Hecke's} indefinite theta series},
  $q$-Series with Applications to Combinatorics, Number Theory, and Physics,
  Univ.\ of Illinois at Urbana-Champaign, October 26--28, 2000, Contemporary
  Mathematics, vol. 291, Amer. Math. Soc., 2001, pp.~183--191.

\bibitem{lost}
S.~Ramanujan, \emph{The lost notebook and other unpublished papers}, Narosa
  Publishing House, New Delhi, 1987.

\bibitem{rogers}
L.J. Rogers, \emph{On two theorems of combinatory analysis and some allied
  identities}, Proc. London Math. Soc. (2) \textbf{16} (1917), 316--336.

\bibitem{selberg}
A.~Selberg, \emph{{\"U}ber die {M}ock-{T}hetafunktionen siebenter {O}rdnung},
  Arch. Math. og Naturvidenskab \textbf{41} (1938), 3--15.

\bibitem{siegel2}
C.L. Siegel, \emph{{\"U}ber {R}iemanns {N}achla{\ss} zur analytischen
  {Z}ahlentheorie}, Quellen und Studien zur Geschichte der Mathematik,
  Astronomie und Physik \textbf{2} (1933), 45--80.

\bibitem{siegel}
\bysame, \emph{Indefinite quadratische {Formen} und {Funktionentheorie}. {I}.},
  Math. Ann. \textbf{124} (1951), 17--54.

\bibitem{watson1}
G.N. Watson, \emph{The final problem: an account of the mock theta functions},
  J. London Math. Soc. \textbf{11} (1936), 55--80.

\bibitem{watson2}
\bysame, \emph{The mock theta functions (2)}, Proc. London Math. Soc. (2)
  \textbf{42} (1937), 274--304.

\bibitem{whit}
E.T. Whittaker and G.N. Watson, \emph{{Modern Analysis}}, fourth ed., Cambridge
  at the University Press, 1927.

\bibitem{zagier}
D.B. Zagier, \emph{Nombres de classes et formes modulaires de poids $3/2$}, C.
  R. Acad. Sci. Paris S\'er. A-B \textbf{281} (1975), no.~21, Ai, A883--A886.

\bibitem{ikke}
S.P. Zwegers, \emph{Mock $\theta$-functions and real analytic modular forms},
  $q$-Series with Applications to Combinatorics, Number Theory, and Physics,
  Univ.\ of Illinois at Urbana-Champaign, October 26--28, 2000, Contemporary
  Mathematics, vol. 291, Amer. Math. Soc., 2001, pp.~269--277.

\end{thebibliography}
\clearemptydoublepage
\addcontentsline{toc}{chapter}{\numberline{}Samenvatting}
\chapter*{Samenvatting}
\markboth{Samenvatting}{}

Dit proefschrift gaat over mock thetafuncties. Deze mock thetafuncties zijn een ``uitvinding'' van de Indiase wiskundige Srinivasa Ramanujan. Ramanujan leefde van 1887 tot 1920. Hoewel hij nooit een universitaire studie heeft afgemaakt, wordt hij door velen gezien als een wiskundig genie. Ramanujan werkte erg intu\"{\i}tief en bewees zelden z'n beweringen. Dit komt doordat hij door z'n onvolledige wiskundige opleiding nooit de kunst van het bewijzen heeft aangeleerd. Het leven en werk van Ramanujan bevat dan ook een element van mysterie en romantiek.

De mock thetafuncties ontdekte Ramanujan kort voor hij op 32 jarige leeftijd stierf. Hoewel verschillende wiskundigen zich na de dood van Ramanujan met de mock thetafuncties hebben beziggehouden, en ook verschillende resultaten hebben geboekt, is nooit echt duidelijk geworden wat er nu werkelijk aan de hand is. Daardoor vormen ze nog steeds een bron van raadsels voor hedendaagse wiskundigen. In dit proefschrift plaats ik de voorbeelden die Ramanujan gaf van mock thetafuncties in een diepere achterliggende theorie (namelijk die van re\"eel analytische modulaire vormen). Hiermee kan een natuurlijke verklaring worden gegeven voor de eigenschappen van mock thetafuncties, zoals Ramanujan die beschreef. 

Het engelse werkwoord ``to mock'' betekent overigens zoiets als bespotten of spottend na\"apen. We zouden dus ook kunnen spreken over spottende thetafuncties of nepthetafuncties (thetafuncties waren in Ramanujans tijd reeds uitgebreid bestudeerd en ook Ramanujan was hier vertrouwd mee), maar dit dekt niet helemaal de lading. In ``Alice in Wonderland'' van Lewis Carroll komt een nepschildpad (mock turtle) voor. Hiervan wordt nepschildpadsoep (mock turtle soup) gemaakt. In werkelijkheid is nepschildpadsoep natuurlijk nagemaakte schildpadsoep (gemaakt van kalfshoofd, kalfsvlees, etc.)

\bigskip

In hoofdstuk 1 bekijk ik de som
\begin{equation*}
\sum_{n\in\Z} \frac{(-1)^n e^{\pi i(n^2+n)\tau+2\pi inv}}{1-e^{2\pi in\tau + 2\pi iu}}\qquad (\tau\in\Ha, v\in\C, u\in\C\setminus (\Z\tau+\Z)).
\end{equation*}
Ik noem dit een Lerch som, omdat deze som ook al werd bestudeerd door Lerch.
Deze Lerch som transformeert bijna als een Jacobivorm onder substituties in $(u,v,\tau)$. Ik laat zien dat het transformatiegedrag precies dat van een Jacobivorm wordt als we er een (relatief eenvoudige) correctieterm bij optellen. Deze correctieterm blijkt niet holomorf te zijn in $(u,v,\tau)$, alleen re\"eel analytisch. Voor bepaalde waarden van $(u,v)$ zouden we de Lerch som als functie van $\tau$ een mock thetafunctie kunnen noemen, hoewel deze niet expliciet bij Ramanujan voorkomen. 

In hoofdstuk 2 bekijk ik thetafuncties bij indefiniete kwadratische vormen. Deze thetafuncties zijn een aanpassing van een klasse van thetafuncties ge\"{\i}ntroduceerd door G\"ottsche en Zagier, en lijken op de thetafuncties ingevoerd door Siegel. Voor deze indefiniete thetafuncties vind ik elliptisch en modulair transformatiegedrag, analoog aan het transformatiegedrag van thetafuncties behorende bij positief definiete kwadratische vormen. In het geval van positief definiete kwadratische vormen zijn de thetafuncties holomorf. De thetafuncties in dit hoofdstuk zijn dat niet. Door specialisatie van de parameters is het mogelijk om mock thetafuncties te verkrijgen uit deze indefiniete thetafuncties.

In hoofdstuk 3 bestudeer ik de modulariteit van de Fourier co\"effici\"enten van meromorfe Jacobivormen. 

In hoofdstuk 4 gebruik ik het verband tussen mock thetafuncties en de indefiniete thetafuncties uit hoofdstuk 2, om deze mock thetafuncties in verband te brengen met re\"eel analytische modulaire vormen. We hadden echter even goed de methoden uit hoofdstuk 1 of hoofdstuk 3 kunnen gebruiken om tot hetzelfde resultaat te komen. Niet alle mock thetafuncties vallen in dit kader: Voor een tweetal mock thetafuncties van orde 5 zou een uitbreiding van mijn theorie gewenst zijn tot een andere klasse van indefiniete kwadratische vormen.

\addcontentsline{toc}{chapter}{\numberline{}Dankwoord}
\chapter*{Dankwoord}

Op deze plaats wil ik een aantal mensen noemen die een rol hebben gespeeld bij het tot stand komen van dit proefschrift.

Allereerst ben ik veel dank verschuldigd aan Don Zagier en Roelof Bruggeman, mijn promotoren en begeleiders. Don: hoewel onze gesprekken niet erg veelvuldig waren, waren ze altijd erg stimulerend en inspirerend. De belangrijkste idee\"en in dit proefschrift zijn dan ook uit deze gesprekken voortgekomen. Verder ben ik je veel dank verschuldigd voor het helpen bij het (enigszins) leesbaar maken van m'n proefschrift. Roelof: het was erg fijn om iemand te hebben om wekelijks de zaken door te spreken. Zeker in de eindfase, omdat je iedere keer weer bereid was om commentaar te leveren op hetgeen ik geschreven had. Hoewel ik je in het begin misschien iets te precies vond, ben ik je daar achteraf alleen maar dankbaar voor. 

I would like to thank my thesis examiners Bruce Berndt, Hans Duistermaat, Joop Kolk and Lo\"{\i}c Merel, for reading the manuscript. 

Verder ben ik m'n ouders dankbaar voor hun voortdurende belangstelling en hun vertrouwen. 

Tenslotte bedank ik Therese, voor alles.

\addcontentsline{toc}{chapter}{\numberline{}Curriculum Vitae}
\chapter*{Curriculum Vitae}
\thispagestyle{empty}

Sander Zwegers is op 16 april 1975 geboren te Oosterhout (N.B.). Hij groeide op in 's Gravenmoer en haalde in 1993 het VWO-diploma aan het Dr. Schaepmancollege te Dongen. Hierna ging hij studeren aan de Universiteit Utrecht, en in augustus 1994 haalde hij zijn propedeuse Wiskunde (cum laude) en Natuurkunde (cum laude). In augustus 1998 studeerde hij af in de Wiskunde (cum laude) bij Dr.\ J.A.C. Kolk, met de scriptie ``Theta functions, certain modular forms arising from them and Rogers--Ramanujan type identities''. Tijdens zijn studie was Sander lid van de faculteitsraad van de Faculteit Wiskunde en Informatica, gaf hij als studentassistent werkcolleges en vervulde hij diverse functies bij de studievereniging A-Eskwadraat, waaronder het voorzitterschap. 

Vanaf oktober 1998 was hij werkzaam als assistent in opleiding (AiO) bij het Mathematisch Instituut van de faculteit Wiskunde en Informatica van de Universiteit Utrecht, met Prof.\ dr.\ D.B. Zagier als promotor, en Dr.\ R.W. Bruggeman als copromotor. Het onderzoek dat hij daar deed heeft geresulteerd in dit proefschrift. In het kader van zijn onderzoek heeft hij voordrachten geven op enkele conferenties en instituten: in Urbana-Champaign (VS), Bonn (Duitsland) en Luminy (Frankrijk). Daarnaast gaf hij onderwijs, was hij lid van de faculteitsraad van de Faculteit Wiskunde en Informatica en was hij lid van de benoemingsadviescommissie algebra.

Vanaf 1 januari 2003 zal hij voor \'e\'en jaar als postdoc aan het Max-Planck-Institut f\"ur Mathematik (Bonn) verbonden zijn.

\end{document}